\documentclass[aos]{imsart}
\makeatletter\def\journal@name{}\makeatother

\RequirePackage{amsthm,amsmath,amsfonts,amssymb}
\RequirePackage[numbers]{natbib}
\RequirePackage[colorlinks,citecolor=blue,urlcolor=blue]{hyperref}
\RequirePackage{graphicx}

\startlocaldefs
\usepackage{upgreek}
\usepackage[T1]{fontenc}
\usepackage{mathtools}
\usepackage{enumitem,microtype}
\usepackage{capt-of,booktabs}
\usepackage{tikz}
\usetikzlibrary{arrows.meta}
\usepackage{fix-cm}
\usepackage{xurl}
\DeclareFontShape{T1}{ptm}{m}{scit}{<->ssub*ptm/m/n}{}

\hypersetup{
  pdftitle={Breaking Resonance: Square-root acceleration for Metropolis-adjusted randomized Hamiltonian Monte Carlo},
  pdfauthor = {Nawaf Bou-Rabee}
}
\numberwithin{equation}{section}
\newtheorem{theorem}{Theorem}
\newtheorem{corollary}{Corollary}
\newtheorem{lemma}{Lemma}[section]
\newtheorem{proposition}[lemma]{Proposition}
\newtheorem{mainproposition}{Proposition}

\theoremstyle{definition}
\newtheorem{definition}[lemma]{Definition}
\newtheorem{remark}[lemma]{Remark}
\newtheorem{algorithm}{Algorithm}[section]
\newcommand{\R}{\mathbb R}
\newcommand{\E}{\mathbb E}
\newcommand{\Pbb}{\mathbb P}

\newcommand{\dd}{\,\mathrm d}

\DeclareMathOperator{\KL}{KL}
\DeclareMathOperator{\TV}{TV}

\DeclareMathOperator{\Law}{Law}

\allowdisplaybreaks[2]
\setlist{itemsep=3pt,topsep=5pt}
\endlocaldefs

\begin{document}

\begin{frontmatter}
\title{Breaking Resonance: Square-root acceleration for Metropolis-adjusted randomized Hamiltonian Monte Carlo}
\runtitle{Breaking resonance in Metropolis-adjusted randomized HMC}

\begin{aug}
\author[A,B]{\fnms{Nawaf}~\snm{Bou-Rabee}\ead[label=e1]{nawaf.bourabee@rutgers.edu}}

\address[A]{Department of Mathematical Sciences, Rutgers University\printead[presep={,\ }]{e1}}
\address[B]{Center for Computational Mathematics, Flatiron Institute}
\end{aug}

\begin{abstract}
Hamiltonian Monte Carlo with a fixed integration time can fail to mix
because of resonances, and randomizing the integration time removes
them. We prove that Metropolis-adjusted HMC with velocity Verlet and
random integration times mixes from a warm start in a number of
transitions that does not depend on the condition number, for every
integration-time law satisfying a Fourier condition that we call
non-resonant; uniform, exponential and triangular laws are examples.
For a target density proportional to $e^{-U}$ on $\mathbb R^d$, where
$U$ is $m$-strongly convex with $L$-Lipschitz gradient and
$L_H$-Lipschitz Hessian and $\kappa=L/m$, the number of gradient
evaluations needed to reach a fixed total variation accuracy is, up to
constants and logarithmic factors,
$\sqrt{\kappa d}+\sqrt{L_H}\,(d/m)^{3/4}$. This is the first square-root
condition-number bound for this sampler beyond Gaussian targets. The
proof converts relative-entropy contraction of the unadjusted proposal
into a conductance bound for the adjusted chain, and shows that the
integration-time law enters only through its characteristic function.
Hessian regularity cannot be dropped for a fixed step size: without it
there are targets in the class for which the cost is at least of order
$\kappa^{2/3}$, up to logarithmic factors, even from a warm start.
\end{abstract}

\begin{keyword}[class=MSC]
\kwdgroup[type=primary]{\kwd{65C05}
\kwd{60J05}}
\kwdgroup[type=secondary]{\kwd{65C40}}
\end{keyword}

\begin{keyword}
\kwd{Hamiltonian Monte Carlo}
\kwd{randomized integration time}
\kwd{Metropolis adjustment}
\kwd{mixing time}
\kwd{condition number}
\kwd{restricted conductance}
\end{keyword}

\end{frontmatter}

\section{Introduction}\label{sec:intro}

Let $\pi$ be a probability measure on $\R^d$ with density proportional
to $e^{-U}$, where $U$ is $m$-strongly convex with $L$-Lipschitz
gradient and $L_H$-Lipschitz Hessian, and let $\kappa=L/m$ be the
condition number. We prove that Metropolis-adjusted Hamiltonian Monte
Carlo (HMC), using velocity Verlet and randomized integration times
obtained from a non-resonant duration law
(Definition~\ref{def:nonresonant}) of mean order $m^{-1/2}$,
reaches total variation error $\varepsilon\in(0,1]$ from an $L^2$-warm
start after
\begin{equation}\label{intro:cost}
 \mathcal O\Bigl(\bigl[\sqrt{\kappa D}+\sqrt{L_H}\,(D/m)^{3/4}\bigr]
  \log\tfrac{2M}{\varepsilon}\Bigr),
   \qquad D=1+d+\log(M/\varepsilon),
\end{equation}
gradient evaluations in expectation, where
$M\ge1$ bounds $\|\mathrm d\nu_0/\mathrm d\pi\|_{L^2(\pi)}$ for the
initial law $\nu_0\ll\pi$ (Theorem~\ref{up:main} and
Corollary~\ref{up:gradient-complexity}). Under the duration-law
conditions of Theorem~\ref{lower:thm}, we also show that strong
convexity and a Lipschitz gradient alone do not suffice for a
gradient cost of order $\sqrt\kappa$, even from a warm start:
for every step size prescribed before the target is chosen, including
the step sizes used in Theorem~\ref{up:main},
there is a target on which the cost is at least of order
$\kappa^{2/3}$ up to logarithmic factors.

For Metropolized HMC with one Verlet step per
transition, a known mixing bound grows linearly in $\kappa$
\cite{LSTupper}. Longer Hamiltonian trajectories can improve this
dependence, but a fixed integration time can produce resonances. On a
Gaussian target, consider a mode whose half-period divides the
integration time: the flow maps the position in that coordinate to
plus or minus its initial value, regardless of the momentum drawn, so
the chain never mixes in that coordinate. Suitable randomization of
the integration time avoids this obstruction.
Mitra, Srinivasan, Wang and Wibisono
\cite[Theorem~1 and Corollary~1]{MSSW} prove that exact-flow
randomized HMC with triangular or exponential durations of mean order
$m^{-1/2}$ contracts relative entropy geometrically per transition,
giving an expected total integration time of order $m^{-1/2}$ up to
logarithmic factors. 

With momentum $v\in\R^d$, the Hamiltonian is
$H(x,v)=U(x)+|v|^2/2$, and $\Pi=\pi\otimes N(0,I_d)$ is the Boltzmann-Gibbs law on $\R^{2d}$.
Each transition draws a Gaussian momentum and an independent duration
$\tau$ from a fixed law $\rho$, then takes
$J=\lfloor\tau/h\rfloor+B$ Verlet steps of step size $h$, where $B$ is an
independent fair Bernoulli variable. The endpoint is accepted or
rejected by a Metropolis test based on the change in $H$; on
rejection, the chain stays at its current position. Write $P$ for the position kernel. The duration law $\rho$ is non-resonant in the sense of
Definition~\ref{def:nonresonant}: for constants $b_\rho>0$ and
$0<\delta_\rho\le1$,
$\mathbb E\tau\le b_\rho/\sqrt m$ and
$\operatorname{Re}\widehat\rho(\omega)\le1-\delta_\rho$ for
$\omega\ge2\sqrt m$.

\paragraph*{Main results}
Theorem~\ref{up:main} gives $\TV(\nu_0P^n,\pi)\le\varepsilon$ after
$n=\mathcal O(\log(2M/\varepsilon))$ transitions, with constants
depending only on $b_\rho,\delta_\rho$, so the condition number enters
only through the cost per transition. The second term in
\eqref{intro:cost} is at most the first precisely when
$L_H\le L\sqrt m/\sqrt D$. It comes from the cubic Verlet energy
remainder $L_Hh^3|w_n|^3/12$ of Lemma~\ref{ent:shell}, where $w_n$ is
the half-step velocity: with $|w_n|^2$ of order the energy cutoff
$R=\mathcal O(D)$, the error accumulated over a trajectory of length of
order $m^{-1/2}$ is of order $L_Hh^2D^{3/2}m^{-1/2}$, and keeping it
bounded requires $h$ of order $L_H^{-1/2}D^{-3/4}m^{1/4}$, hence
$\sqrt{L_H}(D/m)^{3/4}$ gradient evaluations per trajectory. The $\sqrt\kappa$ dependence cannot be improved
for this sampler on Gaussian targets when $h\sqrt L\le1$, even from a
warm start (Proposition~\ref{gaussian:cost-lower}).

Theorem~\ref{lower:thm} gives an obstruction to dropping Hessian
regularity. For each prescribed step size $0<h\le1$, it constructs a
strongly convex potential with $m=\kappa^{-1}$, $L=1$ and discontinuous
Hessian, and an initial law of density at most twice $\pi$, for which
the duration laws satisfying \eqref{lower:time-law} require at least
order $\kappa^{2/3}$ expected gradient evaluations, up to logarithmic
factors, to reach total variation error $1/4$. It does not assert
optimality of the exponent or of the $L_H$-dependence in the upper
bound.
   
\subsection{Contributions}\label{sec:contrib}

Theorem~\ref{up:main} is the first square-root condition-number bound
for Metropolis-adjusted HMC with randomized integration times beyond
Gaussian targets; the only previous bound of this kind is the Gaussian
analysis of \cite{AGS}.

The first contribution, Lemma~\ref{up:reference-conductance}, turns a
relative-entropy contraction for the unadjusted chain into a restricted
conductance bound for the adjusted chain. Relative-entropy
contraction is the mechanism behind the square-root bounds of
\cite{MSSW} for exact randomized HMC and of \cite{PP} for
unadjusted OBABO; a Metropolis step does not
preserve it, and \cite{PP} circumvents this by a coupling that requires
entropy estimates from point starts. Proposition~\ref{ent:raw-entropy}
establishes the contraction for the unadjusted proposal; it is a
discrete-time version of \cite[Lemma~3]{MSSW}, in which the
free-transport inequality of \cite{PP} controls the velocity entropy
between momentum kicks and Hessian regularity controls the accumulated
energy error. We apply it to a single initial law, the conditional law of $\Pi$
given $\{X\in A,H\le R\}$, whose relative entropy with respect to $\Pi$
is close to $\log(1/\pi(A))$. A proposal that remained in $A$ with high probability
would retain almost all of this entropy, so the contraction forces an
exit probability bounded below, and energy control makes these exits
likely to be accepted. The lemma is stated for a general
Metropolis-adjusted kernel whose proposal satisfies such a contraction;
with the lower spectral bound it gives mixing through the weak
Poincar\'e inequality of \cite{ALPW}.

The second contribution is the observation that the integration-time
law enters the Dirichlet form
of the adjusted kernel only through its characteristic function
(Lemma~\ref{up:law-comparison}), because the Verlet map preserves
phase-space volume. Restricted conductance is therefore proved for the
triangular reference law only and transfers to every non-resonant law
with an explicit rejection error; the same representation yields the
lower spectral bound. This is why the Fourier condition in
Definition~\ref{def:nonresonant} is the natural hypothesis. 

The third contribution is Theorem~\ref{lower:thm}, whose proof
identifies numerical instability, rather than resonance, as the
obstruction: Verlet orbits expand in the fast coordinates of the
constructed potential, and the resulting rejections block motion in its
slow coordinate.

\subsection{Outline of the proof of Theorem~\ref{up:main}}

The proof has three parts.
\emph{Restricted conductance and mixing.}
For a Borel set $A\subseteq\R^d$, the stationary flow out of $A$ is
$Q_P(A,A^c)=\int_A P(x,A^c)\,\pi(\mathrm dx)$.
Proposition~\ref{up:warm-mixing}, based on a weak Poincar\'e inequality
of Andrieu, Lee, Power and Wang \cite{ALPW}, applies to a reversible
kernel whose flow satisfies $Q_P(A,A^c)\ge\phi_P\pi(A)$ whenever
$s\le\pi(A)\le1/2$, and whose spectrum lies above $-1+\beta_P$.
For a suitable threshold $s$ depending on $M,\varepsilon$, it gives
mixing from an $L^2$-warm start in
$\mathcal O((\beta_P\phi_P^2)^{-1}\log(2M/\varepsilon))$ transitions.

\emph{Restricted conductance for a reference law.}
For decreasing triangular durations on $[0,4/\sqrt m]$ with the same
rounding, Proposition~\ref{ent:raw-entropy} bounds the relative entropy
of the position proposal by a fixed fraction of the initial joint
relative entropy, plus a discretization error, for initial laws of
bounded energy and finite entropy with $\mathbb E[V\mid X]=0$.
The proof adapts \cite[Lemma~3]{MSSW}: Lemma~\ref{ent:free-flight}
controls the second variation of the squared Wasserstein distance
along free drifts, and Lemma~\ref{ent:weighted-entropy} combines it
with the momentum kicks and convexity of $U$.
Lemma~\ref{up:reference-conductance} turns this estimate into a
restricted conductance bound for the adjusted kernel, using the energy
bounds of Lemma~\ref{ent:shell} to control acceptance.

\emph{Transfer and the lower spectral bound.}
Lemma~\ref{up:law-comparison} represents the Dirichlet forms of the
adjusted kernels using the unitary operator induced by one Verlet step,
integrated over density levels, and a rejection term, so Fourier
comparison of integration-time laws compares their Dirichlet forms
with an explicit rejection error.
We condition the duration on $\tau\le\tau_*$, with $\tau_*$ of order
$m^{-1/2}$ and retained probability bounded below, so that physical
integration times are at most $\tau_*+h$ and the energy estimates
control this error.
Lemma~\ref{up:rounded-fourier} gives comparison with the triangular
reference law uniformly over admissible step sizes; evaluation at
$\mathbf1_A$ and inclusion of the conditional kernel in the full
mixture give Proposition~\ref{up:restricted-conductance}.
A separate Fourier bound for the full duration law gives the lower
spectral bound of Proposition~\ref{up:spectral-bound}.

\subsection{Relation to \emph{Persistence Pays}}\label{sec:persistence}

The companion paper \cite{PP} obtains the two
leading terms of the gradient cost \eqref{intro:cost}, with different
logarithmic factors, for Metropolis-adjusted OBABO, a kinetic Langevin scheme with partial
momentum refreshment and a non-reversibly updated accept/reject
uniform \cite{NealUniform}. The essential difference is reversibility. The sampler in
\cite{PP} is not reversible, so its analysis proceeds by coupling: two
copies of the adjusted chain are coupled through their unadjusted
epoch endpoints, which requires entropy estimates from point starts
and control of rejections within each epoch.

The sampler studied here is reversible with respect to $\pi$. This
makes the spectral route available: restricted conductance and a lower
spectral bound give mixing through a weak Poincar\'e inequality, and
the Fourier comparison of Lemma~\ref{up:law-comparison} transfers
conductance from the decreasing triangular law to every non-resonant
law. The
entropy estimate is applied only to $\Pi$ conditioned on a set of
positions and an energy shell, and no control of the unadjusted chain
from point starts is needed. The two papers share the free-transport
entropy inequality and the exact expression for the energy error of
velocity Verlet (Lemma~\ref{ent:shell}).

Oberd\"orster \cite{Oberdorster} has since combined the entropy argument
of \cite{PP} with the energy-error analysis of Chen, Gatmiry and Jiang
\cite{CGJ} and a stability result for high-probability domains, obtaining
for adjusted OBABO with persistent uniforms a gradient complexity of order
$(\sqrt\kappa+L_H^{1/2}m^{-3/4})\,d^{1/4}$, uniformly from deterministic
starts in a high-probability set, under Hessian-Lipschitz continuity in
the Frobenius norm. The present paper keeps the spectral-norm assumption
\ref{ass:rhmc-a3} and the warm start. Its dimension dependence comes from
the energy-control lemma, and a Frobenius-norm estimate of the kind used
in \cite{Oberdorster} could be substituted there without affecting the
treatment of the integration-time law.

\subsection{Related work}\label{sec:related}

The sensitivity of HMC to its integration time has long been
recognized. Mackenzie \cite{Mackenzie} proposed duration randomization
to avoid resonances; Neal \cite[Section~5.4.2.3]{NealHMC} discusses
this remedy and the choice between short and unnecessarily long
trajectories. The No-U-Turn Sampler of Hoffman and Gelman \cite{NUTS}
chooses the trajectory length locally using a U-turn criterion. The
GIST framework of Bou-Rabee, Carpenter and Marsden \cite{GIST} places
randomized HMC and NUTS in a common auxiliary-variable construction,
giving an interpretation of NUTS as a locally adaptive form of
randomized HMC. WALNUTS \cite{WALNUTS} extends the local adaptation to
the leapfrog step size along a trajectory.

For exact-flow randomized HMC with exponential durations, Bou-Rabee
and Sanz-Serna \cite[Theorem~3.9]{BRSS} established geometric
ergodicity through a Foster-Lyapunov drift condition, minorization and
Harris's theorem. Lu and Wang \cite{LW} obtained explicit $L^2$
convergence rates by a variational hypocoercivity argument. For
$m$-strongly convex potentials with Lipschitz gradients, their rate is
of order $\sqrt m$ in physical time when the refreshment rate is of
that order. These exact-flow sampling results do not account for
numerical integration or its gradient cost.

The entropy argument we discretize is connected to optimization in
the space of probability measures, as developed in Wibisono
\cite{Wibisono18}. Fu and Wibisono \cite{FW} prove accelerated
Hamiltonian descent with randomized integration times and velocity
reset to zero. Wang, Srinivasan, Fu, Mitra, Wilson and Wibisono
\cite{WSFMW} obtain deterministic acceleration using weighted averages
along Hamiltonian trajectories. Both works include discrete
optimization algorithms with accelerated complexity guarantees.
Mitra, Srinivasan, Wang and Wibisono \cite[Appendix~B]{MSSW} carry the
averaging argument over to relative entropy, with Gaussian momentum
refreshment in place of the zero-velocity restart.

For Gaussian targets, Wang and Wibisono \cite{WW} prove acceleration
for exact HMC with a deterministic Chebyshev schedule of integration
times. Apers, Gribling and Szil\'agyi \cite{AGS} analyze
Metropolis-adjusted randomized HMC with Verlet and obtain total
variation error $\varepsilon$ at a gradient cost of order
$\sqrt\kappa\,d^{1/4}\log(1/\varepsilon)$, up to additional logarithmic
factors. Their proof uses Hanson-Wright concentration and restricted
conductance; they also construct a warm start using unadjusted HMC.
Their dimension dependence is sharper than ours.

Lee, Shen and Tian \cite{LSTupper} prove mixing of Metropolized HMC
with one leapfrog step per transition in order $\kappa d$ iterations,
up to logarithmic factors, by gradient concentration and restricted
conductance. Their lower bound \cite[Theorem~4]{LSTlower} gives, for
every fixed step size and fixed number $j$ of leapfrog steps per
transition, a Gaussian target with relaxation time at least a constant
multiple of $\kappa\sqrt d/(j\sqrt{\log d})$ when
$\kappa=\mathcal O(d^4)$. The resonant mode in that construction
depends on the fixed duration, so the bound does not apply to
randomized durations.

\subsection{Organization of the paper}\label{sec:organization}

Section~\ref{sec:setting} gives the assumptions, sampler, duration laws
and main results. Section~\ref{sec:upper} proves
Theorem~\ref{up:main} from three propositions, established in
Sections~\ref{sec:rhmc-criterion}, \ref{sec:rhmc-invariance}
and~\ref{sec:rhmc-conductance}: the mixing criterion, the lower
spectral bound, and the restricted conductance bound.
Section~\ref{sec:lower} proves Theorem~\ref{lower:thm}.
Appendix~\ref{app:regularization} contains the regularization argument
used in Section~\ref{sec:rhmc-proposal-entropy}.
Appendices~\ref{app:gaussian-lower} and~\ref{app:duration-fourier}
give the Gaussian lower bound and explicit Fourier constants for the
three duration laws, respectively.

\subsection{Formal verification}\label{sec:lean}

The statements of Theorems~\ref{up:main} and \ref{lower:thm} and
Proposition~\ref{gaussian:cost-lower}, with the definitions they rest on,
have been formalized and their proofs checked in Lean~4 against the
Mathlib library, with no axioms beyond the standard three of Lean's
kernel. The formalization follows the proofs given here, including the
explicit constants, and covers Propositions~\ref{up:warm-mixing},
\ref{up:spectral-bound}, \ref{up:restricted-conductance} and
\ref{ent:raw-entropy} and Lemmas~\ref{up:reversibility},
\ref{up:law-comparison}, \ref{up:rounded-fourier},
\ref{up:example-fourier}, \ref{ent:modified-energy},
\ref{up:reference-conductance}, \ref{lower:cone} and
\ref{lower:phase-event}. Two results from the literature enter the proof
of Theorem~\ref{up:main} as hypotheses: Talagrand's transport-entropy
inequality (Lemma~\ref{lem:talagrand}) and the interpolation estimate of
\cite{PP} (Lemma~\ref{lem:interpolation}). The Lean sources, a table
matching each formal statement to its counterpart in this paper, and the
build log are available at \url{https://github.com/nawafbourabee/breaking-resonance-lean} (archived at
\url{https://doi.org/10.5281/zenodo.24483261}).

   \section{Setting and main results}\label{sec:setting}\label{sec:prelim}

\subsection{Assumptions and notation}
Let $U:\R^d\to\R$ be the target potential, with $d\ge1$. Write $|\cdot|$
for the Euclidean norm and $\|M\|_{\mathrm{op}}=\sup_{|z|=1}|Mz|$
for its induced operator norm on matrices. We impose the following
assumptions as needed.
\begin{enumerate}[label=(A\arabic*),ref=(A\arabic*),leftmargin=*]
\item\label{ass:rhmc-a1} $U$ is $m$-strongly convex, with $m>0$.
\item\label{ass:rhmc-a2} $U$ is differentiable and $\nabla U$ is
$L$-Lipschitz, with $L>0$.
\item\label{ass:rhmc-a3} $U$ is twice continuously differentiable and
$\|\nabla^2U(x)-\nabla^2U(y)\|_{\mathrm{op}}\le L_H|x-y|$ for every
$x,y\in\R^d$, with $L_H\ge0$.
\end{enumerate}
Write $\kappa=L/m$ for the condition number.

Let $k\ge1$ and let $\nu,\mu$ be probability measures on $\R^k$. We use
$k=d$ for position laws and
$k=2d$ for joint laws of position and velocity. Their relative entropy is
\[
 \KL(\nu\Vert\mu)=
 \begin{cases}
  \displaystyle\int_{\R^k} f(z)\log f(z)\,\mu(\mathrm dz),&\nu\ll\mu,\\[2pt]
  +\infty,&\text{otherwise},
 \end{cases}
\]
where $f=\mathrm d\nu/\mathrm d\mu$ is the Radon-Nikodym derivative in the
absolutely continuous case, with $0\log0=0$.
Their total variation distance is
\[
 \TV(\nu,\mu)=\sup_{A}|\nu(A)-\mu(A)|,
\]
the supremum over Borel sets $A\subseteq\R^k$.
Their 2-Wasserstein distance is
\[
 W_2(\nu,\mu)=\left(\inf_{\gamma\in\Gamma(\nu,\mu)}
       \int_{\R^k\times\R^k}|x-y|^2\,\gamma(\mathrm dx,\mathrm dy)\right)^{1/2},
\]
where $\Gamma(\nu,\mu)$ is the set of probability measures on
$\R^k\times\R^k$ with first marginal $\nu$ and second marginal $\mu$.

Write $Z=\int_{\R^d}e^{-U(x)}\,\mathrm dx$, assumed finite, and define
\begin{equation}\label{pre:target}
 \pi(\dd x)=Z^{-1}e^{-U(x)}\dd x,\qquad
 \gamma_d=N(0,I_d),\qquad \Pi=\pi\otimes\gamma_d.
\end{equation}

\subsection{Metropolis-adjusted randomized HMC}\label{sec:rhmc-sampler}

Under \ref{ass:rhmc-a1}, denote the unique minimizer of $U$ by $x_*$.
Fix a step size $h>0$ and a duration law $\rho$.
The Hamiltonian is
\begin{equation}\label{up:energy}
 H(x,v)=U(x)-U(x_*)+\tfrac12|v|^2.
\end{equation}

\par\medskip\begingroup\setlength{\fboxsep}{6pt}
\noindent\fbox{\begin{minipage}{\dimexpr\linewidth-2\fboxsep-2\fboxrule\relax}
\begin{algorithm}[Metropolis-adjusted randomized HMC]\label{alg:rhmc-transition}
Given the current position $X$, return the next position $X^+$ as follows.
\begin{enumerate}[label=\arabic*.,leftmargin=*]
\item \textbf{Duration and momentum randomization.} Draw $V\sim\gamma_d$,
$\tau\sim\rho$, and a fair Bernoulli $B$ independently of one another
and of $X$. Set the integration time
\begin{equation}\label{up:integration-time}
 J=\lfloor\tau/h\rfloor+B.
\end{equation}
\item \textbf{Verlet integration.} With $(x_0,v_0)=(X,V)$,
apply the following update $J$ times.
\begin{equation}\label{up:verlet}
 \begin{split}
 w_k&=v_k-\tfrac h2\nabla U(x_k),\\
 x_{k+1}&=x_k+hw_k,\\
 v_{k+1}&=w_k-\tfrac h2\nabla U(x_{k+1}).
 \end{split}
\end{equation}
Write $(Q,W)=(x_J,v_J)$.
\item \textbf{Metropolis test.} Draw $Z_{\mathrm{acc}}\sim\operatorname{Unif}(0,1)$
independently and compute
\begin{equation}\label{up:acceptance}
 \alpha(X,V;Q,W)
   =1\wedge\exp\{H(X,V)-H(Q,W)\}.
\end{equation}
Return $X^+=Q$ if $Z_{\mathrm{acc}}\le\alpha(X,V;Q,W)$, and $X^+=X$ otherwise.
\end{enumerate}
\end{algorithm}
\end{minipage}}\par\endgroup\medskip

\paragraph*{Position kernel}
Let $\Phi_h$ be the Verlet map in \eqref{up:verlet}, and write
$\Phi_h^{\circ j}$ for its $j$-fold composition, with $\Phi_h^{\circ0}$
the identity on $\R^{2d}$. For $z=(x,v)$ and
$j\ge0$, let $Q_j(z)$ be the position component of $\Phi_h^{\circ j}(z)$
and set $\alpha_j(z)=1\wedge\exp\{H(z)-H(\Phi_h^{\circ j}(z))\}$.
For a measurable set $A\subseteq\R^d$, the position kernels are
\begin{align*}
 K_j(x,A)
 &=\int_{\R^d}
   \left[\alpha_j(x,v)\mathbf1_A(Q_j(x,v))+
         (1-\alpha_j(x,v))\mathbf1_A(x)\right]\,\gamma_d(\mathrm dv),\\
 P(x,A)&=\sum_{j\ge0}\mathbb P(J=j)K_j(x,A).
\end{align*}
Here $K_j$ uses $j$ Verlet steps, and $P$ is the position kernel of
Algorithm~\ref{alg:rhmc-transition}. In particular, the law after $n$
iterations from $\nu_0$ is $\nu_0P^n$.

\subsection{Non-resonant integration-time laws}
\label{sec:nonresonant}

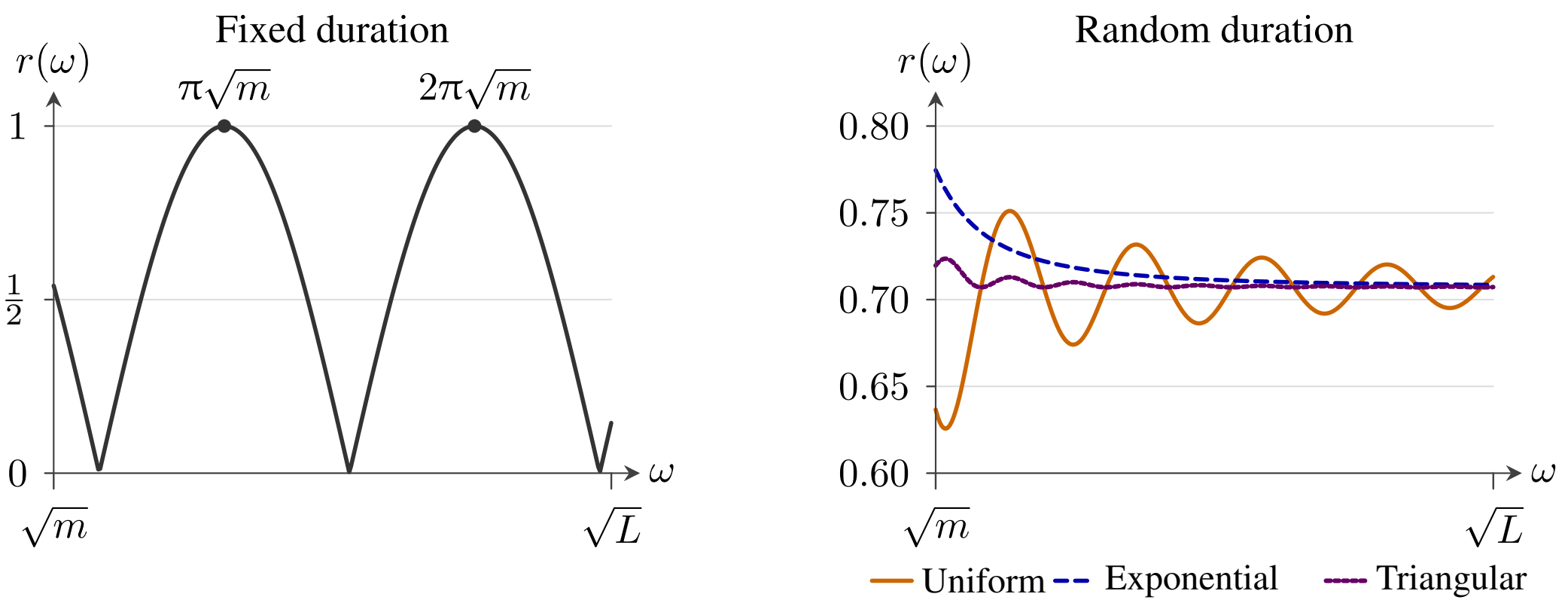
\begin{figure}[!t]
 \centering
\begin{tikzpicture}[
  x=.70cm,y=3.05cm,font=\small,line join=round,line cap=round,
  axis/.style={-{Stealth[length=4pt,width=4pt]},draw=black!75},
  tick/.style={draw=black!75,line width=.45pt},
  grid/.style={draw=black!13,line width=.4pt},
  fixed law/.style={draw=black!80,line width=1.05pt},
  uniform law/.style={draw=orange!80!black,line width=1.05pt},
  exponential law/.style={draw=blue!68!black,line width=1.05pt,densely dashed},
  triangular law/.style={draw=violet!80!black,line width=1.15pt,densely dotted}
]
  \begin{scope}
    \node at (4.5,1.28) {Fixed duration};
    \foreach \level/\ticktext in {0/0,.5/{\frac12},1/1}{
      \draw[grid] (1,\level) -- (8,\level);
      \draw[tick] (1,\level) -- (.89,\level) node[left=1pt] {$\ticktext$};
    }
    \draw[fixed law,domain=1:8,samples=301,variable=\w]
      plot (\w,{abs(cos(deg(\w)))});
    \draw[axis] (1,0) -- (8.36,0) node[right,inner sep=2pt] {$\omega$};
    \draw[axis] (1,0) -- (1,1.10) node[above,inner sep=2pt] {$r(\omega)$};
    \draw[tick] (1,0) -- (1,-.045) node[below=1pt] {$\sqrt m$};
    \draw[tick] (8,0) -- (8,-.045) node[below=1pt] {$\sqrt L$};
    \fill[black!80] (3.141592654,1) circle[radius=1.7pt];
    \fill[black!80] (6.283185307,1) circle[radius=1.7pt];
    \node[above=3pt,inner sep=1pt] at (3.141592654,1) {$\uppi\sqrt m$};
    \node[above=3pt,inner sep=1pt] at (6.283185307,1) {$2\uppi\sqrt m$};
  \end{scope}

  \begin{scope}[shift={(7.75cm,0)}]
    \node at (4.5,1.28) {Random duration};
    \foreach \level/\ticktext in {0/0.60,.25/0.65,.5/0.70,.75/0.75,1/0.80}{
      \draw[grid] (1,\level) -- (8,\level);
      \draw[tick] (1,\level) -- (.89,\level) node[left=1pt] {$\ticktext$};
    }
    \draw[uniform law,domain=1:8,samples=281,variable=\w]
      plot (\w,{5*(sqrt(.5+sin(deg(4*\w))/(8*\w))-.60)});
    \draw[exponential law,domain=1:8,samples=121,variable=\w]
      plot (\w,{5*(sqrt((1+2*\w^2)/(1+4*\w^2))-.60)});
    \draw[triangular law,domain=1:8,samples=281,variable=\w]
      plot (\w,{5*(sqrt(.5+(1-cos(deg(8*\w)))/(64*\w^2))-.60)});
    \draw[axis] (1,0) -- (8.36,0) node[right,inner sep=2pt] {$\omega$};
    \draw[axis] (1,0) -- (1,1.10) node[above,inner sep=2pt] {$r(\omega)$};
    \draw[tick] (1,0) -- (1,-.045) node[below=1pt] {$\sqrt m$};
    \draw[tick] (8,0) -- (8,-.045) node[below=1pt] {$\sqrt L$};

    \begin{scope}[font=\footnotesize]
      \draw[uniform law] (.2,-.31) -- (.7,-.31) node[right=2pt,inner sep=0pt] {Uniform};
      \draw[exponential law] (2.5,-.31) -- (3,-.31) node[right=2pt,inner sep=0pt] {Exponential};
      \draw[triangular law] (5.9,-.31) -- (6.4,-.31) node[right=2pt,inner sep=0pt] {Triangular};
    \end{scope}
  \end{scope}
\end{tikzpicture}

 \caption{Root mean square contraction factor $r(\omega)$ of exact Gaussian HMC on
 \([\sqrt m,\sqrt L]\), \(L/m=64\): fixed duration (left) and random durations (right).}
 \label{fig:integration-time-resonance}
\end{figure}

The sampler draws a duration \(\tau\) from a law \(\rho\). With Verlet
step size \(h\), it converts \(\tau\) into an integer-valued integration time
\(J\), measured in units of \(h\). The physical integration time is \(Jh\). We first explain
the choice of \(\rho\) using exact Gaussian HMC, then define the
rounding rule that determines \(J\).

For a centered Gaussian target with precision matrix \(\Omega\), exact HMC
draws a fresh velocity \(v\sim N(0,I_d)\), follows the Hamiltonian equations
\(\dot x=v\), \(\dot v=-\Omega x\) for a duration \(\tau\) independent of the
position and velocity, and returns the endpoint position. The flow is
computed exactly and preserves the Hamiltonian, so every proposal is
accepted.

\paragraph*{Resonance in exact Gaussian HMC}
Figure~\ref{fig:integration-time-resonance} shows the effect of duration
randomization on contraction in exact Gaussian HMC.
For a Gaussian target, let \(\lambda_i\in[m,L]\) be a precision eigenvalue and
write \(\omega_i=\sqrt{\lambda_i}\). Start two exact HMC trajectories from
different positions with the same initial momentum and the same
duration \(\tau\).
If \(\Delta_i\) is their initial position difference in mode \(i\),
the difference after one transition is
\(\Delta_i^+=\cos(\omega_i\tau)\Delta_i\). For a fixed duration \(s\),
the ratio of the absolute differences is \(|\cos(\omega_i s)|\),
which equals one at a resonance \(\omega_i s\in\uppi\mathbb Z\).

For \(\Delta_i\ne0\), let \(r(\omega_i)\) be the root mean square of
\(\Delta_i^+\) divided by \(|\Delta_i|\). Since the two trajectories
share the same initial momentum and the same duration,
\(\Delta_i^+=\cos(\omega_i\tau)\Delta_i\) is random only through
\(\tau\sim\rho\), and
\[
 r(\omega_i)^2=\E_\rho[\cos^2(\omega_i\tau)]
 =\frac{1+\operatorname{Re}\widehat\rho(2\omega_i)}{2},
\]
where \(\widehat\rho(\omega)=\E_\rho[e^{i\omega\tau}]\) is the
characteristic function of \(\rho\), so that
\(\operatorname{Re}\widehat\rho(\omega)=\E_\rho\cos(\omega\tau)\).
The position difference in mode \(i\) contracts in mean square when
\(r(\omega_i)<1\), that is, when
\(\operatorname{Re}\widehat\rho(2\omega_i)<1\). The doubled frequencies
\(2\omega_i=2\sqrt{\lambda_i}\) all lie in \([2\sqrt m,2\sqrt L]\).
Requiring \(\operatorname{Re}\widehat\rho\le1-\delta_\rho\) on the whole
half-line \([2\sqrt m,\infty)\) gives contraction in every mode by a
factor independent of \(L\), and hence of \(\kappa\); this is the
second condition in Definition~\ref{def:nonresonant}.

\begin{definition}[Non-resonant duration law]\label{def:nonresonant}
Let \(b_\rho>0\) and \(0<\delta_\rho\le1\). A probability law \(\rho\) on \((0,\infty)\)
is non-resonant with constants \(b_\rho,\delta_\rho\) if
\begin{equation}\label{pre:nonresonant}
 \E_\rho[\tau]\le \frac{b_\rho}{\sqrt m},
 \qquad
 \sup_{\omega\ge 2\sqrt m}\operatorname{Re}\widehat\rho(\omega)\le 1-\delta_\rho.
\end{equation}
\end{definition}

The root mean square contraction factor \(r(\omega_i)\) is then at most
\(\sqrt{1-\delta_\rho/2}\) in every mode. The two conditions force $\mathbb E_\rho\tau$ to be of order
$m^{-1/2}$: the first gives $\mathbb E_\rho\tau\le b_\rho/\sqrt m$, and
the second at $\omega=2\sqrt m$, with $1-\cos x\le|x|$, gives
$\delta_\rho\le\mathbb E_\rho[1-\cos(2\sqrt m\,\tau)]\le2\sqrt m\,\mathbb E_\rho\tau$.

\paragraph*{Examples}
The three laws below satisfy Definition~\ref{def:nonresonant} with
\(b_\rho=2\) and \(\delta_\rho=3/4\). The last column uses \(u=\omega/\sqrt m\ge2\).
\begin{center}
\setlength{\tabcolsep}{10pt}
\renewcommand{\arraystretch}{1.15}
\begin{tabular}{@{}lcc@{}}
\toprule
Law of \(\tau\) & \(\mathbb E_\rho\tau\) &
  \(\operatorname{Re}\widehat\rho(\omega)\) \\
\midrule
Uniform on \([0,2/\sqrt m]\) & \(1/\sqrt m\) &
  \(\sin(2u)/(2u)\le 1/4\) \\
\addlinespace[3pt]
Exponential & \(1/\sqrt m\) &
  \(1/(1+u^2)\le 1/5\) \\
\addlinespace[3pt]
Decreasing triangular on \([0,4/\sqrt m]\) & \(4/(3\sqrt m)\) &
  \(2(1-\cos(4u))/(4u)^2\le 1/16\) \\
\bottomrule
\end{tabular}
\end{center}

\begin{figure}[!t]
 \centering
\begin{tikzpicture}[
  font=\small,
  axis/.style={-{Stealth[length=4pt,width=4pt]},draw=black!70},
  highlight/.style={blue!65!black},
  tick/.style={draw=black!70}
]
  \begin{scope}[x=4.7cm,y=1.55cm]
    \node at (0.5,1.43) {Sampled duration};
    \node[anchor=west] at (0.01,1.16) {$f_\tau(t)$};
    \fill[highlight,opacity=.2] (0.25,0) -- (0.25,0.75)
      -- (0.375,0.625) -- (0.375,0) -- cycle;
    \draw[highlight,densely dashed] (0.25,0) -- (0.25,0.75)
      (0.375,0) -- (0.375,0.625);
    \draw[axis] (-0.02,0) -- (1.10,0) node[below] {$t$};
    \draw[axis] (0,-0.055) -- (0,1.15);
    \draw[black!75,line width=.85pt] (0,1) -- (1,0);
    \draw[tick] (-0.018,1) -- (0,1) node[left=2pt] {$2/T$};
    \node[below=2pt] at (0,0) {$0$};
    \foreach \time/\labeltext in {0.25/2h,0.375/3h,1/T}{
      \draw[tick] (\time,0) -- (\time,-0.055)
        node[below=1pt] {$\labeltext$};
    }
  \end{scope}
  \begin{scope}[shift={(6.7cm,0)},x=.53cm,y=1.55cm]
    \node at (4,1.43) {Integration time};
    \node[anchor=west] at (-0.55,1.16) {$\mathbb P(J=j)$};
    \draw[axis] (-0.62,0) -- (8.85,0) node[below] {$j$};
    \draw[axis] (-0.55,-0.055) -- (-0.55,1.15);
    \draw[tick] (-0.68,1) -- (-0.55,1) node[left=1pt] {$7/32$};
    \foreach \j/\mass in {0/15,1/28,2/24,3/20,4/16,5/12,6/8,7/4,8/1}{
      \pgfmathsetmacro{\height}{\mass/28}
      \pgfmathtruncatemacro{\emphasis}{(\j==2)||(\j==3)}
      \ifnum\emphasis=1\relax
        \draw[highlight,line width=1pt] (\j,0) -- (\j,\height);
        \fill[highlight] (\j,\height) circle[radius=2.15pt];
      \else
        \draw[black!64,line width=.7pt] (\j,0) -- (\j,\height);
        \fill[black!64] (\j,\height) circle[radius=1.85pt];
      \fi
      \node[below=3pt] at (\j,0) {$\j$};
    }
  \end{scope}
\end{tikzpicture}

 \caption{Triangular durations on \([0,T]\) and integration times
 \(J\) for \(h=T/8\). A draw in the shaded interval gives \(J=2\) or \(J=3\)
 with equal probability.}
 \label{fig:integration-time}
\end{figure}
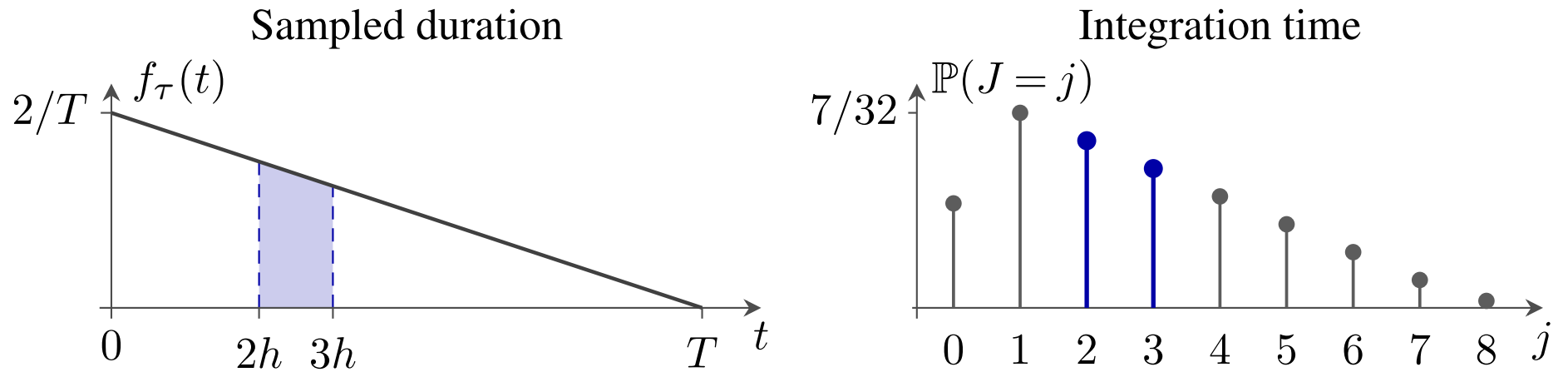

\paragraph*{Randomized rounding}
Rounding down sends every duration in \([jh,(j+1)h)\) to the same
integration time \(j\). On a quadratic potential, \(j\) Verlet steps act
as a rotation in each eigendirection of the Hessian; if \(j\) steps
complete a full period in some eigendirection, every duration in this
interval returns that coordinate to its initial value. We instead choose between
\(\lfloor\tau/h\rfloor\) and \(\lfloor\tau/h\rfloor+1\) with equal
probability. The extra step advances the
phase. In Fourier terms, this step multiplies the characteristic
function of \(\lfloor\tau/h\rfloor\) by \((1+e^{i\theta})/2\), whose modulus
is strictly less than one away from \(2\uppi\mathbb Z\).
For a non-resonant duration law with gap \(\delta_\rho\), we use a step size
satisfying \(h^2m\le\delta_\rho/16\), and call the resulting law of \(J\)
a \emph{non-resonant integration-time law}.
Figure~\ref{fig:integration-time} illustrates the rounding rule for a
triangular duration law.

\subsection{Mixing and gradient complexity}

An initial position law $\nu_0$ is an $L^2$-warm start with bound
$M\ge1$ if it has density $f_{\mathrm{init}}=\mathrm d\nu_0/\mathrm d\pi$ with
\begin{equation}\label{up:warm}
 \|f_{\mathrm{init}}\|_{L^2(\pi)}\le M.
\end{equation}

\begin{theorem}[Mixing from a warm start]\label{up:main}
Suppose $U:\R^d\to\R$ satisfies \ref{ass:rhmc-a1} to \ref{ass:rhmc-a3},
$\nu_0$ is an $L^2$-warm start with bound $M$ as in \eqref{up:warm}, and
$\rho$ is a non-resonant duration law with constants $b_\rho,\delta_\rho$
as in Definition~\ref{def:nonresonant}.
Let $\varepsilon\in(0,1]$ and $D=1+d+\log(M/\varepsilon)$. There are
constants $c_1,c_2>0$ depending only on $b_\rho,\delta_\rho$ such that
if $h=4/(N\sqrt m)$ for an integer
\begin{equation}\label{up:main-N}
 N\ge c_1\bigl[\sqrt{\kappa D}+\sqrt{L_H}\,(D/m)^{3/4}\bigr]
\end{equation}
and $n\ge c_2\log(2M/\varepsilon)$, then the position kernel $P$ with
step size $h$ satisfies
\begin{equation}\label{up:tv-conclusion}
 \operatorname{TV}(\nu_0P^n,\pi)\le\varepsilon.
\end{equation}
\end{theorem}

Explicit values of $c_1$ and $c_2$ are given in
\eqref{up:explicit-constants-c}. Since $h=4/(N\sqrt m)$, condition \eqref{up:main-N} gives
$h\le4/(c_1\sqrt{LD})$, which is stronger than the linear stability
condition $h\sqrt L<2$ for Verlet on quadratic potentials; the proof
also requires small energy error over a full trajectory to control
Metropolis rejections. The number of
transitions is independent of the condition number, which enters only
through the number of gradient evaluations per transition.

\begin{corollary}[Gradient complexity]\label{up:gradient-complexity}
Under the assumptions of Theorem~\ref{up:main}, let $C_{\mathrm{init}}$
be the number of gradient evaluations used to sample the initial law
$\nu_0$ and to set the parameters of the sampler, and let $C_\nabla$ be
the total number of gradient evaluations for initialization and the $n$
transitions. If $N$ and $n$ are the smallest integers permitted by
Theorem~\ref{up:main}, then
\begin{equation}\label{up:complexity-order}
 \mathbb E C_{\nabla}\le\mathbb E C_{\mathrm{init}}+1+
 \mathcal O\Bigl(\bigl[\sqrt{\kappa D}+\sqrt{L_H}\,(D/m)^{3/4}\bigr]
       \log\tfrac{2M}{\varepsilon}\Bigr),
\end{equation}
where the implicit constant depends only on $b_\rho$ and $\delta_\rho$.
\end{corollary}

\begin{proof}
Evaluate the initial gradient once. Each Verlet step reuses the gradient
at its starting point and evaluates one new gradient at its endpoint.
After the Metropolis test, retain the gradient at the returned position.
With stored gradients reused in this way,
\begin{equation}\label{up:exact-count}
 C_{\nabla}=C_{\mathrm{init}}+1+\sum_{i=1}^n J_i.
\end{equation}
Randomized rounding and Definition~\ref{def:nonresonant} give
\[
 \mathbb E J_i\le\frac{\mathbb E_\rho\tau}{h}+\frac12
       =\frac{N\sqrt m}{4}\mathbb E_\rho\tau+\frac12
       \le\frac{b_\rho N}{4}+\frac12,
\]
and the choices of $N$ and $n$ give \eqref{up:complexity-order}.
\end{proof}

\begin{remark}[Convexity and a transport-entropy inequality suffice]\label{up:talagrand-remark}
Strong convexity enters the proof of Theorem~\ref{up:main} only through
Talagrand's inequality \eqref{pre:smooth-bound} and through the choice of
the integration time $T=4/\sqrt m$ in the proof. Elsewhere only
convexity of $U$ is used, for the sign of the transport terms in
Lemma~\ref{ent:weighted-entropy} and for the tail bound
\eqref{up:gibbs-tail}. Consequently Theorem~\ref{up:main} holds, with the
same constants, for every convex $U$ with a Lipschitz gradient and Hessian
whose Gibbs law satisfies \eqref{pre:smooth-bound} for some $m>0$, with
this $m$ in place of the strong-convexity constant throughout. By the
Holley-Stroock perturbation lemma \cite{HolleyStroock} and the
Otto-Villani theorem \cite{OttoVillani}, this includes bounded
perturbations of strongly convex potentials that remain convex with a
Lipschitz gradient and Hessian. The lower bounds of
Theorem~\ref{lower:thm} and Proposition~\ref{gaussian:cost-lower} use
strong convexity as stated.
\end{remark}
\subsection{Lower bound}

For $T>0$, consider duration laws $\rho$ satisfying
\begin{equation}\label{lower:time-law}
 \mathbb P(\tau<t)\leq \frac{2t}{T}\quad(t\geq0),
 \qquad \frac{T}{3}\leq\mathbb E\tau<\infty .
\end{equation}
The first condition bounds the probability of drawing a short duration,
and the second bounds the mean duration from below. With $m=\kappa^{-1}$,
the uniform and exponential laws of Section~\ref{sec:nonresonant}
satisfy \eqref{lower:time-law} with $T=2\sqrt\kappa$, and the decreasing
triangular law with $T=4\sqrt\kappa$; this is the duration scale of
Theorem~\ref{up:main}, where $\mathbb E_\rho\tau$ is of order
$m^{-1/2}=\sqrt\kappa$.

\begin{theorem}[Necessity of Hessian regularity for fixed step sizes]\label{lower:thm}
Set
\[
q_*=10^{-3},\qquad
d_f=\left\lceil 6q_*^{-1}\log\kappa\right\rceil,\qquad d=d_f+2,
\qquad
\Lambda_\kappa=\log(e+\log\kappa).
\]
For every $\kappa\ge64$, $0<h\le1$, and $T\in[2\sqrt\kappa,4\sqrt\kappa]$, there is a
potential $U_h:\R^d\to\R$ satisfying \ref{ass:rhmc-a1} and
\ref{ass:rhmc-a2} with optimal constants $m=\kappa^{-1}$ and $L=1$,
and a Borel set $A_+$ with $\pi(A_+)=1/2$, such that the following holds.
Let $\nu=2\mathbf1_{A_+}\pi$, and let $P$ be the position kernel of
Algorithm~\ref{alg:rhmc-transition} with step size $h$ and a duration law
satisfying \eqref{lower:time-law}. For every integer $n\ge0$,
$\TV(\nu P^n,\pi)\le1/4$ implies
\begin{equation}\label{lower:cost}
 \mathbb E N_\nabla\ge
 c_3^{2/3}c_4^{1/3}\frac{\kappa^{2/3}}{\Lambda_\kappa^{1/3}},
\end{equation}
where $N_\nabla$ is the number of gradient evaluations in these $n$ transitions,
and $c_3,c_4>0$ are absolute constants.
\end{theorem}

Since $d=d_f+2$ grows only logarithmically in $\kappa$ and
$d\nu/d\pi\le2$,
\eqref{lower:cost} rules out a uniform $\sqrt\kappa$ gradient bound
multiplied by any fixed polynomial in $d$ and $\log\kappa$ for this
algorithm.

\section{Proof of Theorem~\ref{up:main}}\label{sec:upper}

We first fix the constants and parameters used in the proof, then
state a general mixing criterion (Proposition~\ref{up:warm-mixing}),
and establish the sampler's lower spectral bound
(Proposition~\ref{up:spectral-bound}) and restricted conductance lower
bound (Proposition~\ref{up:restricted-conductance}).

Fix a duration cutoff $\tau_*>0$ such that
$p_S:=\mathbb P(\tau\le\tau_*)\ge1-\delta_\rho/4$, and set
$T_*=\tau_*+4/\sqrt m$. Markov's inequality permits
$\tau_*=4b_\rho/(\delta_\rho\sqrt m)$. For a bounded duration law one may
instead take its upper support endpoint, giving $p_S=1$.

For each integer $N\ge16/\sqrt{\delta_\rho}$, put $h=4/(N\sqrt m)$.
Let $\mu_J$ be the full law of $J$, let $w$ be its law conditional on
$\tau\le\tau_*$, and let $q$ be the law obtained by the same rounding
rule from decreasing triangular durations on $[0,4/\sqrt m]$.
Choose $0<c,\beta_P\le1$, uniformly in these $N$, such that for every
$\theta\in[-\uppi,\uppi]$,
\begin{equation}\label{up:law-constants}
 1-\operatorname{Re}\widehat w(\theta)
 \ge c\bigl(1-\operatorname{Re}\widehat q(\theta)\bigr),
 \qquad 1+\operatorname{Re}\widehat\mu_J(\theta)\ge\beta_P.
\end{equation}
Here the hats denote characteristic functions.
Lemma~\ref{up:rounded-fourier} gives explicit positive choices for every
non-resonant duration law; Lemma~\ref{up:example-fourier} gives the
choices for the three examples. Define
\begin{equation}\label{up:explicit-constants}
 \eta=\frac c{600},\qquad \phi_P=\frac{p_Sc}{60}.
\end{equation}
The first quantity bounds the energy increase and the second is the
restricted conductance bound.

\paragraph*{Explicit constants}
In the setting of Theorem~\ref{up:main}, take the cutoff
$\tau_*=4b_\rho/(\delta_\rho\sqrt m)$, so that $T_*=C_T/\sqrt m$ with
$C_T=4(b_\rho+\delta_\rho)/\delta_\rho$, set
\begin{equation}\label{up:parameters-a}
 s=\frac{\varepsilon^4}{64M^4},\qquad T=\frac4{\sqrt m},\qquad
 R=3d\log2+2\log\frac1{\eta s},
\end{equation}
and put $C_R=2\log(64/\eta)$. The constants in Theorem~\ref{up:main} are
\begin{equation}\label{up:explicit-constants-c}
 c_1=\max\Bigl\{\frac{16}{\sqrt{\delta_\rho}},\;
   \sqrt{\frac{48}{\eta}}\,\max\bigl\{\sqrt{C_R},\,\sqrt{C_T}\,C_R^{3/4}\bigr\}\Bigr\},
 \qquad
 c_2=\frac{6}{\beta_P\phi_P^2}.
\end{equation}
By Lemma~\ref{up:rounded-fourier}, $c$ and $\beta_P$ are bounded
below by positive quantities depending only on $b_\rho,\delta_\rho$,
and $p_S\ge1-\delta_\rho/4$ by the choice of $\tau_*$. The constant
$c_2$ depends on $p_S$ through $\phi_P$; replacing $p_S$ by this lower
bound gives values of $\eta$, $\phi_P$, $c_1$ and $c_2$ that depend only
on $b_\rho,\delta_\rho$. The proof below uses
these choices through three consequences, which we now derive: the
tail bound $\Pi(H>R)\le\eta s$, and the conditions
\begin{align}
 N&\ge\frac{16}{\sqrt{\delta_\rho}},\qquad
 h^2(LR+L_HT_*R^{3/2})\le\frac\eta3,\label{up:parameters-b}\\
 n&\ge\frac{2}{\beta_P\phi_P^2}\log\frac{8M^2}{\varepsilon^3}.
 \label{up:iterations}
\end{align}

Under \ref{ass:rhmc-a1}, the Boltzmann-Gibbs law satisfies the tail
bound
\begin{equation}\label{up:gibbs-tail}
 \Pi(H>R)\le2^{3d/2}e^{-R/2}\qquad(R\ge0).
\end{equation}
To see this, translate $x_*$ to the origin and put $u(x)=U(x)-U(x_*)$.
For $0<t\le1$, convexity gives
$u(tx)\le t u(x)$, so $Z(t)=\int_{\R^d} e^{-t u(x)}\,\mathrm dx$ satisfies
\[
 Z(t)\le\int_{\R^d} e^{-u(tx)}\,\mathrm dx=t^{-d}Z(1).
\]
Consequently $\mathbb E_\pi e^{u(X)/2}\le2^d$.
The independent kinetic energy contributes
$\mathbb E_{\gamma_d}e^{|V|^2/4}=2^{d/2}$.
Markov's inequality gives \eqref{up:gibbs-tail}, and with $R$ as in
\eqref{up:parameters-a} it yields $\Pi(H>R)\le\eta s$. Moreover
$R\le C_RD$, since $\eta\le1/600$.

For $N$, $n$ and $h=T/N$ as in Theorem~\ref{up:main},
$\kappa D\ge1$ gives $N\ge c_1\ge16/\sqrt{\delta_\rho}$; the
bounds $R\le C_RD$, $T_*=C_T/\sqrt m$ and $\sqrt{a}+\sqrt{b}\ge\sqrt{a+b}$
give
\[
 N^2\ge\frac{48}{\eta}\Bigl(C_R\kappa D+L_HC_TC_R^{3/2}\frac{D^{3/2}}{m^{3/2}}\Bigr)
 \ge\frac{48}{\eta m}\bigl(LR+L_HT_*R^{3/2}\bigr),
\]
which is the second condition in \eqref{up:parameters-b} with
$h^2=16/(N^2m)$; and \eqref{up:iterations} follows from
$\log(8M^2/\varepsilon^3)\le3\log(2M/\varepsilon)$, valid for $M\ge1$.
The remainder of this section proves \eqref{up:tv-conclusion} under
\eqref{up:parameters-b} and \eqref{up:iterations}.

For a probability measure $\pi$ and a $\pi$-reversible kernel $P$,
write $Pf(x)=\int_{\R^d}f(y)P(x,\mathrm dy)$ and
\begin{align*}
 Q_P(A,A^c)&=\int_A P(x,A^c)\,\pi(\mathrm dx),\\
 \mathcal E_P(f)&=\frac12\int_{\R^d}\int_{\R^d}
           (f(x)-f(y))^2\,P(x,\mathrm dy)\pi(\mathrm dx)
           =\langle f,(I-P)f\rangle_{L^2(\pi)},
\end{align*}
where $A$ is a Borel subset of $\R^d$, $f\in L^2(\pi)$, and $I$ is the
identity operator. 

\begin{mainproposition}[Mixing from restricted conductance]\label{up:warm-mixing}
Let $P$ be a Markov kernel on $\R^d$, reversible with respect to a
probability measure $\pi$. Suppose that, for some $0<s\le1/2$ and
$0<\phi_P,\beta_P\le1$,
\begin{equation}\label{up:criterion-conductance}
 Q_P(A,A^c)\ge\phi_P\pi(A)
 \quad\text{for every Borel }A\text{ with }s\le\pi(A)\le\tfrac12,
\end{equation}
\begin{equation}\label{up:criterion-spectral}
 \|f\|_{L^2(\pi)}^2-\|Pf\|_{L^2(\pi)}^2
 \ge \beta_P\mathcal E_P(f),\qquad f\in L^2(\pi).
\end{equation}
If a probability measure $\nu_0\ll\pi$ has density $f_{\mathrm{init}}$ with
$\|f_{\mathrm{init}}\|_{L^2(\pi)}\le M$, then for $0<\varepsilon\le1$
and every integer $k\ge0$,
\[
 \operatorname{TV}(\nu_0P^k,\pi)
 \le\frac{\varepsilon}{2}
   +\frac12\left(e^{-k\beta_P\phi_P^2/2}M_\infty
                         +2sM_\infty^2\right)^{1/2},
 \qquad M_\infty=\frac{4M^2}{\varepsilon}.
\]
\end{mainproposition}

Proposition~\ref{up:warm-mixing} is proved in
Section~\ref{sec:rhmc-criterion}. For the sampler, the following proposition
establishes \eqref{up:criterion-spectral}.

\begin{mainproposition}[Lower spectral bound]\label{up:spectral-bound}
Assume \ref{ass:rhmc-a1} and \ref{ass:rhmc-a2}. Let $\rho$ be a
non-resonant duration law with constants $b_\rho,\delta_\rho$.
Let $N\ge1$ be an integer and $h=4/(N\sqrt m)$.
If $h^2m\le\delta_\rho/16$, then $P$ is reversible with respect to $\pi$ and,
with $\beta_P$ satisfying \eqref{up:law-constants},
\[
 \|f\|_{L^2(\pi)}^2-\|Pf\|_{L^2(\pi)}^2
 \ge \beta_P\mathcal E_P(f),\qquad f\in L^2(\pi).
\]
\end{mainproposition}

Proposition~\ref{up:spectral-bound} is proved in
Section~\ref{sec:rhmc-spectral-proof}. It remains to verify the restricted conductance hypothesis
\eqref{up:criterion-conductance}.

\begin{mainproposition}[Restricted conductance of the sampler]\label{up:restricted-conductance}
Suppose $U$ satisfies \ref{ass:rhmc-a1}--\ref{ass:rhmc-a3} and $\rho$ is a
non-resonant duration law with constants $b_\rho,\delta_\rho$.
Choose $\tau_*,c$ as in \eqref{up:law-constants}, and put
$\eta=c/600$, $\phi_P=p_Sc/60$.
Let $0<s\le1/2$ and $R\ge1$ satisfy $\Pi(H>R)\le\eta s$.
If $h=4/(N\sqrt m)$ for an integer $N\ge16/\sqrt{\delta_\rho}$ and
\[
 h^2(LR+L_HT_*R^{3/2})\le\eta/3,
\]
then
\[
 Q_P(A,A^c)\ge\phi_P\pi(A)
 \quad\text{for every Borel }A\subset\R^d
       \text{ with }s\le\pi(A)\le\tfrac12.
\]
\end{mainproposition}

Proposition~\ref{up:restricted-conductance} is proved in
Section~\ref{sec:rhmc-conductance-proof}.

Table~\ref{tab:proof-parameters} records the roles of the constants and
cutoffs. The remaining quantities are the expressions
$T_*=\tau_*+4/\sqrt m$, $\eta=c/600$, and $\phi_P=p_Sc/60$.

\begin{center}
\begin{minipage}{\linewidth}
\centering
\setlength{\tabcolsep}{8pt}
\renewcommand{\arraystretch}{1.15}
\captionof{table}{Constants and cutoffs in the proof of Theorem~\ref{up:main}.}
\label{tab:proof-parameters}
\medskip
\begin{tabular}{@{}lll@{}}
\toprule
Symbol & Condition & Role \\
\midrule
$c$ & \eqref{up:law-constants} & Fourier comparison \\
$\beta_P$ & \eqref{up:law-constants} & Lower spectral bound \\
$p_S$ & $\mathbb P(\tau\le\tau_*)$ & Retained duration probability \\
$\tau_*$ & $p_S\ge1-\delta_\rho/4$ & Duration cutoff \\
$R$ & $R\ge1$, $\Pi(H>R)\le\eta s$ & Energy cutoff \\
\bottomrule
\end{tabular}
\end{minipage}
\end{center}

For the three duration laws, Table~\ref{tab:explicit-examples} gives
choices valid for every admissible $N$. They are proved in
Lemma~\ref{up:example-fourier}.

\begin{center}
\begin{minipage}{\linewidth}
\centering
\setlength{\tabcolsep}{12pt}
\renewcommand{\arraystretch}{1.15}
\captionof{table}{Choices for the three duration laws.}
\label{tab:explicit-examples}
\medskip
\begin{tabular}{@{}lcccc@{}}
\toprule
Law & $c$ & $\beta_P$ & $\sqrt m\,\tau_*$ & $p_S$ \\
\midrule
Triangular & $1$ & $1$ & $4$ & $1$ \\
Uniform & $1/2$ & $3/4$ & $2$ & $1$ \\
Exponential & $3/5$ & $1$ & $8$ & $1-e^{-8}$ \\
\bottomrule
\end{tabular}
\end{minipage}
\end{center}

\edef\theoremproofindent{\the\parindent}
\par\medskip\noindent\begin{minipage}{\linewidth}
\setlength{\parindent}{\theoremproofindent}
\begin{proof}[Proof of Theorem~\ref{up:main}]
We apply the restricted conductance criterion of
Proposition~\ref{up:warm-mixing} with $s$ as in \eqref{up:parameters-a}.
Assumptions \ref{ass:rhmc-a1}--\ref{ass:rhmc-a3}, the non-resonance of
$\rho$, and $M\ge1$, $0<\varepsilon\le1$ give the target and duration-law
hypotheses of Propositions~\ref{up:spectral-bound}
and~\ref{up:restricted-conductance}.
The cutoff satisfies $\Pi(H>R)\le\eta s$ by the choice of $R$ in
\eqref{up:parameters-a}, and
\eqref{up:parameters-b} is the energy condition in
Proposition~\ref{up:restricted-conductance}.
It also gives $N\ge16/\sqrt{\delta_\rho}$, hence
$h^2m=16/N^2\le\delta_\rho/16$ as required by
Proposition~\ref{up:spectral-bound}.
Thus both propositions apply. Lemma~\ref{up:reversibility} gives
reversibility; Propositions~\ref{up:spectral-bound} and
\ref{up:restricted-conductance} give
\eqref{up:criterion-spectral} and \eqref{up:criterion-conductance}, respectively.
The constants satisfy $0<\phi_P,\beta_P\le1$.

Since $0<s\le1/64$ and $\nu_0$ is an $L^2$-warm start with bound $M$,
all hypotheses of Proposition~\ref{up:warm-mixing} are satisfied.
With $M_\infty=4M^2/\varepsilon$ and
$2sM_\infty^2=\varepsilon^2/2$, it gives
\[
 \operatorname{TV}(\nu_0P^n,\pi)
 \le\frac{\varepsilon}{2}
   +\frac12\left(e^{-n\beta_P\phi_P^2/2}M_\infty
                         +\frac{\varepsilon^2}{2}\right)^{1/2}.
\]
The right side is at most $\varepsilon$ when
\begin{equation}\label{up:proof-iterations}
 n\ge\frac{2}{\beta_P\phi_P^2}\log\frac{8M^2}{\varepsilon^3},
\end{equation}
as required by \eqref{up:iterations}.
\end{proof}
\end{minipage}\par

\section{Mixing from restricted conductance}\label{sec:rhmc-criterion}

We prove Proposition~\ref{up:warm-mixing}. The argument applies to a
general reversible Markov kernel. For bounded real $f$, write
$\operatorname{osc}(f)$ for its essential supremum minus its essential
infimum under $\pi$.

\begin{proof}[Proof of Proposition~\ref{up:warm-mixing}]
By \cite[Remark~41 and Lemma~42]{ALPW},
\eqref{up:criterion-conductance} gives the weak Poincar\'e inequality
\begin{equation}\label{up:weak-poincare-equation}
 \operatorname{Var}_\pi(f)
    \le\frac2{\phi_P^2}\mathcal E_P(f)
                    +2s\operatorname{osc}(f)^2
\end{equation}
for every bounded real measurable $f$.
Indeed, in the notation of \cite[Definition~34]{ALPW},
\eqref{up:criterion-conductance} says that
$\Phi_P^{\mathrm W}(s)\ge\phi_P$.
The cited results bound the weak Poincar\'e coefficient at parameter
$r$ by $2[\Phi_P^{\mathrm W}(r/2)]^{-2}$; take $r=2s$.
To meet their assumption of a positive profile at every threshold,
apply them first to $(1-\delta)P(x,\mathrm dy)+\delta\pi(\mathrm dy)$
and let $\delta\downarrow0$.

First let $f_0$ be a probability density relative to $\pi$ with
$0\le f_0\le M_\infty$. Reversibility identifies the density after
$k$ steps with $f_k=P^kf_0$, and the Markov property gives
$0\le f_k\le M_\infty$.  Write
$V_k=\|f_k-1\|_{L^2(\pi)}^2$.
By \eqref{up:criterion-spectral},
\[
 V_k-V_{k+1}
  =\|f_k-1\|_{L^2(\pi)}^2-\|P(f_k-1)\|_{L^2(\pi)}^2
  \ge \beta_P\mathcal E_P(f_k).
\]
Using \eqref{up:weak-poincare-equation} and
$\operatorname{osc}(f_k)\le M_\infty$ yields
\[
 V_{k+1}\le(1-\beta_P\phi_P^2/2)V_k+\beta_P\phi_P^2sM_\infty^2.
\]
Since $V_0=\int_{\R^d} f_0^2\,\mathrm d\pi-1\le M_\infty-1\le M_\infty$,
iteration and Cauchy-Schwarz give
\begin{equation}\label{up:bounded-warm-mixing}
 V_k\le e^{-k\beta_P\phi_P^2/2}M_\infty+2sM_\infty^2,
 \qquad
 \operatorname{TV}(f_0\pi P^k,\pi)\le\tfrac12\sqrt{V_k}.
\end{equation}
For the $L^2$-warm initial law, truncate and renormalize
$f_{\mathrm{init}}$ as follows.
\[
 M_{\mathrm{cap}}=\frac{2M^2}{\varepsilon},\qquad
 f_{\mathrm{clip}}=\min\{f_{\mathrm{init}},M_{\mathrm{cap}}\},\qquad
 Z_{\mathrm{cap}}=\int_{\R^d} f_{\mathrm{clip}}\,\mathrm d\pi,\qquad
 f_0^{\mathrm{cap}}=\frac{f_{\mathrm{clip}}}{Z_{\mathrm{cap}}}.
\]
The $L^2$ assumption gives
\begin{align}
 1-Z_{\mathrm{cap}}
  &=\int_{\R^d}(f_{\mathrm{init}}-M_{\mathrm{cap}})_+\,\mathrm d\pi
       \le\frac{M^2}{M_{\mathrm{cap}}}=\frac\varepsilon2,
       \label{up:cap-mass}\\
 \operatorname{TV}(f_{\mathrm{init}}\pi,f_0^{\mathrm{cap}}\pi)
  &\le1-Z_{\mathrm{cap}}\le\varepsilon/2,\label{up:cap-tv}\\
 0\le f_0^{\mathrm{cap}}
  &\le\frac{M_{\mathrm{cap}}}{1-\varepsilon/2}
       \le\frac{4M^2}{\varepsilon}=:M_\infty.
       \label{up:cap-density}
\end{align}
For \eqref{up:cap-mass}, use
$(f_{\mathrm{init}}-M_{\mathrm{cap}})_+\le f_{\mathrm{init}}^2/M_{\mathrm{cap}}$.
For \eqref{up:cap-tv}, the two probability densities dominate the common
subdensity $f_{\mathrm{clip}}$, whose mass is $Z_{\mathrm{cap}}$.
Apply \eqref{up:bounded-warm-mixing} to $f_0^{\mathrm{cap}}$ and use
\eqref{up:cap-tv} together with total variation contraction under $P$
to obtain the stated bound.
\end{proof}

\section{Fourier comparison and the lower spectral bound}\label{sec:rhmc-invariance}

The triangular law $q$ supplies the reference kernel for the entropy
argument. The conditional law $w$ restricts trajectory lengths so that
the energy error can be controlled. The full law $\mu_J$ defines the
sampler and is used to establish its lower spectral bound.

\subsection{Reversibility}

\begin{lemma}[Reversibility]\label{up:reversibility}
Suppose $U$ is continuously differentiable and
$Z<\infty$. Then each $K_j$ is reversible
with respect to $\pi$, as is every state-independent mixture of the
kernels $K_j$.
\end{lemma}

\begin{proof}
The additive constant in $H$ can be chosen arbitrarily.
Write $\mathsf{Kick}_s(x,v)=(x,v-s\nabla U(x))$ for the kick and
$\mathsf{Drift}_s(x,v)=(x+sv,v)$ for the drift. Both are invertible shears that
preserve Lebesgue volume. The Verlet update \eqref{up:verlet} is
$\Phi_h=\mathsf{Kick}_{h/2}\circ \mathsf{Drift}_h\circ \mathsf{Kick}_{h/2}$. For the momentum flip
$F(x,v)=(x,-v)$,
\[
 F\circ \mathsf{Kick}_s\circ F=\mathsf{Kick}_{-s},\qquad F\circ \mathsf{Drift}_s\circ F=\mathsf{Drift}_{-s}.
\]
The symmetric order of the kick and drift therefore gives
$F\circ\Phi_h\circ F=\Phi_h^{-1}$. It follows that
$\Psi_j=F\circ\Phi_h^{\circ j}$ is a volume-preserving involution for
every integer $j\ge0$.

Let $p_\Pi$ be the Lebesgue density of $\Pi$ and define
\[
 a_j(z)=p_\Pi(z)\alpha_j(z)
       =\min\{p_\Pi(z),p_\Pi(\Phi_h^{\circ j}(z))\}
       =\min\{p_\Pi(z),p_\Pi(\Psi_j(z))\}.
\]
The last equality holds because the momentum flip preserves $H$.
It also preserves position, so $Q_j(z)$ is the position component
of $\Psi_j(z)$. For measurable $A,A'\subseteq\R^d$, the accepted part
of $\int_A K_j(x,A')\,\pi(\mathrm dx)$ is
\[
 \int_{\R^{2d}}\mathbf1_A(x)\mathbf1_{A'}(Q_j(z))a_j(z)\,\mathrm dz.
\]
The substitution $z\mapsto\Psi_j(z)$ exchanges the two position
arguments and preserves $a_j(z)\,\mathrm dz$, since
$a_j(\Psi_j(z))=a_j(z)$. The rejected part is
supported on the diagonal. Hence $\pi(\mathrm dx)K_j(x,\mathrm dy)$ is
symmetric. The mixture $P$ is therefore reversible with respect to $\pi$.
\end{proof}

\subsection{Comparison of adjusted kernels}

For a law $w=(w_j)_{j\ge0}$ on integration times, write $P_w=\sum_jw_jK_j$
and $\widehat w(\theta)=\sum_jw_je^{ij\theta}$. For real $f\in L^2(\pi)$,
write $f^\sharp(x,v)=f(x)$
and define
\[
 \operatorname{Rej}_w(f)=\sum_jw_j\int_{\R^{2d}}
       (p_\Pi(z)-a_j(z))f^\sharp(z)^2\,\mathrm dz.
\]

\begin{lemma}[Comparison of adjusted kernels]\label{up:law-comparison}
Let $w$ and $q$ be laws on the nonnegative integers. Suppose for
every $\theta\in[-\uppi,\uppi]$ that
\begin{equation}\label{up:fourier-comparison}
 1-\operatorname{Re}\widehat w(\theta)
 \ge c\big(1-\operatorname{Re}\widehat q(\theta)\big),
 \qquad \operatorname{Re}\widehat w(\theta)\ge -1+\beta,
\end{equation}
where $c>0$ and $0<\beta\le1$. The second inequality alone gives
$P_w\succeq(-1+\beta)I$ on $L^2(\pi)$, and the first alone gives
\begin{equation}\label{up:dirichlet-comparison}
 \langle f,(I-P_w)f\rangle_\pi
 \ge c\langle f,(I-P_q)f\rangle_\pi-\operatorname{Rej}_w(f).
\end{equation}
\end{lemma}

\begin{proof}
On phase-space Lebesgue $L^2$, the operator $\mathcal U g=g\circ\Phi_h$
is unitary. Set $S_w=\sum_j w_j(\mathcal U^j+\mathcal U^{-j})/2$.
For $f\in L^2(\pi)$ and $u>0$, let $g_u=f^\sharp\mathbf1_{\{p_\Pi\ge u\}}$.
On this superlevel set, $1\le p_\Pi/u$, so
\[
 \|g_u\|_{L^2(\mathrm dz)}^2
 \le u^{-1}\int_{\R^{2d}} |f^\sharp(z)|^2p_\Pi(z)\,\mathrm dz
 =u^{-1}\|f\|_{L^2(\pi)}^2.
\]
In particular, the spectral theorem applies to $g_u$. The identity
\[
 a_j(z)=\int_0^\infty\mathbf1_{\{p_\Pi(z)\ge u\}}
                  \mathbf1_{\{p_\Pi(\Phi_h^{\circ j}(z))\ge u\}}
                  \,\mathrm du
\]
writes the accepted part of the quadratic form as an integral over
density levels. Adding the rejection part gives
\begin{equation}\label{up:layer-cake}
 \langle f,P_wf\rangle_\pi
  =\int_0^\infty\langle g_u,S_wg_u\rangle\,\mathrm du+\operatorname{Rej}_w(f).
\end{equation}
Indeed, the integral for each $j$ is the accepted contribution
$\int_{\R^{2d}} a_j(z) f^\sharp(z) f^\sharp(\Phi_h^{\circ j}(z))\,\mathrm dz$.
It is absolutely integrable with absolute integral at most $\|f\|_{L^2(\pi)}^2$,
by $2|ab|\le a^2+b^2$, the two bounds on $a_j$, and volume preservation.
Also $\int_0^\infty\|g_u\|_{L^2(\mathrm dz)}^2\,\mathrm du
=\|f\|_{L^2(\pi)}^2$ by Tonelli.

The spectral theorem for $\mathcal U$ turns \eqref{up:fourier-comparison}
into $I-S_w\succeq c(I-S_q)$ and $S_w\succeq(-1+\beta)I$.
The latter and $\operatorname{Rej}_w(f)\ge0$ give the claimed lower spectral
bound. Subtracting \eqref{up:layer-cake} from $\|f\|_{L^2(\pi)}^2$ and using the
same identity with $q$ gives \eqref{up:dirichlet-comparison}, since
$c\operatorname{Rej}_q(f)\ge0$.
\end{proof}

\subsection{Fourier bounds after rounding}

For comparison, use the decreasing triangular density $2(T-t)/T^2$
on $[0,T]$. For $T=Nh$, use the
randomized rounding in Algorithm~\ref{alg:rhmc-transition} with this
density. An interval $[jh,(j+1)h)$ has probability
$[2(N-j)-1]/N^2$, split equally between $J=j$ and $J=j+1$.
Thus the probabilities $q_j=\mathbb P(J=j)$ are
\begin{equation}\label{up:triangular}
 q_0=\frac{2N-1}{2N^2},\qquad
 q_j=\frac{2(N-j)}{N^2}\quad(1\le j<N),\qquad
 q_N=\frac1{2N^2}.
\end{equation}
Write $P_{\mathrm{ref}}=P_q$ for this kernel. The proof of Theorem~\ref{up:main}
uses $T=4/\sqrt m$ for $P_{\mathrm{ref}}$ and compares it with the kernel of the
sampler for the given law $\rho$.

\begin{lemma}[Fourier bounds after rounding]\label{up:rounded-fourier}
Let $\rho$ be non-resonant with constants $b_\rho,\delta_\rho$. Choose
$\tau_*>0$ with $\mathbb P(\tau>\tau_*)\le\delta_\rho/4$, and let
$w,q,\mu_J$ be as in \eqref{up:law-constants}.
Suppose $h=4/(N\sqrt m)$ and $h^2m\le\delta_\rho/16$.
Set
\[
 S=\sqrt m\,\tau_*,\;\; u_0=\frac1{S+1},\;\;
 k=2\lceil S+1\rceil,\;\; \zeta=\frac{\delta_\rho}{4k^2},\quad
 c_0=\min\left\{\frac{\delta_\rho}{64},\,\zeta,\,1-\cos(\zeta/2)\right\}.
\]
For the full law, put
\[
 k_\rho=2\lceil2(b_\rho+1)\rceil,\qquad
 \zeta_\rho=\frac{\delta_\rho}{2k_\rho^2}.
\]
Then \eqref{up:law-constants} holds with
$c=c_0/8$ and $\beta_P=1-\cos(\zeta_\rho/2)$, both positive.
\end{lemma}

\begin{proof}
Put $Y=\sqrt m\,\tau$. The omitted mass is at most $\delta_\rho/4$.
This omitted mass contributes at most $\delta_\rho/2$ to
$\mathbb E[1-\cos(uY)]$.
Thus under the conditional law,
$\mathbb E[1-\cos(uY)]\ge\delta_\rho/2$ whenever $|u|\ge2$.
For the conditional estimates below, expectations are under $\tau\le\tau_*$.
Put $u_0=1/(S+1)$, and let $k$ be the least even integer
with $ku_0\ge2$. For $|u|\ge u_0$, the inequalities
\[
 1-\cos(kx)\le k^2(1-\cos x),\qquad
 1-\cos(kx)\le k^2(1+\cos x)
\]
give $1\pm\mathbb E\cos(uY)\ge\delta_\rho/(2k^2)$. The second inequality
follows from the first by shifting $x$ by $\uppi$.

We first prove
\[
 1-\operatorname{Re}\widehat w(\theta)
 \ge c_0\min\{\theta^2/(h^2m),1\},\quad -\uppi\le\theta\le\uppi,
\]
by considering three frequency ranges.
The step-size condition and $\delta_\rho\le1$ imply $h\sqrt m\le1/4$.
For small frequencies, we use a lower bound on the second moment of $J$.
If $|\theta|\le u_0h\sqrt m$, then $h\sqrt m\,J\le Y+h\sqrt m\le S+1$,
so $|\theta J|\le1$. Since $1-\cos(2Y)\le2Y^2$, the gap at $u=2$
gives $\mathbb EY^2\ge\delta_\rho/4$.
Also $|h\sqrt m\,J-Y|\le h\sqrt m$ implies
$Y^2\le2h^2mJ^2+2h^2m$. Thus
$h^2m\,\mathbb EJ^2\ge\delta_\rho/8-h^2m\ge\delta_\rho/16$.
Using $1-\cos x\ge x^2/4$ for $|x|\le1$ now gives
\[
 1-\operatorname{Re}\widehat w(\theta)
 \ge\frac{\delta_\rho}{64}\frac{\theta^2}{h^2m},
 \qquad \operatorname{Re}\widehat w(\theta)\ge\cos 1.
\]
For intermediate frequencies, non-resonance gives a gap that exceeds
the rounding error. If $|\theta|\ge u_0h\sqrt m$ and
$|\theta|\le\zeta=\delta_\rho/(4k^2)$, then
\[
 \big|\mathbb E\cos(\theta J)-\mathbb E\cos(\theta Y/(h\sqrt m))\big|\le|\theta|.
\]
Consequently, $1-\operatorname{Re}\widehat w(\theta)\ge\zeta$ and
$1+\operatorname{Re}\widehat w(\theta)\ge\zeta$.

For the remaining frequencies,
$\zeta\le|\theta|\le\uppi$, the independent Bernoulli step gives the
factor $|(1+e^{i\theta})/2|\le\cos(\zeta/2)<1$ in
$|\widehat w(\theta)|$.

With the stated $c_0$, these estimates give
$1-\operatorname{Re}\widehat w(\theta)\ge
c_0\min\{\theta^2/(h^2m),1\}$ and
$\operatorname{Re}\widehat w(\theta)\ge-1+[1-\cos(\zeta/2)]$.
Since $Nh\sqrt m=4$, the reference weights satisfy
\[
 1-\operatorname{Re}\widehat q(\theta)
 \le \min\{N^2\theta^2/2,2\}
 \le 8\min\{\theta^2/(h^2m),1\}.
\]
This proves the first inequality in \eqref{up:law-constants} with $c=c_0/8$.

For the full law, put $u_1=1/[2(b_\rho+1)]$.
If $|\theta|\le u_1h\sqrt m$, then
\[
 \operatorname{Re}\widehat\mu_J(\theta)
 \ge1-|\theta|\mathbb EJ
 \ge1-u_1(b_\rho+h\sqrt m)\ge\tfrac12.
\]
For $|u|\ge u_1$, the same even-integer inequalities, now without
conditioning, give $1\pm\mathbb E\cos(uY)\ge\delta_\rho/k_\rho^2$.
Since $|J-Y/(h\sqrt m)|\le1$, for
$u_1h\sqrt m\le|\theta|\le\zeta_\rho$ we obtain
$1+\operatorname{Re}\widehat\mu_J(\theta)\ge\zeta_\rho$.
For $\zeta_\rho\le|\theta|\le\uppi$, the Bernoulli factor gives
$|\widehat\mu_J(\theta)|\le\cos(\zeta_\rho/2)$.
These three ranges prove the stated choice of $\beta_P$.
For the general cutoff, $S=4b_\rho/\delta_\rho$, so the constants
depend only on $b_\rho,\delta_\rho$.
\end{proof}

\subsection{Proof of Proposition~\ref{up:spectral-bound}}
\label{sec:rhmc-spectral-proof}

\begin{proof}[Proof of Proposition~\ref{up:spectral-bound}]
Reversibility follows from Lemma~\ref{up:reversibility}.
Apply the spectral assertion of Lemma~\ref{up:law-comparison} directly
to the full law $\mu_J$. It gives $P\succeq(-1+\beta_P)I$ with
$\beta_P$ satisfying \eqref{up:law-constants}; no energy cutoff is used.
Lemma~\ref{up:rounded-fourier} makes this constant positive uniformly in
admissible $N$.
The spectrum of the reversible Markov contraction $P$ lies in
$[-1+\beta_P,1]$. On this interval,
$1-\lambda^2=(1-\lambda)(1+\lambda)\ge\beta_P(1-\lambda)$.
The spectral theorem gives
\begin{equation}\label{up:actual-squared}
 \|f\|_{L^2(\pi)}^2-\|Pf\|_{L^2(\pi)}^2
 \ge\beta_P\mathcal E_P(f),\qquad f\in L^2(\pi).
\end{equation}
\end{proof}

\begin{lemma}[The three duration laws]\label{up:example-fourier}
Suppose $N\ge16$ and $h=4/(N\sqrt m)$. For each duration law in
Table~\ref{tab:explicit-examples}, the displayed choices of
$\tau_*,c,\beta_P$ satisfy \eqref{up:law-constants}, with the stated
retained probability $p_S$.
\end{lemma}

The proof of Lemma~\ref{up:example-fourier}, including the explicit constants
for the three laws, is given in Appendix~\ref{app:duration-fourier}.

\section{Restricted conductance}\label{sec:entropy}\label{sec:rhmc-conductance}

We prove Proposition~\ref{up:restricted-conductance} in four steps.
First, Lemmas~\ref{ent:modified-energy} and \ref{ent:shell} control the
Verlet energy error along a trajectory, which bounds the rejection
probability. Second, for the triangular duration law \eqref{up:triangular}
with integration time $T=4/\sqrt m$, Proposition~\ref{ent:raw-entropy}
shows that the proposal step contracts relative entropy: the law of the
proposed position $x_J$ is closer to $\pi$ than the initial law, by a
fixed factor. The argument is a discrete-time version of the weighted
entropy estimate of \cite{PP}, carried out in
Section~\ref{sec:rhmc-proposal-entropy}. Third,
Lemma~\ref{up:reference-conductance} converts this contraction into a
restricted conductance bound for the kernel $P_{\mathrm{ref}}$ built on
the triangular law. Fourth, Section~\ref{sec:rhmc-conductance-proof}
transfers the bound to the kernel $P$ of a general non-resonant duration
law, using the comparison Lemma~\ref{up:law-comparison} of
Section~\ref{sec:rhmc-invariance}. The triangular law is used only inside
the proof; the theorem holds for every law satisfying
Definition~\ref{def:nonresonant}.

Normalize $U(x_*)=0$ and write $g=\nabla U$. For $h>0$ and $T=Nh$, apply
the Verlet update \eqref{up:verlet} for $N$ steps. Between kicks the
position moves freely: within step $n$ it is $x_n+\sigma w_n$ for
$0\le\sigma\le h$, where $w_n$ is the half-step velocity. The entropy
estimate of Section~\ref{sec:rhmc-proposal-entropy} is proved along these
free-flight segments and then evaluated at the step endpoints.

\subsection{Energy control}\label{sec:energy-control}

The energy error of a Verlet step splits into the increment of a
modified energy, which telescopes along the trajectory, and a remainder
that is cubic in the step size when the Hessian is Lipschitz. This is
where Hessian regularity enters. The identity and the cubic bound are
\cite[Lemma~2.1]{PP}; the proof is short, and we include it.

\begin{lemma}[Modified energy of velocity Verlet]\label{ent:modified-energy}
Assume \ref{ass:rhmc-a3}. Write $H_n=H(x_n,v_n)$, and
define the modified energy $F_n$ and the remainder $r_n$ by
\[
 F_n=H_n-\frac{h^2}{8}|g(x_n)|^2,
 \qquad
 r_n=U(x_{n+1})-U(x_n)
   -\tfrac12\langle g(x_n)+g(x_{n+1}),x_{n+1}-x_n\rangle.
\]
Then
\begin{equation}\label{ent:exact-energy}
 F_{n+1}-F_n=r_n,
 \qquad |r_n|\le \frac{L_Hh^3}{12}|w_n|^3.
\end{equation}
\end{lemma}

\begin{proof}
Put $a=x_{n+1}-x_n=hw_n$. Substituting $v_n=w_n+(h/2)g(x_n)$ and
$v_{n+1}=w_n-(h/2)g(x_{n+1})$ gives
\[
 \tfrac12(|v_{n+1}|^2-|v_n|^2)
 =-\tfrac12\langle g(x_n)+g(x_{n+1}),a\rangle
  +\tfrac{h^2}{8}(|g(x_{n+1})|^2-|g(x_n)|^2).
\]
Adding the potential change proves $F_{n+1}-F_n=r_n$. Applying the
fundamental theorem of calculus to $t\mapsto\langle g(x_n+ta),a\rangle$
gives
\[
 r_n=\int_0^1(\tfrac12-t)
       a^{\mathsf T}\nabla^2U(x_n+ta)a\,dt.
\]
Subtract $\nabla^2U(x_n+a/2)$ inside the integral.  Its coefficient integrates
to zero, and $\int_0^1(t-1/2)^2dt=1/12$.  This proves the second assertion in
\eqref{ent:exact-energy}.
\end{proof}

The energy-shell bounds that follow are specific to the present setting.
Small $a_T$ controls the accumulated change in $F_n$ and keeps the
trajectory in a bounded energy region; $b_h$ controls the endpoint
correction from $F_n$ to $H_n$.

\begin{lemma}[Energy control]\label{ent:shell}
Assume \ref{ass:rhmc-a1}--\ref{ass:rhmc-a3}. Suppose $R\ge1$,
$h\sqrt L\le1$, and
\begin{equation}\label{ent:errors}
 a_T=L_HT h^2R^{3/2},\qquad b_h=h^2LR,\qquad
 2a_T+\frac12b_h\le R.
\end{equation}
For every input satisfying $H(x_0,v_0)\le R$, all $0\le n\le N$ satisfy
\begin{equation}\label{ent:shell-bounds}
 F_n\le R+2a_T,\qquad H_n\le2R,\qquad |H_n-H_0|\le 2a_T+b_h.
\end{equation}
For $0\le n<N$, also $|w_n|^2\le8R$ and
$\sum_{i<n}|r_i|\le2a_T$.  In particular every possible endpoint of this
trajectory is accepted by the Metropolis test with probability
at least $e^{-(2a_T+b_h)}$.
\end{lemma}

\begin{proof}
The bound $|g(x)|^2\le2LU(x)$ follows by applying the
descent lemma at $x-L^{-1}g(x)$ and using $U\ge0$.  Consequently
\begin{equation}\label{ent:energy-comparison}
 (1-h^2L/4)H_n\le F_n\le H_n,
 \qquad |w_n|^2\le2(1+h^2L/4)H_n\le4F_n.
\end{equation}
Since $h^2L\le1$, the first bound implies $F_n\ge3H_n/4$.
The middle bound follows by Cauchy-Schwarz applied to
$v_n-(h/2)g(x_n)$, with the weights $1$ and $h\sqrt L/2$.
If $F_i\le2R$ for preceding indices, then $|w_i|^2\le8R$ for those
indices, and the cubic bound in \eqref{ent:exact-energy} gives
\[
 \sum_{i<n}|r_i|
 \le \frac23 L_HT h^2(2R)^{3/2}\le2a_T.
\]
By the identity in \eqref{ent:exact-energy}, $F_n\le F_0+2a_T\le R+2a_T\le2R$.
Induction, starting from $F_0\le R$, proves this for every index.
Moreover, \eqref{ent:errors} gives
\[
 H_n\le\frac{R+2a_T}{1-h^2L/4}\le2R,
\]
since the last inequality is equivalent to $2a_T+b_h/2\le R$.
Also $|w_n|^2\le8R$. Finally,
\[
 |H_n-H_0|
 \le 2a_T+\frac{h^2L}{4}(H_n+H_0)
 \le2a_T+\tfrac34 b_h\le2a_T+b_h.
\]
The lower bound on acceptance follows directly from its definition.
\end{proof}

We now apply Lemma~\ref{ent:shell} to the trajectories that arise in
the proof of Proposition~\ref{up:restricted-conductance}, with the
parameters of Table~\ref{tab:proof-parameters}. Both the reference
trajectories and those conditional on $\tau\le\tau_*$ have physical
integration time at most $T_*$. Indeed, $4/\sqrt m\le T_*$ and
\[
 Jh\le\tau+h\le\tau_*+\frac1{4\sqrt m}\le T_*,
\]
since $N\ge16/\sqrt{\delta_\rho}$ gives $h\sqrt m\le1/4$.
Set $a_*=L_HT_*h^2R^{3/2}$ and $b_h=h^2LR$. The step-size hypothesis
in Proposition~\ref{up:restricted-conductance} is exactly
\begin{equation}\label{up:proof-step}
 a_*+b_h=h^2(LR+L_HT_*R^{3/2})\le\eta/3.
\end{equation}
Since $R\ge1$ and $\eta=c/600\le1/600$, this implies
$h^2L=b_h/R<1$ and $2a_*+b_h/2\le2\eta/3<R$.
Thus the hypotheses of Lemma~\ref{ent:shell} hold.

Let $N_*=\lfloor T_*/h\rfloor$. Apply Lemma~\ref{ent:shell} with physical integration time
$N_*h$. Its accumulated error parameter $a_{N_*h}$ is at most $a_*$,
and its endpoint error parameter is $b_h$. Thus for $H(x_0,v_0)\le R$,
every endpoint $(Q,W)$ reached by time $T_*$ has an energy increase
at most $b_h+2a_*\le\eta$, so
\begin{equation}\label{up:shell-work}
 \alpha(X,V;Q,W)\ge e^{-\eta}\ge1-\eta.
\end{equation}

\subsection{Entropy of the Verlet proposal}\label{sec:rhmc-proposal-entropy}

For a joint law $p$ of position and velocity with position marginal $\nu$, define
\begin{equation}\label{pre:entropy-chain}
 E_{\mathrm{joint}}(p)=\KL(p\Vert\Pi),\qquad
 E_{\mathrm{pos}}(\nu)=\KL(\nu\Vert\pi),\qquad
 E_v(p)=\KL(p\Vert\nu\otimes\gamma_d).
\end{equation}
The entropy chain rule gives $E_{\mathrm{joint}}(p)=E_{\mathrm{pos}}(\nu)+E_v(p)$.
When $\nu$ and $\pi$ have finite differential entropy and finite expectations of $U$,
write $\nu(x)$ and $\pi(x)$ for their Lebesgue densities and define the
difference of negative differential entropies by
\begin{equation}\label{pre:deltaL}
 \Delta\mathcal L(\nu)=\int_{\R^d}\nu(x)\log\nu(x)\,\mathrm dx
       -\int_{\R^d}\pi(x)\log\pi(x)\,\mathrm dx
 =E_{\mathrm{pos}}(\nu)-\bigl(\E_\nu U-\E_\pi U\bigr).
\end{equation}

We first recall Talagrand's transport-entropy inequality \cite{Villani}.
In the weighted entropy argument, Talagrand's inequality bounds the initial
squared Wasserstein distance to the target by the initial relative entropy.

\begin{lemma}[Talagrand's transport-entropy inequality]\label{lem:talagrand}
Assume \ref{ass:rhmc-a1} and let $\pi$ be given by \eqref{pre:target}.
For every probability measure $\nu$ on $\R^d$ with
$\KL(\nu\Vert\pi)<\infty$,
\begin{equation}\label{pre:smooth-bound}
 \frac12W_2(\nu,\pi)^2\le\frac1m\KL(\nu\Vert\pi).
\end{equation}
\end{lemma}

The free-transport entropy inequality of \cite{PP} bounds the velocity
relative entropy during a free flight by the sum of its endpoint values.
Together with the Verlet energy identities, it allows us to bound $\KL(\Law(x_J)\Vert\pi)$, where $x_J$ is the position
proposed for the Metropolis test and $J$ has the rounded triangular law.

\begin{lemma}[{\cite[Proposition~3.1 and Corollary~3.5]{PP}}]\label{lem:interpolation}
For the step size $h>0$, let $X,V$ be random vectors in $\R^d$, write
\[
 p_t=\Law(X+tV,V),\qquad \nu_t=\Law(X+tV),
\]
and define $e(t)=\KL(p_t\,\Vert\,\nu_t\otimes\gamma_d)$ for $0\le t\le h$.
If $e(0)+e(h)<\infty$, then
\begin{equation}\label{pre:endpoint}
 e(t)\le e(0)+e(h),\qquad 0\le t\le h.
\end{equation}
In particular all intermediate entropies are finite.
\end{lemma}
This is \cite[Corollary~3.5]{PP}.

The free-transport estimate below bounds the second derivative of
$\frac12W_2(\nu_t,\pi)^2$ by the first and second moments of the
conditional velocity law. The following identity converts those
moments into relative entropy with respect to $\gamma_d$; for $G=\Sigma^{-1}$ it
reduces to the maximum-entropy property of the Gaussian. 

\begin{lemma}[Relative entropy to a Gaussian]\label{ent:gaussian-inequality}
Let $\mu$ be a probability measure on $\R^d$ with finite second moment,
mean $u$ and covariance $\Sigma$, and let $G$ be symmetric positive
definite. Write $\Sigma:G=\operatorname{tr}(\Sigma G)$ for the Frobenius
inner product. If $\KL(\mu\|\gamma_d)<\infty$, then
\[
 2\KL(\mu\|\gamma_d)
 =|u|^2+\operatorname{tr}\Sigma-\Sigma:G+\log\det G
    +2\KL\bigl(\mu\|N(u,G^{-1})\bigr).
\]
In all cases,
\[
 |u|^2+\operatorname{tr}\Sigma-\Sigma:G+\log\det G
 \le 2\KL(\mu\|\gamma_d),
\]
with equality if and only if $\mu=N(u,G^{-1})$.
\end{lemma}

\begin{proof}
If $\KL(\mu\|\gamma_d)=\infty$, the inequality is immediate and
strict. Otherwise, put $\widetilde\gamma=N(u,G^{-1})$. The logarithm
of the density ratio of $\widetilde\gamma$ and $\gamma_d$ is quadratic,
hence integrable under $\mu$. Thus $\KL(\mu\|\widetilde\gamma)$ is
finite, and
\[
 2\KL(\mu\|\gamma_d)
 =2\KL(\mu\|\widetilde\gamma)
   +\int_{\R^d}\bigl[|v|^2-(v-u)^{\mathsf T}G(v-u)+\log\det G\bigr]\,\mu(dv).
\]
Evaluating the integral gives the identity. Nonnegativity of
$\KL(\mu\|\widetilde\gamma)$ gives the inequality and the equality case.
\end{proof}

Let $p=p_0$ be the joint law of position and velocity $(X,V)$, with
$p_t$ and $\nu_t$ as in Lemma~\ref{lem:interpolation}. For a position
law $\nu$ having a density, let $T_\nu$ be its Brenier map to $\pi$. The
current of $p$ is
\[
 \mathcal J(p)=\int_{\R^{2d}}\langle v,x-T_\nu(x)\rangle\,p(dx,dv).
\]
Use $E_v$ and $\Delta\mathcal L$ from \eqref{pre:entropy-chain}
and \eqref{pre:deltaL}.

\begin{lemma}[Free-transport estimate]\label{ent:free-flight}
Assume \ref{ass:rhmc-a1} and \ref{ass:rhmc-a2}. Let $p$ be compactly supported and satisfy
$\KL(p\|\Pi)<\infty$. For
$0\le a\le b<\infty$,
\begin{equation}\label{ent:current-bound}
 \mathcal J(p_b)-\mathcal J(p_a)
 \le \int_a^b\{2E_v(p_t)-\Delta\mathcal L(\nu_t)\}\,dt.
\end{equation}
Moreover $\mathcal W(t)=\frac12W_2(\nu_t,\pi)^2$ is continuously differentiable and
$\mathcal W'(t)=\mathcal J(p_t)$.  Thus, in the distributional sense,
\begin{equation}\label{ent:free-second}
 \mathcal W''(t)\le 2E_{\mathrm{joint}}(t)-3E_{\mathrm{pos}}(t)+U_\Delta(t),
\end{equation}
where
\[
 E_{\mathrm{joint}}(t)=\KL(p_t\|\Pi),\qquad
 E_{\mathrm{pos}}(t)=\KL(\nu_t\|\pi),\qquad
 U_\Delta(t)=\mathbb E U(X+tV)-\int_{\R^d} U\,d\pi.
\]
\end{lemma}

\begin{proof}
\widowpenalty=10000
We first identify the derivative of $\mathcal W$ with the current
$\mathcal J(p_t)$. For the position-smoothed laws, we then bound its second
differences using the conditional velocity covariance and the Gaussian
entropy inequality. Integrating this bound and removing the smoothing
gives \eqref{ent:current-bound}. Assumptions \ref{ass:rhmc-a1} and
\ref{ass:rhmc-a2} enter only through $\pi$: by \ref{ass:rhmc-a1} it has
a positive density with finite second moment, so optimal couplings to
$\pi$ and Brenier maps exist, and by \ref{ass:rhmc-a2} the potential
grows at most quadratically, so $\int_{\R^d}U\,d\pi$ and
$\Delta\mathcal L(\nu_t)$ are finite. Both assumptions are also used in
Lemma~\ref{app:smoothed-regularity}, on which Steps 2 to 5 rely.

\emph{Step 1: first variation.}
The identity $\mathcal W'=\mathcal J(p_t)$ is the first-variation
formula for $W_2$ along the curve $t\mapsto\nu_t$
\cite[Theorem~8.4.7]{AGSflows}; we give the short argument because we
need the derivative at every $t$, together with its continuity.
Finite relative entropy gives a density for $(X,V)$. The free shear
preserves absolute continuity, so each position law $\nu_t$ has a density.
Let $M_v=\mathbb E|V|^2$. Couple $X+tV$ optimally to $\pi$ and keep its
target coordinate fixed while moving the source by $\epsilon V$. The
resulting transport cost bounds $\mathcal W(t+\epsilon)$ from above.
Reversing the roles of $t$ and $t+\epsilon$ gives the lower bound
\[
 \epsilon \mathcal J(p_{t+\epsilon})-\tfrac12\epsilon^2M_v
 \le \mathcal W(t+\epsilon)-\mathcal W(t)
 \le \epsilon \mathcal J(p_t)+\tfrac12\epsilon^2M_v,\qquad \epsilon>0.
\]
We claim that $\mathcal J(p_s)\to\mathcal J(p_t)$ as $s\to t$. The laws
$p_s$ converge to $p_t$ in $W_2$. Lift the optimal coupling of $\nu_s$
with $\pi$ to a coupling of $p_s$ with $\pi$ by carrying the velocity
along. These lifted couplings are tight, and every limit point is an
optimal coupling of $\nu_t$ with $\pi$ by \cite[Theorem~5.20]{Villani}.
Since $\nu_t$ has a density, that optimal coupling is unique and is
induced by the Brenier map $T_t$. Hence the lifted couplings converge weakly to the law of
$(X_t,V_t,T_t(X_t))$, and $\mathcal J(p_s)$ is the integral of
$\langle v,x-y\rangle$ against the lifted coupling at time $s$. This
integrand is bounded by $|v|^2+|x|^2+|y|^2$, whose integrals converge, so
it is uniformly integrable and the integrals converge to
$\mathcal J(p_t)$.
Applying this to $t+\epsilon\to t$ in the two-sided estimate gives
$\mathcal W'(t)=\mathcal J(p_t)$ and continuity of the derivative.

\emph{Step 2: second variation.}
Fix a smoothing parameter $\iota\in(0,1]$, to be sent to zero in
Step~5. Add independent Gaussian noise of covariance $\iota I_d$ to $X$
and leave $V$ fixed. Write $p^\iota$ for the resulting joint law of
position and velocity and suppress the superscript in this calculation.
The coupling of Step~1 gives the first derivative, but its
second-difference bound is only $M_v$. To obtain the entropy bound,
we represent the same position laws using the flow of the conditional
mean velocity. Let $n_t$ be the position density and write
\begin{equation}\label{ent:conditional-mixture}
 u_t(y)=\mathbb E[V\mid X+tV=y],\qquad
 \Sigma_t(y)=\operatorname{Cov}(V\mid X+tV=y),
\end{equation}
with $X$ now denoting the smoothed position. The smoothing makes $n_t$
positive and the spatial moment fields smooth. On compact time
intervals, $u_t,\Sigma_t$ and their first spatial derivatives are
bounded, and $\nabla\log n_t$ has at most linear growth, by
Lemma~\ref{app:smoothed-regularity}\ref{app:regularity-moments}
and~\ref{app:regularity-score}.

The zeroth and first velocity moments of the free-transport equation
give the continuity and momentum equations
\[
 \partial_t n_t+\operatorname{div}(n_tu_t)=0,\qquad
 \partial_t(n_tu_t)+\operatorname{div}\bigl(n_t(u_t\otimes u_t+\Sigma_t)\bigr)=0.
\]
Expand the second equation and use the first. The acceleration of the
flow of $u_t$ is therefore
\begin{equation}\label{ent:fluid}
 a_t=\partial_tu_t+(Du_t)u_t
 =-\frac{\operatorname{div}(n_t\Sigma_t)}{n_t}
 =-\operatorname{div}\Sigma_t-\Sigma_t\nabla\log n_t.
\end{equation}
In particular $a_t$ is continuous with at most linear growth in space.

Define the flow $\Phi_{t,\sigma}$ by
\[
 \partial_\sigma\Phi_{t,\sigma}(y)
 =u_{t+\sigma}\bigl(\Phi_{t,\sigma}(y)\bigr),\qquad
 \Phi_{t,0}(y)=y.
\]
The continuity equation gives
$n_{t+\sigma}(\Phi_{t,\sigma}(y))\det D_y\Phi_{t,\sigma}(y)=n_t(y)$, by differentiating
the left side in $\sigma$. Thus $\Phi_{t,\sigma}$ transports $\nu_t$ to
$\nu_{t+\sigma}$. Its acceleration is $a_{t+\sigma}(\Phi_{t,\sigma}(y))$, so
\[
 \Phi_{t,\sigma}(y)=y+\sigma u_t(y)+\tfrac12\sigma^2a_t(y)+o_{L^2(n_t\,dy)}(\sigma^2).
\]
For the remainder, integrate the acceleration twice. The difference
$a_{t+r}(\Phi_{t,r}(y))-a_t(y)$ tends to zero pointwise as $r\to0$ and
is bounded by $C(1+|y|)$ for small $r$, since $u_t$ is bounded and
$a_t$ has linear growth. Hence the difference tends to zero in
$L^2(n_t\,dy)$.

Set $T_t=T_{\nu_t}$. For $Y\sim\nu_t$, the pair
$(\Phi_{t,\sigma}(Y),T_t(Y))$ couples $\nu_{t+\sigma}$ to $\pi$.
Its cost therefore bounds $\mathcal W(t+\sigma)$ from above.
For $\sigma=\pm\epsilon$, expanding this cost gives
\begin{align*}
 \mathcal W(t+\sigma)
 &\le \tfrac12\int_{\R^d}|\Phi_{t,\sigma}(y)-T_t(y)|^2 n_t(y)\,dy\\
 &=\mathcal W(t)+\sigma\int_{\R^d}\langle y-T_t(y),u_t(y)\rangle n_t(y)\,dy\\
 &\quad+\tfrac12\sigma^2\int_{\R^d}
   \bigl(|u_t(y)|^2+\langle y-T_t(y),a_t(y)\rangle\bigr)n_t(y)\,dy+o(\sigma^2).
\end{align*}
The error in the squared norm is $o(\sigma^2)$ by Cauchy-Schwarz.
Write $D^2_\epsilon f(t)=[f(t+\epsilon)+f(t-\epsilon)-2f(t)]/\epsilon^2$.
Adding the two estimates and using \eqref{ent:fluid} gives
\begin{equation}\label{ent:frozen-second}
 \limsup_{\epsilon\to0}D^2_\epsilon\mathcal W(t)
 \le \int_{\R^d} |u_t|^2 n_t\,dy
 -\int_{\R^d}\langle y-T_t(y),\operatorname{div}(n_t\Sigma_t)(y)\rangle\,dy .
\end{equation}

Write $T_t=\nabla\varphi_t$, where $\varphi_t$ is convex. The Brenier
potential need not be twice continuously differentiable, so we
integrate by parts using its distributional Hessian. Its singular
part contributes a nonpositive term. Decompose
$D^2\varphi_t=G_t\,dy+D^s\varphi_t$.
Both parts are positive semidefinite, and $G_t$ is the approximate
derivative of $T_t$; see the discussion following
\cite[Theorem~11.1]{Villani}. To integrate the divergence
term by parts, let $0\le\chi_k\uparrow1$ be smooth cutoffs, equal to
one on the ball of radius $k$, with $|\nabla\chi_k|\le C/k$. Then
\begin{align*}
 &-\int_{\R^d}\chi_k\langle y-T_t(y),\operatorname{div}(n_t\Sigma_t)(y)\rangle\,dy\\
 &\quad=\int_{\R^d}\langle y-T_t(y),n_t\Sigma_t\nabla\chi_k\rangle\,dy
       +\int_{\R^d}\chi_k n_t\Sigma_t:(I-G_t)\,dy
       -\int_{\R^d}\chi_k n_t\Sigma_t:dD^s\varphi_t.
\end{align*}
The last term is nonpositive because both $\Sigma_t$ and the singular
Hessian measure are positive semidefinite. The first term tends to
zero since $\Sigma_t$ is bounded and $y-T_t(y)\in L^2(n_t\,dy)$.
The left side converges since $a_t$ has linear growth. Dropping the
nonpositive term also bounds the increasing integrals
$\int_{\R^d}\chi_k n_t\Sigma_t:G_t$, so monotone convergence gives
\begin{equation}\label{ent:pressure}
 -\int_{\R^d}\langle y-T_t(y),\operatorname{div}(n_t\Sigma_t)(y)\rangle\,dy
 \le \int_{\R^d} n_t\Sigma_t:(I-G_t)\,dy.
\end{equation}

\emph{Step 3: entropy bound.}
Both position laws have densities and finite second moments, so
the Brenier change-of-variables formula
\cite[Theorem~11.1 and Example~11.2]{Villani} gives
$n_t(y)=\pi(T_t(y))\det G_t(y)$ almost everywhere. Hence $G_t(y)$
is positive definite almost everywhere. For such $y$, the conditional
law $\mu_y=\mathcal L(V\mid X+tV=y)$ has mean $u_t(y)$ and covariance
$\Sigma_t(y)$, and its relative entropy with respect to $\gamma_d$ is
finite for almost every $y$ by the chain rule and
Lemma~\ref{app:smoothed-regularity}\ref{app:regularity-entropies}.
Lemma~\ref{ent:gaussian-inequality} with $\mu=\mu_y$ and $G=G_t(y)$ gives
\[
 |u_t(y)|^2+\operatorname{tr}\Sigma_t(y)-\Sigma_t(y):G_t(y)+\log\det G_t(y)
 \le2\KL(\mu_y\Vert\gamma_d).
\]
Integrate against $n_t(y)\,dy$. By the chain rule the right side
integrates to $2E_v(p_t)$, and for the last term on the left we use
\[
 \int_{\R^d} n_t\log\det G_t\,dy
 =\int_{\R^d} n_t\log n_t\,dy-\int_{\R^d}\pi\log\pi\,dy
 =\Delta\mathcal L(\nu_t).
\]
The Gaussian bounds in
Lemma~\ref{app:smoothed-regularity}\ref{app:regularity-gaussian}
give integrability of $\log n_t$,
and $\log\pi(T_t)$ is integrable since $T_t$ transports $n_t$ to $\pi$.
Thus the integrated inequality is
\begin{equation}\label{ent:entropy-second}
 \int_{\R^d} n_t\bigl(|u_t|^2+\operatorname{tr}\Sigma_t-\Sigma_t:G_t\bigr)\,dy
 \le 2E_v(p_t)-\Delta\mathcal L(\nu_t).
\end{equation}

Combining \eqref{ent:frozen-second}, \eqref{ent:pressure} and
\eqref{ent:entropy-second} gives
\begin{align*}
 \limsup_{\epsilon\to0}D^2_\epsilon\mathcal W(t)
 &\le\int_{\R^d} n_t\bigl(|u_t|^2+\operatorname{tr}\Sigma_t-\Sigma_t:G_t\bigr)\,dy\\
 &\le 2E_v(p_t)-\Delta\mathcal L(\nu_t).
\end{align*}
\emph{Step 4: integration in time.}
We now pass from this bound on second differences to the integrated
inequality \eqref{ent:current-bound}.
The upper bound in the two-sided estimate of Step 1, applied at the
same $t$ with both signs of $\epsilon$, gives
$D^2_\epsilon\mathcal W(t)\le M_v$. For a smooth nonnegative test
function $\zeta$ compactly supported in $(0,\infty)$, move
$D^2_\epsilon$ from $\mathcal W$ to $\zeta$ under the integral.
Reverse Fatou, using the upper bound $M_v$, then gives
\[
 \langle\mathcal W'',\zeta\rangle
 =\lim_{\epsilon\to0}\int_{(0,\infty)}\zeta D^2_\epsilon\mathcal W\,dt
 \le\int_{(0,\infty)}\zeta\,[2E_v(p_t)-\Delta\mathcal L(\nu_t)]\,dt.
\]
The right side is locally integrable by
Lemma~\ref{app:smoothed-regularity}\ref{app:regularity-entropies}.
Since $\mathcal W'=\mathcal J(p_t)$ is continuous, the function
\[
 t\longmapsto \mathcal J(p_t)
 -\int_a^t\bigl[2E_v(p_r)-\Delta\mathcal L(\nu_r)\bigr]\,dr
\]
is continuous and has nonpositive distributional derivative, hence
is nonincreasing. Comparing its values at $a$ and $b$ gives
\eqref{ent:current-bound} for the smoothed law.

\emph{Step 5: removing the smoothing.}
We pass to the limit using the convergence claims in
Appendix~\ref{app:regularization}:
\ref{app:convergence-entropies} and \ref{app:convergence-energies}
give convergence of the integral, and \ref{app:convergence-currents}
gives convergence of the endpoint currents. This proves
\eqref{ent:current-bound} for the original law. Finally, the chain
rule and \eqref{pre:deltaL} give
\[
 2E_v(p_t)-\Delta\mathcal L(\nu_t)
 =2\bigl(E_{\mathrm{joint}}(t)-E_{\mathrm{pos}}(t)\bigr)
  -\bigl(E_{\mathrm{pos}}(t)-U_\Delta(t)\bigr).
\]
Together with $\mathcal W'=\mathcal J$, this yields \eqref{ent:free-second}.
\end{proof}

For the weighted estimate, assume \ref{ass:rhmc-a1} and \ref{ass:rhmc-a2}.
Let $T=Nh$ and let $p_0$ be the initial joint law of position and velocity. Assume
it has finite relative entropy $E_0=\KL(p_0\|\Pi)$, is supported on
$\{H\le R\}$, and satisfies
\begin{equation}\label{ent:centered}
 \mathbb E[v_0\mid x_0]=0.
\end{equation}
Condition \eqref{ent:centered} holds when $\Pi$ is conditioned on
$\{x_0\in A,H\le R\}$.

For $t=nh+\sigma$, $0\le \sigma\le h$, define
\begin{align*}
 E_{\mathrm{pos}}(n,\sigma)&=\KL(\mathcal L(x_n+\sigma w_n)\|\pi),&
 E_{\mathrm{joint}}(n,\sigma)&=\KL(\mathcal L(x_n+\sigma w_n,w_n)\|\Pi),\\
 E_n^x&=\KL(\mathcal L(x_n)\|\pi),&
 U_\Delta(t)&=\mathbb E U(x_n+\sigma w_n)-\mathbb E_\pi U,\\
 \mathcal W(t)&=\tfrac12W_2(\mathcal L(x_n+\sigma w_n),\pi)^2.&
\end{align*}
Here $E_0$ denotes the initial joint entropy, whereas $E_0^x$ denotes
the initial position entropy; in particular $E_0^x\le E_0$.
Write $E_{\mathrm{pos}}(t)=E_{\mathrm{pos}}(n,\sigma)$ and
$E_{\mathrm{joint}}(t)=E_{\mathrm{joint}}(n,\sigma)$ for $nh<t<(n+1)h$, and set
\[
 I_{\mathrm{pos}}=\int_0^T(T-t)E_{\mathrm{pos}}(t)\,dt,\qquad
 I_{\mathrm{joint}}=\int_0^T(T-t)E_{\mathrm{joint}}(t)\,dt.
\]
Define
\begin{align*}
 \omega_n&=\int_0^h(T-nh-\sigma)\,d\sigma=h(T-nh)-\tfrac12h^2,\\
 b_n&=\mathbb E[F_n-H(x_0,v_0)],\\
 Q_n&=\mathbb E[U(x_{n+1})-U(x_n)-h\langle g(x_n),w_n\rangle]\ge0,\\
 I_{\mathrm{end}}&=\sum_n\omega_n(E_n^x+E_{n+1}^x),\qquad
 I_{\mathrm{mod}}=\sum_n\omega_nb_n,\\
 Q_\Sigma&=\sum_n\omega_nQ_n,\qquad
 \Delta_U=\mathbb EU(x_N)-\mathbb EU(x_0),\\
 Q_{\mathrm{int}}
 &=\sum_n\int_0^h(T-nh-\sigma)
   \mathbb E[U(x_n+\sigma w_n)-U(x_n)-\sigma\langle g(x_n),w_n\rangle]\,d\sigma\ge0.
\end{align*}
All sums over $n$ run from $0$ to $N-1$.
\begin{lemma}[Weighted entropy estimate]\label{ent:weighted-entropy}
Assume \ref{ass:rhmc-a1} and \ref{ass:rhmc-a2}, and let the initial joint
law $p_0$ have finite relative entropy $E_0$, be supported on
$\{H\le R\}$, and satisfy \eqref{ent:centered}. Then
\begin{equation}\label{ent:signed-final}
 3I_{\mathrm{end}}\le \mathcal W(0)+\tfrac52E_0T^2+5I_{\mathrm{mod}}+3Q_\Sigma
     +\frac{h^2}{4}\Delta_U
     -\frac{h^2}{12}\sum_nQ_n-Q_{\mathrm{int}}.
\end{equation}
\end{lemma}

\begin{proof}
\emph{Step 1: free transport on each drift.}
The interpolated position is continuous at each time \(nh\). At the
start of each drift, the joint law has compact support and finite
relative entropy with respect to \(\Pi\). The preceding maps are
continuous and preserve Lebesgue volume, and \(H\) is bounded on
the image of the initial support.
Thus Lemma~\ref{ent:free-flight}, which uses \ref{ass:rhmc-a1} and
\ref{ass:rhmc-a2}, applies to each drift.

\emph{Step 2: the kicks, and integration against the weight $T-t$.}
At time $nh$ for $1\le n<N$, the interpolation velocity changes by
$-hg(x_n)$.
Let $T_n$ be the Brenier map from $\mathcal L(x_n)$ to $\pi$, and put
\[
 A_n=\mathbb E\langle g(x_n),x_n-T_n(x_n)\rangle
     =U_\Delta(nh)+C_n^U,
\]
where
\[
 C_n^U=\mathbb E\big[U(T_n(x_n))-U(x_n)
                    -\langle g(x_n),T_n(x_n)-x_n\rangle\big]\ge0.
\]
The equality uses $T_n(x_n)\sim\pi$ and the inequality is the convexity
of $U$ in \ref{ass:rhmc-a1}.
The one-sided first variations therefore satisfy
\[
 \mathcal W'(nh+)-\mathcal W'(nh-)=-hA_n\quad(1\le n<N),\qquad
 \mathcal W'(0+)=-\tfrac h2A_0.
\]
The last identity uses \eqref{ent:centered} and the initial half kick.
Integrate \eqref{ent:free-second} with weight $T-t$, adding each jump
of $\mathcal W'$ with the same weight. Because the weight is linear and
vanishes at $T$, integration by parts gives only the endpoint terms
$\mathcal W(T)-\mathcal W(0)-T\mathcal W'(0+)$ on the left. For a general smooth weight $a(t)$,
the same calculation also gives an interior term
$\int_0^T a''(t)\mathcal W(t)\,\mathrm dt$. Substituting the initial half-kick
and the internal jumps gives
\begin{align}\label{ent:integrated-flight}
 \mathcal W(T)+3I_{\mathrm{pos}}
 &\le \mathcal W(0)+2I_{\mathrm{joint}}+\mathcal Q
 -h\left[\frac T2 C_0^U+
           \sum_{n=1}^{N-1}(T-nh)C_n^U\right],\\
 \mathcal Q
 &=\int_0^T(T-t)U_\Delta(t)\,dt
   -h\left[\frac T2U_\Delta(0)+\sum_{n=1}^{N-1}(T-nh)U_\Delta(nh)\right].\nonumber
\end{align}

\emph{Step 3: the quadrature error $\mathcal Q$.}
Define
\[
 \operatorname{Chord}_n(\sigma)
  =(1-\sigma/h)\mathbb EU(x_n)+(\sigma/h)\mathbb EU(x_{n+1})
                   -\mathbb EU(x_n+\sigma w_n)\ge0,
\]
where the inequality is again the convexity of $U$.
On the interval $[nh,(n+1)h]$, the two endpoint coefficients of the
linear interpolant are $h(T-nh)/2-h^2/6$ and $h(T-nh)/2-h^2/3$.
The contributions from adjacent intervals add to $h(T-nh)$ at each
interior time $nh$, leaving the boundary correction and the difference
between the linear interpolant and the potential.
\begin{equation}\label{ent:chord}
 \mathcal Q=\frac{h^2}{6}\,[U_\Delta(T)-U_\Delta(0)]
       -\sum_{n=0}^{N-1}\int_0^h(T-nh-\sigma)
                         \operatorname{Chord}_n(\sigma)\,d\sigma.
\end{equation}
\emph{Step 4: energy bookkeeping along a Verlet step.}
By volume preservation, the change in joint relative entropy with
respect to $\Pi$ equals the expected change in $H$. Since
$v_n=w_n+(h/2)g(x_n)$, expansion gives
\begin{align}\label{ent:cell-energy}
 E_{\mathrm{joint}}(n,\sigma)
 &=E_0+b_n+(\sigma-h/2)\mathbb E\langle g(x_n),w_n\rangle\\
 &\quad+\mathbb E[U(x_n+\sigma w_n)-U(x_n)-\sigma\langle g(x_n),w_n\rangle].\nonumber
\end{align}
The linear terms at the two drift endpoints cancel when added,
as shown in Figure~\ref{fig:verlet-cancellation}. We obtain
\begin{align}
 E_{\mathrm{joint}}(n,0)+E_{\mathrm{joint}}(n,h)&=2E_0+2b_n+Q_n,\label{ent:endpoint-energy}\\
 I_{\mathrm{joint}}-\tfrac12E_0T^2
 &=I_{\mathrm{mod}}-\frac{h^2}{12}\Delta_U
          +\frac{h^2}{12}\sum_nQ_n+Q_{\mathrm{int}}.\label{ent:weighted-energy}
\end{align}
For the latter identity, the linear term integrates to
\[
 \int_0^h(T-nh-\sigma)(\sigma-h/2)\,d\sigma=-h^3/12,
 \qquad
 \sum_nh\mathbb E\langle g(x_n),w_n\rangle=\Delta_U-\sum_nQ_n.
\]
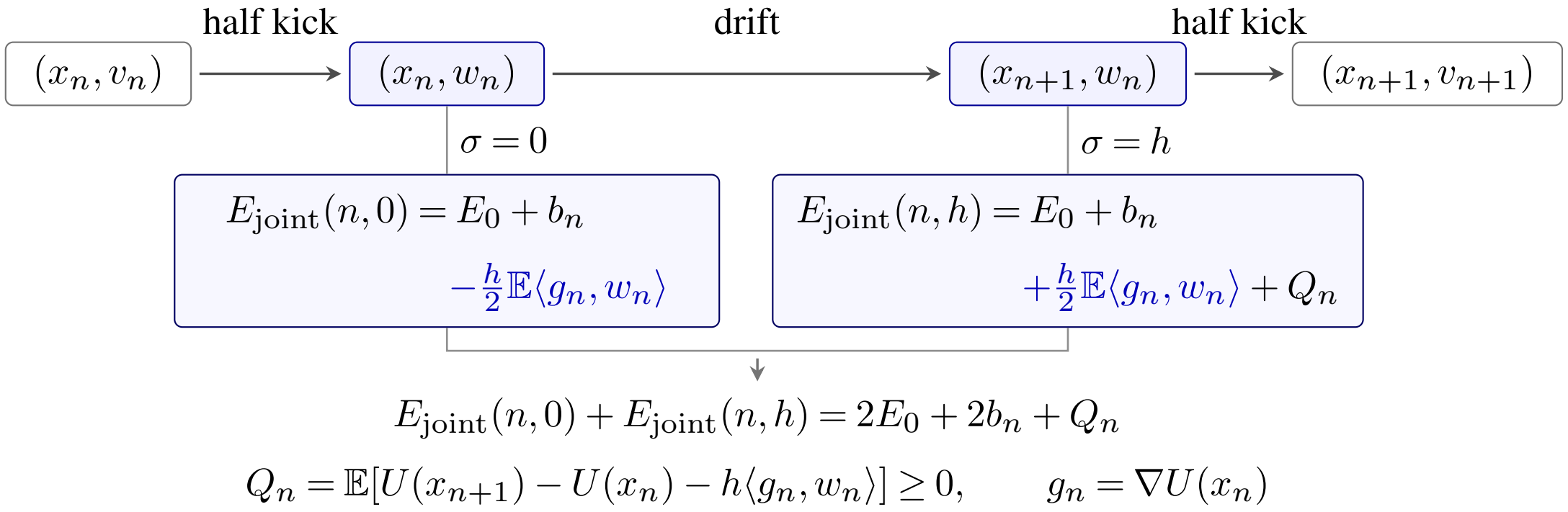
\begin{figure}[!t]
\centering
\begin{tikzpicture}[
  font=\small,
  state/.style={draw=black!55,rounded corners=2pt,
                minimum height=0.58cm,inner xsep=7pt},
  drift state/.style={state,draw=blue!58!black,fill=blue!5},
  cell arrow/.style={-{Stealth[length=4pt,width=4pt]},draw=black!70,
                     line width=0.65pt,shorten <=2pt,shorten >=2pt},
  endpoint/.style={draw=blue!38!black,fill=blue!3,rounded corners=2pt,
                   minimum width=5.0cm,inner xsep=6pt,inner ysep=4pt,
                   align=center,anchor=north},
  guide/.style={draw=black!44,line width=0.5pt}
]
  \node[state] (initial) at (0.8,0) {$(x_n,v_n)$};
  \node[drift state] (start) at (4.0,0) {$(x_n,w_n)$};
  \node[drift state] (end) at (9.7,0) {$(x_{n+1},w_n)$};
  \node[state] (final) at (13.0,0) {$(x_{n+1},v_{n+1})$};

  \draw[cell arrow] (initial.east) -- node[above=7pt] {half kick} (start.west);
  \draw[cell arrow] (start.east) -- node[above=7pt] {drift} (end.west);
  \draw[cell arrow] (end.east) -- node[above=7pt] {half kick} (final.west);

  \node[endpoint] (left) at (4.0,-0.92)
    {$\begin{aligned}
       E_{\mathrm{joint}}(n,0)&=E_0+b_n\\
       &\quad{\color{blue!72!black}-\tfrac h2\mathbb E\langle g_n,w_n\rangle}
     \end{aligned}$};
  \node[endpoint] (right) at (9.7,-0.92)
    {$\begin{aligned}
       E_{\mathrm{joint}}(n,h)&=E_0+b_n\\
       &\quad{\color{blue!72!black}+\tfrac h2\mathbb E\langle g_n,w_n\rangle}+Q_n
     \end{aligned}$};
  \draw[guide] (start.south) -- node[right,inner sep=3pt] {$\sigma=0$} (left.north);
  \draw[guide] (end.south) -- node[right,inner sep=3pt] {$\sigma=h$} (right.north);

  \coordinate (under left) at ([yshift=-6pt]left.south);
  \coordinate (under right) at (right.south |- under left);
  \coordinate (middle) at (6.85,0);
  \coordinate (join) at (middle |- under left);
  \draw[guide] (left.south) -- (under left) -- (under right) -- (right.south);
  \node[anchor=north,inner sep=2pt] (sum) at ([yshift=-10pt]join)
    {$E_{\mathrm{joint}}(n,0)+E_{\mathrm{joint}}(n,h)=2E_0+2b_n+Q_n$};
  \draw[cell arrow,draw=black!54] (join) -- (sum.north);
  \node[anchor=north,inner sep=2pt] at ([yshift=-3pt]sum.south)
    {$Q_n=\mathbb E[U(x_{n+1})-U(x_n)-h\langle g_n,w_n\rangle]\ge0,
     \qquad g_n=\nabla U(x_n)$};
\end{tikzpicture}

\caption{One velocity Verlet step. The linear terms in the two drift
endpoint entropies have opposite signs. Their sum leaves the common
contribution $2(E_0+b_n)$ and the potential remainder $Q_n$.}
\label{fig:verlet-cancellation}
\end{figure}

\emph{Step 5: transfer to the position marginal and conclusion.}
On each interval $[nh,(n+1)h]$, Lemma~\ref{lem:interpolation} gives
\[
 E_{\mathrm{joint}}(n,\sigma)-E_{\mathrm{pos}}(n,\sigma)
 \le E_{\mathrm{joint}}(n,0)-E_n^x+E_{\mathrm{joint}}(n,h)-E_{n+1}^x.
\]
Multiplying by $T-nh-\sigma$ and integrating, we obtain
\begin{equation}\label{ent:node-transfer}
 I_{\mathrm{end}}\le I_{\mathrm{pos}}-I_{\mathrm{joint}}
             +\sum_n\omega_n[E_{\mathrm{joint}}(n,0)+E_{\mathrm{joint}}(n,h)].
\end{equation}
Multiply \eqref{ent:node-transfer} by three and substitute
\eqref{ent:endpoint-energy} and \eqref{ent:integrated-flight}.
Using $\mathcal W(T)\ge0$ and dropping the nonpositive terms involving
$C_n^U$ gives
\[
 3I_{\mathrm{end}}\le \mathcal W(0)-I_{\mathrm{joint}}+\mathcal Q+
 3E_0T^2+6I_{\mathrm{mod}}+3Q_\Sigma.
\]
Substituting \eqref{ent:weighted-energy} and \eqref{ent:chord} and
dropping the nonpositive term involving $\operatorname{Chord}_n$ gives
\eqref{ent:signed-final}.

\end{proof}

The following proposition collects the estimates obtained by adapting
the weighted-entropy argument in the proof of \cite[Lemma~3]{MSSW}
to velocity Verlet.

\begin{proposition}[Relative entropy of the proposal law]\label{ent:raw-entropy}
Assume \ref{ass:rhmc-a1}--\ref{ass:rhmc-a3}. Suppose $T=Nh\ge2/\sqrt m$,
$h\sqrt L\le1$, $R\ge1$, and that $a_T=L_HTh^2R^{3/2}$ and $b_h=h^2LR$
satisfy $2a_T+b_h/2\le R$, as in \eqref{ent:errors}. Let the initial position and velocity $(x_0,v_0)$ have a law supported
on $\{H\le R\}$, with finite relative entropy $E_0$ with respect to
$\Pi$ and with $\mathbb E[v_0\mid x_0]=0$ as in \eqref{ent:centered}.
If $J$ is drawn independently from
\eqref{up:triangular}, the proposed position $x_J$ satisfies
\begin{equation}\label{ent:raw-general}
 \KL(\mathcal L(x_J)\|\pi)
 \le\left(\frac56+\frac{1}{3mT^2}\right)E_0
       +\frac53 a_T+\frac{49}{24}b_h.
\end{equation}
\end{proposition}

\begin{proof}
Lemma~\ref{ent:shell} and convexity give
\[
 b_n\le2a_T,\quad I_{\mathrm{mod}}\le a_TT^2,\quad
 Q_n\le\tfrac12Lh^2\mathbb E|w_n|^2\le4Lh^2R,
 \quad Q_\Sigma\le2b_hT^2,\quad \Delta_U\le2R.
\]
Also $\mathcal W(0)\le E_0/m$ by the transport-entropy inequality
(Lemma~\ref{lem:talagrand}) and entropy data processing to the position
marginal. Divide the estimate \eqref{ent:signed-final} of
Lemma~\ref{ent:weighted-entropy} by $3T^2$, discard its last two
nonpositive terms, and use $LT^2\ge4$.  This yields
\begin{equation}\label{ent:symmetric-average}
 \frac{I_{\mathrm{end}}}{T^2}\le
 \left(\frac56+\frac{1}{3mT^2}\right)E_0
       +\frac53a_T+\left(2+\frac1{24}\right)b_h.
\end{equation}
The triangular law assigns probability $2\omega_n/T^2$ to
$[nh,(n+1)h)$. Randomized rounding chooses $J=n$ or $J=n+1$ with
equal probability on that interval. Hence $I_{\mathrm{end}}/T^2=\sum_{j=0}^N q_j E_j^x$
for the probabilities in \eqref{up:triangular}. Convexity of relative
entropy under mixtures gives
$\KL(\mathcal L(x_J)\|\pi)\le I_{\mathrm{end}}/T^2$, which proves the claim.
\end{proof}

\subsection{Restricted conductance from entropy contraction}

The next lemma converts a relative-entropy contraction for the
proposal of a Metropolis-adjusted kernel into a restricted conductance
bound for the kernel itself. The contraction is required only from a
specific family of initial laws, and it is verified in practice by
proving it for every initial law in a class that contains this family.

\begin{lemma}[Entropy contraction implies restricted conductance]\label{up:reference-conductance}
Let $\bar\pi$ be a probability measure on $\R^d\times\mathcal V$, where
$\mathcal V$ is a measurable space, with first marginal $\pi$ and a
measurable conditional law $\bar\pi(\,\cdot\mid x)$ given the first
coordinate, and let
$K$ be the Metropolis-adjusted kernel with measurable proposal
$\Upsilon:\R^d\times\mathcal V\to\R^d$ and acceptance probability
$\alpha:\R^d\times\mathcal V\to[0,1]$: from $x$, draw
$V\sim\bar\pi(\,\cdot\mid x)$ and move to $\Upsilon(x,V)$ with probability
$\alpha(x,V)$, so that
\begin{equation}\label{up:abstract-kernel}
 K(x,A)=\int_{\mathcal V}
 \bigl[\alpha(x,v)\mathbf1_A(\Upsilon(x,v))+(1-\alpha(x,v))\mathbf1_A(x)\bigr]
 \,\bar\pi(\mathrm dv\mid x).
\end{equation}
Let $s>0$ and $\eta\in(0,1)$, and let $S\subseteq\R^d\times\mathcal V$
be a measurable set with $\bar\pi(S^c)\le\eta s$ on which proposals are
accepted with probability at least $1-\eta$, that is,
$\alpha(x,v)\ge1-\eta$ for all $(x,v)\in S$.
Suppose that, for some $\theta\in[0,1)$ and $\delta_0\ge0$ and every Borel
$A\subseteq\R^d$ with $s\le\pi(A)\le1/2$, the proposal contracts
relative entropy from the conditional law
$p_A:=\bar\pi(\,\cdot\mid\{x\in A\}\cap S)$:
\begin{equation}\label{up:raw-entropy}
 \KL(\Law(\Upsilon(Z))\Vert\pi)\le\theta\,\KL(p_A\Vert\bar\pi)+\delta_0
 \qquad\text{for }Z\sim p_A.
\end{equation}
Then, for every such $A$ and every $\phi_0\in(0,(1-\theta)/2]$ with
\begin{equation}\label{up:entropy-threshold}
 H_2(\phi_0)+\frac{\theta\eta}{1-\eta}+\delta_0\le(1-\theta)\log2,
 \qquad H_2(p)=-p\log p-(1-p)\log(1-p),
\end{equation}
\begin{equation}\label{up:conductance}
 Q_K(A,A^c)\ge(1-\eta)^2\,\phi_0\,\pi(A).
\end{equation}
\end{lemma}

Note that, if $\lambda:=(1-\theta)\log2-\theta\eta/(1-\eta)-\delta_0>0$, then
\eqref{up:entropy-threshold} holds with
$\phi_0=\min\{(1-\theta)/2,\ \lambda/(2\log(2e/\lambda))\}$, since
$H_2(p)\le p\log(e/p)$ on $(0,1/2]$ and
$\log(2e/\lambda)+\log\log(2e/\lambda)\le2\log(2e/\lambda)$. Conversely, no $\phi_0>0$ satisfies
\eqref{up:entropy-threshold} when $\lambda\le0$. 

\begin{proof}
Fix $A$ with $\pi(A)\in[s,1/2]$ and $\phi_0$ as in the statement, and write
$E_A=\{x\in A\}\cap S$ and $m_A=\bar\pi(E_A)$. Since
$m_A\ge\pi(A)-\bar\pi(S^c)$, $\bar\pi(S^c)\le\eta s$ and $s\le\pi(A)$,
\begin{equation}\label{up:conditioned-reference}
 m_A\ge(1-\eta)\pi(A),
 \qquad
 \KL(p_A\Vert\bar\pi)=\log\frac1{m_A}
       \le\log\frac1{\pi(A)}-\log(1-\eta).
\end{equation}
Let $\nu_{\rm prop}$ be the law of the proposal $\Upsilon(Z)$ for $Z\sim p_A$ and
$r=\nu_{\rm prop}(A^c)$. Writing $z=(x,v)$, restricting the integral
below to $E_A$, where $\alpha\ge1-\eta$, and using
$p_A=\mathbf1_{E_A}\bar\pi/m_A$,
\begin{align*}
 Q_K(A,A^c)&=\int_{\R^d\times\mathcal V}\mathbf1_A(x)\,\alpha(z)\,\mathbf1_{A^c}(\Upsilon(z))\,\bar\pi(\mathrm dz)
 \ge(1-\eta)\int_{E_A}\mathbf1_{A^c}(\Upsilon(z))\,\bar\pi(\mathrm dz)\\
 &\ge(1-\eta)\,m_A\int_{\R^d\times\mathcal V}\mathbf1_{A^c}(\Upsilon(z))\,p_A(\mathrm dz)
 =(1-\eta)\,m_A\,r\ge(1-\eta)^2\pi(A)\,r
\end{align*}
by \eqref{up:conditioned-reference}, and it suffices to show $r\ge \phi_0$.

Write $H_A=\log(1/\pi(A))$. Under $x\mapsto\mathbf1_A(x)$, the laws
$\nu_{\rm prop}$ and $\pi$ become Bernoulli laws with success
probabilities $1-r$ and $\pi(A)$, and relative entropy does not
increase, so
\begin{align}
 \KL(\nu_{\rm prop}\Vert\pi)
 &\ge\KL\bigl(\operatorname{Bernoulli}(1-r)\,\big\Vert\,\operatorname{Bernoulli}(\pi(A))\bigr)\notag\\
 &\ge(1-r)\log\frac{1-r}{\pi(A)}+r\log\frac r{1-\pi(A)}\notag\\
 &\ge (1-r)H_A+r\log\frac1{1-\pi(A)}-H_2(r),\label{up:binary-kl}
\end{align}
with $0\log0=0$. On the other hand, by \eqref{up:raw-entropy},
\eqref{up:conditioned-reference}, and $-\log(1-\eta)\le\eta/(1-\eta)$,
\[
 \KL(\nu_{\rm prop}\Vert\pi)
 \le\theta H_A+\frac{\theta\eta}{1-\eta}+\delta_0.
\]
Suppose $r<\phi_0$. Combining the two bounds and moving $\theta H_A$ to the
left gives
\begin{equation}\label{up:squeeze}
 (1-\theta-2r)H_A+r\log\frac1{\pi(A)(1-\pi(A))}-H_2(r)
 \le\frac{\theta\eta}{1-\eta}+\delta_0,
\end{equation}
where we have written $rH_A+r\log(1/(1-\pi(A)))$ as
$r\log(1/(\pi(A)(1-\pi(A))))$ and subtracted $rH_A$ from the first
term. Both coefficients on the left are nonnegative, because
$r<\phi_0\le(1-\theta)/2$, and since $\pi(A)\le1/2$ we have $H_A\ge\log2$
and $\log(1/(\pi(A)(1-\pi(A))))\ge2\log2$. Replacing each by its lower
bound, the first two terms of \eqref{up:squeeze} are at least
$(1-\theta-2r)\log2+2r\log2=(1-\theta)\log2$, so
\[
 H_2(r)\ge(1-\theta)\log2-\frac{\theta\eta}{1-\eta}-\delta_0\ge H_2(\phi_0)
\]
by \eqref{up:entropy-threshold}. But $H_2$ is increasing on $[0,1/2]$
and $r<\phi_0\le1/2$, so $H_2(r)<H_2(\phi_0)$, a contradiction. Hence $r\ge \phi_0$.
\end{proof}

\subsection{Transfer to the sampler}
\label{sec:rhmc-conductance-proof}

We now combine the ingredients of this section. The energy estimates
of Section~\ref{sec:energy-control} show that on $\{H\le R\}$, a set
whose complement has small mass, proposals are accepted with
probability at least $1-\eta$; Proposition~\ref{ent:raw-entropy}
gives the entropy contraction; Lemma~\ref{up:reference-conductance}
converts the two into restricted conductance for the reference kernel
$P_{\mathrm{ref}}$; and Lemma~\ref{up:law-comparison} transfers the
bound to the kernel $P$ of the non-resonant law $\rho$.

\begin{proof}[Proof of Proposition~\ref{up:restricted-conductance}]
The step-size condition gives \eqref{up:proof-step}, hence
$a_*+b_h\le\eta/3$. We apply Lemma~\ref{up:reference-conductance} to $P_{\mathrm{ref}}$ with
$\mathcal V=\R^d\times\{0,\dots,N\}$, $\bar\pi=\pi\otimes\gamma_d\otimes q$,
$\Upsilon(x,v,j)=Q_j(x,v)$, $\alpha(x,v,j)=\alpha_j(x,v)$ and $S=\{H\le R\}$.
The energy estimates above give $\bar\pi(S^c)=\Pi(H>R)\le\eta s$ and
$\alpha_j\ge1-\eta$ on $\{H\le R\}$, as the lemma requires.
Each $p_A$ is supported on $\{H\le R\}$, has finite relative entropy,
and satisfies \eqref{ent:centered} because $H$ is even in $v$; so
Proposition~\ref{ent:raw-entropy} with $T=4/\sqrt m$ gives
\[
 \KL(\Law(x_J)\Vert\pi)
 \le\frac{41}{48}E_0+\frac53a_*+\frac{49}{24}b_h
 \le\frac{41}{48}E_0+\frac{49}{72}\eta,
\]
that is, \eqref{up:raw-entropy} with $\theta=41/48$ and $\delta_0=49\eta/72$.
Take $\phi_0=1/[50(1-\eta)^2]$. Then $\phi_0\le7/96=(1-\theta)/2$, and since
$\eta=c/600\le1/600$,
\[
 H_2\!\left(\frac1{50(1-\eta)^2}\right)
 +\frac{41\eta}{48(1-\eta)}+\frac{49}{72}\eta
 <\frac7{48}\log2,
\]
because the left side is increasing in $\eta$ and is less than
$101/1000$ at $\eta=1/600$, whereas the right side is greater than
$101/1000$. This is \eqref{up:entropy-threshold}, so \eqref{up:conductance}
yields $Q_{P_{\mathrm{ref}}}(A,A^c)\ge\pi(A)/50$ for
$s\le\pi(A)\le1/2$.

Let $w$ be the conditional law of $J$ given $\tau\le\tau_*$.
Its trajectories have physical integration time at most $T_*$.
For $s\le\pi(A)\le1/2$, \eqref{up:shell-work} gives
\begin{align*}
 \operatorname{Rej}_w(\mathbf1_A)
 &\le\eta\pi(A)+\Pi(H>R)\\
 &\le\eta\pi(A)+\eta s\le2\eta\pi(A).
\end{align*}
By \eqref{up:law-constants}, Lemma~\ref{up:law-comparison} gives
\[
 Q_{P_w}(A,A^c)\ge(c/50-2\eta)\pi(A)=\frac c{60}\pi(A),
\]
since $\eta=c/600$. The full kernel $P$ gives this conditional kernel weight
$p_S=\mathbb P(\tau\le\tau_*)$. The remaining component has nonnegative
stationary flow. Therefore
\begin{equation}\label{up:actual-conductance}
 Q_P(A,A^c)\ge\phi_P\pi(A),\qquad
 \phi_P=\frac{p_Sc}{60}>0.
\end{equation}
\end{proof}

\section{Necessity of Hessian regularity for fixed step sizes}\label{sec:lower}

Theorem~\ref{lower:thm} shows that Hessian regularity cannot be omitted
from a uniform square-root gradient bound for Metropolis-adjusted
randomized HMC with a fixed Verlet step size. The counterexamples are
strongly convex potentials with Lipschitz gradients, and the initial
density is at most twice the target density. Low acceptance probabilities
produce the lower bound.

The proof first constructs a set of expanding Verlet trajectories for
one fast coordinate. With $d_f$ independent copies, expansion in at least
one coordinate makes acceptance unlikely after sufficiently many Verlet
steps. The Metropolis test
accepts or rejects all coordinates together, so a large energy increase
in a fast coordinate can keep the slow coordinate on its initial
half-line. We bound the gradient cost of obtaining the first acceptance.

\subsection{Construction of the target}

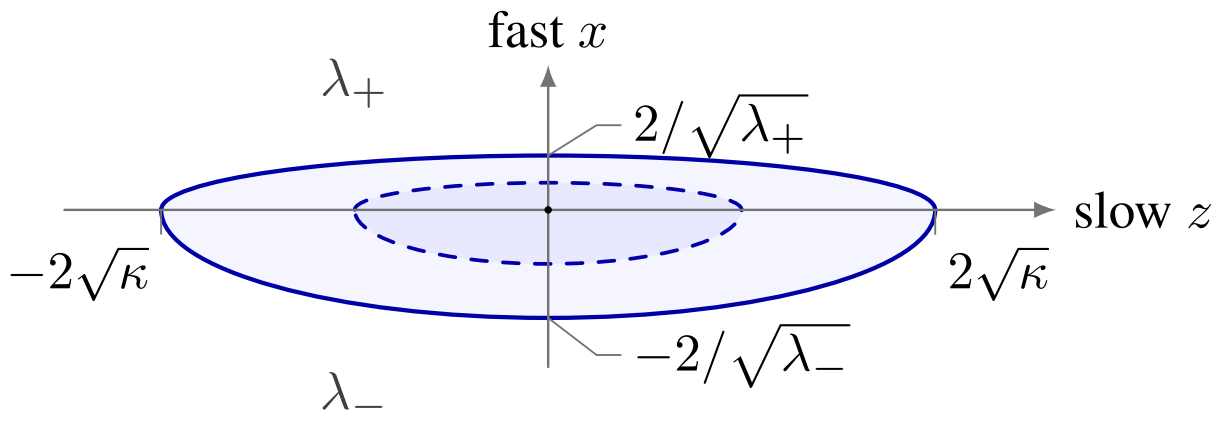
\begin{figure}[!t]
\centering
\begingroup
\begin{tikzpicture}[
  x=0.08cm,y=0.08cm,
  line cap=round,line join=round,font=\small
]
  \path[fill=blue!4]
    (-32,0)
    arc[start angle=180,end angle=0,
        x radius=2.56cm,y radius=0.359488037310cm]
    arc[start angle=0,end angle=-180,
        x radius=2.56cm,y radius=0.714541230709cm]
    --cycle;
  \path[fill=blue!9]
    (-16,0)
    arc[start angle=180,end angle=0,
        x radius=1.28cm,y radius=0.179744018655cm]
    arc[start angle=0,end angle=-180,
        x radius=1.28cm,y radius=0.357270615354cm]
    --cycle;

  \draw[draw=blue!65!black,line width=0.8pt]
    (-32,0)
    arc[start angle=180,end angle=0,
        x radius=2.56cm,y radius=0.359488037310cm]
    arc[start angle=0,end angle=-180,
        x radius=2.56cm,y radius=0.714541230709cm];
  \draw[draw=blue!65!black,line width=0.7pt,dashed]
    (-16,0)
    arc[start angle=180,end angle=0,
        x radius=1.28cm,y radius=0.179744018655cm]
    arc[start angle=0,end angle=-180,
        x radius=1.28cm,y radius=0.357270615354cm];

  \draw[-latex,black!55,line width=0.45pt]
    (-40,0)--(42,0) node[right,text=black] {slow $z$};
  \draw[-latex,black!55,line width=0.45pt]
    (0,-13)--(0,12) node[above,text=black] {fast $x$};

  \draw[black!55,line width=0.35pt] (-32,0)--(-32,-2);
  \draw[black!55,line width=0.35pt] (32,0)--(32,-2);
  \node[below left,inner sep=2pt] at (-32,-2)
    {$-2\sqrt\kappa$};
  \node[below right,inner sep=2pt] at (32,-2)
    {$2\sqrt\kappa$};

  \draw[black!55,line width=0.35pt]
    (0,4.493600466380)--(4,7)--(6,7);
  \node[anchor=west,inner sep=2pt] at (6,7)
    {$2/\sqrt{\lambda_+}$};
  \draw[black!55,line width=0.35pt]
    (0,-8.931765383858)--(4,-12)--(6,-12);
  \node[anchor=west,inner sep=2pt] at (6,-12)
    {$-2/\sqrt{\lambda_-}$};

  \fill (0,0) circle[radius=0.7pt];
  \node[anchor=south,text=black!75] at (-16,7) {$\lambda_+$};
  \node[anchor=north,text=black!75] at (-16,-12) {$\lambda_-$};
\end{tikzpicture}
\endgroup

\caption{One fast coordinate and the slow coordinate for $h=1$ and
$\kappa=256$. The curves are
$\psi_h(x)+z^2/(2\kappa)=1/2$ and $2$; their densities relative to
the mode are $e^{-1/2}$ and $e^{-2}$, respectively.}
\label{lower:geometry}
\end{figure}

The example separates the role of the fast and slow coordinates.
For a prescribed step $0<h\leq1$, choose positive curvatures
$\lambda_->0$ and $\lambda_+>\lambda_-$ as specified below, and let
\begin{equation}\label{lower:psi-shape}
\psi_h(x)=
\begin{cases}
\lambda_+x^2/2,&x\geq0,\\
\lambda_-x^2/2,&x<0.
\end{cases}
\end{equation}
For $d_f$ fast coordinates and two additional coordinates $w,z$, the
full potential is
\begin{equation}\label{lower:target}
U_h(x_1,\ldots,x_{d_f},w,z)
   =\sum_{i=1}^{d_f}\psi_h(x_i)+\frac{w^2}{2}
                          +\frac{z^2}{2\kappa}.
\end{equation}
The target position coordinates are independent.  Each $x_i$ has a
Gaussian density on each half-line, with the wider half on the negative
side.  The slow coordinate is exactly $z\sim N(0,\kappa)$; the extra
coordinate $w\sim N(0,1)$ makes the smallest Lipschitz constant of
$\nabla U_h$ equal to one. Figure~\ref{lower:geometry} shows one fast
coordinate together with the slow coordinate.

\subsection{An expanding Verlet orbit}

We now give the promised choice of $\lambda_\pm$ and prove the
theorem.  First define the integer cycle lengths from the prescribed
step:
\[
p=\left\lceil\frac{2\uppi}{h}\right\rceil\geq7,
\qquad p_{\mathrm{cyc}}=3p.
\]
The orbit will spend $p$ steps on the positive side and $2p$ on the
negative side in each cycle. At the unperturbed angles
$\alpha=\uppi/p$ and $\beta=\uppi/(2p)$, these are exact half-periods,
so a full cycle returns to its starting pair of positions.
Perturb the angles by a parameter $0<\delta\leq(80p^2)^{-1}$, whose
value we fix at the end of the proof.

\begin{equation}\label{lower:angles}
\alpha=\frac\uppi p,\quad
\beta=\frac\uppi{2p},\quad
a=\alpha+\frac\delta p,\quad
b=\beta-\frac\delta{4p}.
\end{equation}
Choose
\begin{equation}\label{lower:psi}
\lambda_+=\frac{4\sin^2(a/2)}{h^2},\qquad
\lambda_-=\frac{4\sin^2(b/2)}{h^2},
\end{equation}
and abbreviate $\psi_h$ by $\psi$ in the scalar calculations.
The function \eqref{lower:psi-shape} is continuously differentiable.
Every secant slope of $\psi_h'$ lies in $[\lambda_-,\lambda_+]$:
a secant across zero is a weighted average of the two slopes. Thus
$\psi_h$ is $\lambda_-$-strongly convex and $\lambda_+$-smooth.
We will use the bounds
\begin{equation}\label{lower:curvatures}
\frac1{16}\leq\lambda_+<1,\qquad
\frac1{64}\leq\lambda_-<\lambda_+,
\qquad 2\uppi\leq ph<2\uppi+1.
\end{equation}
Indeed, $a\geq\uppi/p$ and $b\geq3/(2p)$.  The inequality
$\sin t\geq2t/\uppi$ on $[0,\uppi/2]$ gives
$\lambda_+\geq4/(ph)^2\geq1/16$ and
$\lambda_-\geq9/[\uppi^2(ph)^2]\geq1/64$.  The last numerical
inequalities follow, for example, from $\uppi^2<10$ and
$(2\uppi+1)^2<54$.  Also
$\lambda_+\leq[(\uppi+\delta)/(ph)]^2<1$.

Write $\Phi$ for one scalar velocity-Verlet step for $\psi$.
The position-velocity coordinates $(x_n,v_n)$ can be replaced by consecutive positions
$(x_n,x_{n-1})$, since
\begin{equation}\label{lower:native}
x_{n-1}=x_n-hv_n-\frac{h^2}{2}\psi'(x_n).
\end{equation}
The Verlet recurrence is
\[
x_{n+1}=2x_n-x_{n-1}-h^2\psi'(x_n).
\]
Consequently the pair evolves by the matrix
\begin{equation}\label{lower:A}
\begin{pmatrix}x_{n+1}\\x_n\end{pmatrix}
=A(\theta)\begin{pmatrix}x_n\\x_{n-1}\end{pmatrix},
\qquad
A(\theta)=\begin{pmatrix}2\cos\theta&-1\\1&0\end{pmatrix},
\end{equation}
where $\theta=a$ when $x_n>0$ and $\theta=b$ when $x_n<0$.

\begin{figure}[!t]
\centering
\begingroup
\newcommand{\ReturnConeAxes}{%
  \path[fill=blue!7] (0,0)--(2.23,-6.69)--(2.23,-8.92)--cycle;
  \draw[black!56,line width=0.55pt] (0,0)--(2.23,-6.69)
    node[right,inner sep=2pt] {$x_{n-1}=-3x_n$};
  \draw[black!56,line width=0.55pt] (0,0)--(2.23,-8.92)
    node[right,inner sep=2pt] {$x_{n-1}=-4x_n$};
  \draw[-{Latex[length=1.6mm]},black!62,line width=0.45pt]
    (-0.33,0)--(2.55,0) node[right,inner sep=2pt] {$x_n$};
  \draw[-{Latex[length=1.6mm]},black!62,line width=0.45pt]
    (0,-9.28)--(0,1.13) node[above,inner sep=2pt] {$x_{n-1}$};
  \node[blue!55!black,inner sep=1pt] at (2.04,-7.55) {$\mathcal C$};
}
\begin{tikzpicture}[
  x=1.5cm,y=0.37cm,font=\small,line join=round,line cap=round,
  return guide/.style={black!35,densely dotted,line width=0.5pt},
  return leader/.style={black!45,line width=0.4pt}
]
  \begin{scope}
    \ReturnConeAxes
    \node[font=\small\bfseries] at (1.38,2.8)
      {Unperturbed return, $\delta=0$};
    \draw[return guide] (1,0)--(1,-3.5);
    \node[above,inner sep=2pt] at (1,0) {$x_0=x_{p_{\mathrm{cyc}}}$};
    \draw[return leader] (1,-3.5)--(1.14,-2.63);
    \node[right,inner sep=2pt] at (1.14,-2.63)
      {$(x_0,x_{-1})=(x_{p_{\mathrm{cyc}}},x_{p_{\mathrm{cyc}}-1})$};
    \fill (1,-3.5) circle[radius=1.6pt];
    \draw[blue!65!black,line width=0.55pt] (1,-3.5) circle[radius=2.7pt];
  \end{scope}
  \begin{scope}[xshift=7.25cm]
    \ReturnConeAxes
    \node[font=\small\bfseries] at (1.38,2.8)
      {Perturbed return, $\delta>0$};
    \draw[return guide] (1,0)--(1,-3.5);
    \draw[return guide] (1.8,0)--(1.8,-6.667716358115);
    \node[above,inner sep=2pt] at (1,0) {$x_0$};
    \node[above,inner sep=2pt] at (1.8,0) {$x_{p_{\mathrm{cyc}}}$};
    \draw[return leader] (1,-3.5)--(1.12,-2.62);
    \node[right,fill=white,inner sep=2pt] at (1.12,-2.62) {$(x_0,x_{-1})$};
    \draw[return leader] (1.8,-6.667716358115)--(1.94,-4.72);
    \node[right,inner sep=2pt] at (1.94,-4.72) {$(x_{p_{\mathrm{cyc}}},x_{p_{\mathrm{cyc}}-1})$};
    \draw[-{Latex[length=2.2mm,width=1.65mm]},
          blue!65!black,line width=0.95pt,shorten <=2pt,shorten >=2pt]
      (1,-3.5)--(1.8,-6.667716358115);
    \fill (1,-3.5) circle[radius=1.6pt];
    \fill[blue!65!black] (1.8,-6.667716358115) circle[radius=1.9pt];
  \end{scope}
\end{tikzpicture}
\endgroup

\caption{Return after $p_{\mathrm{cyc}}=3p$ steps in the position-pair plane. The
unperturbed cycle returns to its starting point. The perturbed cycle
returns to the cone with a larger first coordinate. The displacement
in the right panel is magnified.}
\label{fig:return-map}
\end{figure}
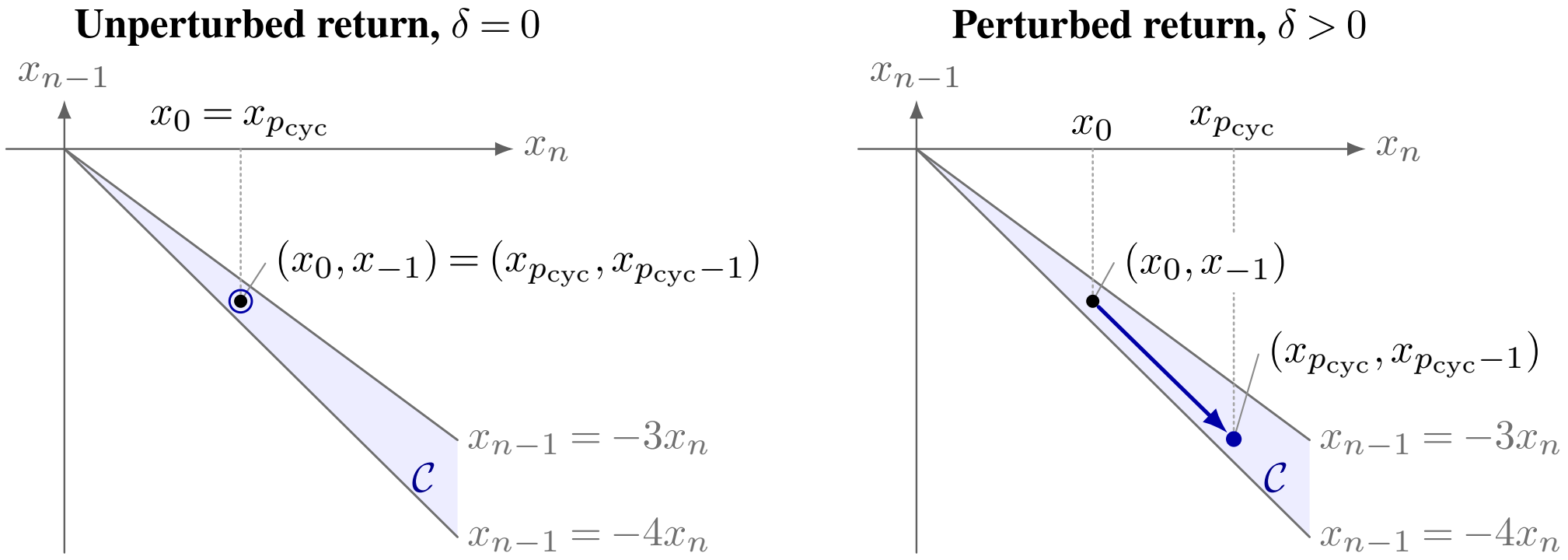

The scalar Verlet recurrence belongs to the family of piecewise-linear
area-preserving maps studied by Lagarias and Rains \cite{LR}.
Their Theorem~2.4 classifies the dynamics for rational rotation numbers,
including cases with expanding rays and divergent orbits. We give an
explicit parameter choice and quantitative estimates on an invariant
cone, the expansion rate, and the intermediate iterates, which we then
use to obtain the sampling lower bound.
Figure~\ref{fig:return-map} shows the return of a point in the cone
after one cycle.

\begin{lemma}[A return cone]\label{lower:cone}
Let
\[
\mathcal C=\{(x,y):x>0,\ -4x\leq y\leq-3x\}.
\]
Starting from any pair in $\mathcal C$, the recurrence
\eqref{lower:A} uses $A(a)$ for $p$ consecutive steps and then $A(b)$
for $2p$ consecutive steps.  It returns to $\mathcal C$ after $p_{\mathrm{cyc}}=3p$
steps.  Moreover, for every integer $k\geq0$,
\begin{equation}\label{lower:expansion}
(x_{kp_{\mathrm{cyc}}},x_{kp_{\mathrm{cyc}}-1})\in\mathcal C,\qquad
x_{kp_{\mathrm{cyc}}}\geq \chi^k x_0,\qquad
\chi=1+\frac45\delta b_0,
\end{equation}
where
$b_0=\tan(\beta/2)/(2\cos\beta)$.
Every intermediate position in a cycle has absolute value at least
half the positive position at the start of that cycle.
\end{lemma}

\begin{proof}
We first verify the signs of the intermediate positions, then prove
invariance and expansion of the cone.

\paragraph*{The unperturbed sign sequence}
For $0<\theta<\uppi$, define
\[
C(\theta)=
\begin{pmatrix}\cot\theta&-\csc\theta\\
                 \csc\theta&-\cot\theta\end{pmatrix}.
\]
Since $A(\theta)=\cos\theta\,I+\sin\theta\,C(\theta)$ and
$C(\theta)^2=-I$,
\begin{equation}\label{lower:powers}
A(\theta)^j=\cos(j\theta)I+\sin(j\theta)C(\theta).
\end{equation}
At $\delta=0$, the first coordinate after $j$ positive steps is
\[
\frac{\sin((j+1)\alpha)x-\sin(j\alpha)y}{\sin\alpha}.
\]
It is positive for $0\leq j<p$, and $A(\alpha)^p=-I$.
Starting from $(-x,-y)$, the analogous formula with $\beta$ is
negative for the next $2p$ iterates, and $A(\beta)^{2p}=-I$.
All these positions have magnitude at least $x$. To check this last
assertion, use $-y\geq x$ and
\[
\sin((j+1)\theta)+\sin(j\theta)\geq\sin\theta,
\qquad 0\leq j<\uppi/\theta,
\]
for $\theta=\alpha$ or $\beta$. The left-hand side is concave in $j$
on $[0,\uppi/\theta-1]$ and equals $\sin\theta$ at both endpoints.

\paragraph*{Preservation of the sign sequence}
Use the induced infinity norm for matrices. For $0<\theta<\uppi/2$,
\[
 \|C(\theta)\|_\infty=\cot(\theta/2)\leq\frac2\theta,
 \qquad
 \|C'(\theta)\|_\infty=\frac1{2\sin^2(\theta/2)}
                    \leq\frac{\uppi^2}{2\theta^2}.
\]
Since $a\geq\alpha=\uppi/p$ and $b\geq3/(2p)$, these give
\begin{equation}\label{lower:C-bounds}
\begin{aligned}
 &\|C(\theta)\|_\infty\leq2p/3,\quad
   \|C'(\theta)\|_\infty\leq p^2/2
       &&(\alpha\leq\theta\leq a),\\
 &\|C(\theta)\|_\infty\leq4p/3,\quad
   \|C'(\theta)\|_\infty\leq20p^2/9
       &&(b\leq\theta\leq\beta).
\end{aligned}
\end{equation}
By \eqref{lower:powers}, the derivative of $A(\theta)^j$ has norm at most
\[
 j\bigl(1+\|C(\theta)\|_\infty\bigr)+\|C'(\theta)\|_\infty.
\]
Multiplying by $\delta/p$ for $j\leq p$, or by $\delta/(4p)$ for
$j\leq2p$, shows that each partial block changes by at most $2p\delta$.
The two full blocks are
\begin{align*}
\mathsf M_a&=A(a)^p=-\cos\delta\,I-\sin\delta\,C(a),\\
\mathsf M_b&=A(b)^{2p}=-\cos(\delta/2)I+\sin(\delta/2)C(b).
\end{align*}
In particular, $\|\mathsf M_a+I\|_\infty\leq p\delta$ and
$\|\mathsf M_b+I\|_\infty\leq p\delta$, by $1-\cos t\leq t^2/2$ and
$|\sin t|\leq t$ for $t\geq0$. During the negative block, the
comparison also includes the preceding positive block. Using
$\|A(b)^j\|_\infty\leq1+4p/3\leq2p$ gives
\[
\|A(b)^j\mathsf M_a+A(\beta)^j\|_\infty
\leq(2p)(p\delta)+2p\delta\leq\frac52p^2\delta.
\]
Since $\|(x,y)\|_\infty\leq4x$ on $\mathcal C$, the perturbation
changes the position at any iterate by at most
\begin{equation}\label{lower:sign-margin}
10p^2\delta x\leq\frac{x}{8}<x/2.
\end{equation}
Thus all prescribed signs remain strict and every intermediate
magnitude is at least $x/2$.  The matrix products below consequently
describe the trajectory generated by \eqref{lower:A}.

\paragraph*{Expansion of the return map}
Put $\mathsf M_{\mathrm{lin}}=C(\alpha)-\tfrac12C(\beta)$. We will show that
\[
 \mathsf M_b\mathsf M_a=I+\delta\mathsf M_{\mathrm{lin}}+\mathsf E_\delta,
 \qquad \|\mathsf E_\delta\|_\infty<\frac{\delta b_0}{20}.
\]
We will use this bound to show that both cone boundaries move inward
and that the first coordinate increases.
Write $\Delta_a=\mathsf M_a+I$ and $\Delta_b=\mathsf M_b+I$. Using
$|\sin t-t|\leq t^3/6$ for $t\geq0$ and the derivative
bound in \eqref{lower:C-bounds} gives
\begin{align*}
\|\Delta_a+\delta C(\alpha)\|_\infty
 &\leq\frac{\delta^2}{2}+\frac{p\delta^3}{9}
                         +\frac{p\delta^2}{2}
 \leq\frac23p\delta^2,\\
\|\Delta_b-\tfrac\delta2 C(\beta)\|_\infty
 &\leq\frac{\delta^2}{8}+\frac{p\delta^3}{36}
                         +\frac5{18}p\delta^2
 \leq\frac13p\delta^2.
\end{align*}
To bound the product term, multiply the matrices before taking norms:
\[
 C(b)C(a)=\frac1{\sin a\sin b}
 \begin{pmatrix}
   \cos a\cos b-1&\cos a-\cos b\\
   \cos a-\cos b&\cos a\cos b-1
 \end{pmatrix}.
\]
Since $0<b<a\leq\uppi/3$ and $a/b\leq3$,
\[
 \|C(b)C(a)\|_\infty
 =\tan(a/2)\cot(b/2)
 \leq\frac{2a}{\sqrt3\,b}<4.
\]
Expanding $\Delta_b\Delta_a$ with the displayed formulas for
$\mathsf M_a$ and $\mathsf M_b$ therefore gives
\[
 \|\Delta_b\Delta_a\|_\infty
 \leq\frac{\delta^4}{16}+\frac5{12}p\delta^3+2\delta^2
 \leq3\delta^2.
\]
Thus the return matrix $I-\Delta_a-\Delta_b+\Delta_b\Delta_a$ satisfies
\begin{equation}\label{lower:M-expansion}
\mathsf M_b\mathsf M_a=I+\delta\mathsf M_{\mathrm{lin}}+\mathsf E_\delta,
\qquad \|\mathsf E_\delta\|_\infty
 \leq(p+3)\delta^2\leq\frac32p\delta^2.
\end{equation}
Writing $c_0=\cos\beta$ and using $\alpha=2\beta$ gives the exact
identity
\begin{equation}\label{lower:B}
\mathsf M_{\mathrm{lin}}=b_0\begin{pmatrix}-(1+c_0)&-1\\1&1+c_0\end{pmatrix},
\qquad b_0=\frac{1-c_0}{2\sin\beta\,c_0}
          =\frac{\tan(\beta/2)}{2c_0}\geq\frac\uppi{8p}.
\end{equation}
In particular,
\begin{equation}\label{lower:R-small}
\frac{\|\mathsf E_\delta\|_\infty}{\delta b_0}
\leq\frac{12}{\uppi}p^2\delta\leq\frac{12}{80\uppi}<\frac1{20}.
\end{equation}

On the boundary $y=-4x$, the first-order change of $y+4x$ under
$I+\delta\mathsf M_{\mathrm{lin}}$ is $\delta b_0(9-8c_0)x\geq\delta b_0x$.
The error contributed by $\mathsf E_\delta$ is at most
$20\|\mathsf E_\delta\|_\infty x$, so the image lies on the inward side
of this boundary.  On $y=-3x$, the first-order change of $y+3x$ is
$\delta b_0(4-6c_0)x$.  Since
$c_0\geq\cos(\uppi/8)>\sqrt3/2$, this is at most
$-(3\sqrt3-4)\delta b_0x$; the error is at most
$16\|\mathsf E_\delta\|_\infty x<(4/5)\delta b_0x$.
Since $3\sqrt3-4>4/5$, this boundary also moves inward.
Finally, throughout the cone the first coordinate of the image is at
least
\[
\bigl[1+\delta b_0(2-c_0)-4\|\mathsf E_\delta\|_\infty\bigr]x
\geq(1+\tfrac45\delta b_0)x=\chi x>0.
\]
For the return matrix $M=\mathsf M_b\mathsf M_a$, the outgoing slope
as a function of the incoming slope $s=y/x$ is
$(M_{21}+M_{22}s)/(M_{11}+M_{12}s)$. Its derivative is
$\det M/(M_{11}+M_{12}s)^2>0$, since $\det M=1$ and the outgoing
first coordinate is positive. The two boundary inequalities
therefore put every outgoing slope in $[-4,-3]$. The return preserves
$\mathcal C$, and iteration proves \eqref{lower:expansion}.
\end{proof}

\subsection{Energy growth and rejection}\label{sec:rhmc-rejection-cutoff}

The crossing event $E$ below has probability proportional to $h$.
The union of its images over one cycle has probability bounded below
independently of $h$.
Let
\begin{equation}\label{lower:E}
E=\left\{(x,v):\frac h{10}\leq x\leq\frac h5,
       \ -4x\leq\cos(a)x-hv\leq-3x\right\},
\qquad
E_* =\bigcup_{s=0}^{p_{\mathrm{cyc}}-1}\Phi^{\circ s}(E).
\end{equation}
Let $H_f(x,v)=\psi(x)+v^2/2$, and let $\Pi_f$ be the joint law of
position and velocity with density proportional to $e^{-H_f(x,v)}$.

\begin{lemma}[A set of expanding trajectories]\label{lower:phase-event}
The $p_{\mathrm{cyc}}$ images in \eqref{lower:E} are disjoint up to boundaries,
$H_f<2$ on $E_*$, and
\begin{equation}\label{lower:E-prob}
\Pi_f(E_*)\geq q_*.
\end{equation}
For every $z\in E_*$ and every integer $n\geq0$,
\begin{equation}\label{lower:future-energy}
H_f(\Phi^{\circ n}(z))\geq\frac1{20}\chi^{2\lfloor n/p_{\mathrm{cyc}}\rfloor}.
\end{equation}
\end{lemma}

\begin{proof}
The interval of velocities in $E$ has width $x/h$, so its Lebesgue
area is
\[
|E|=\int_{h/10}^{h/5}\frac{x}{h}\,dx=\frac{3h}{200}.
\]
Moreover, $0\leq v\leq1$ there.  Each image has the same area,
because a Verlet step is a composition of volume-preserving shears.
If $E$ met $\Phi^{\circ j}(E)$ for $0<j<p_{\mathrm{cyc}}$, an $E$-trajectory would have a
second crossing with $x_j>0>x_{j-1}$ before the end of its first
cycle.  Lemma~\ref{lower:cone} excludes such a crossing.  Invertibility
then gives disjointness of any two images.

For an iterate with angle $\theta$, introduce the quadratic form
\[
Q_\theta(x,y)=\frac{x^2+y^2-2\cos\theta\,xy}{h^2}.
\]
It is invariant under $A(\theta)$.  From
$v=(\cos\theta\,x-y)/h$ and
$\lambda=4\sin^2(\theta/2)/h^2$, direct expansion gives
\begin{equation}\label{lower:Q-energy}
Q_\theta=v^2+\lambda(1-h^2\lambda/4)x^2,
\qquad
\frac{Q_\theta}{2}\leq H_f
\leq\frac{Q_\theta}{2(1-h^2\lambda/4)}.
\end{equation}
On $E$, $Q_a\leq1+1/25=26/25$.
At the positive-to-negative crossing one has $xy<0$ and
$\cos b>\cos a>0$.  Thus changing from $Q_a$ to $Q_b$ increases
the value by a factor smaller than two, since $x^2+y^2\geq2|xy|$ and
\[
\frac{x^2+y^2+2\cos b\,|xy|}{x^2+y^2+2\cos a\,|xy|}<2.
\]
There is only this one curvature switch before time $p_{\mathrm{cyc}}$.
Throughout $E_*$, therefore, $Q\leq52/25$.  Since
$h^2\lambda<1$, \eqref{lower:Q-energy} gives $H_f<2$.

The normalizing constant of $\Pi_f$ is
\[
Z_f=\uppi(\lambda_+^{-1/2}+\lambda_-^{-1/2})\leq12\uppi.
\]
Using disjointness, area preservation and $p_{\mathrm{cyc}}h\geq6\uppi$,
\[
\Pi_f(E_*)\geq
\frac{3p_{\mathrm{cyc}}h}{200}\frac{e^{-2}}{12\uppi}
\geq\frac9{1200}e^{-2}>10^{-3}=q_*.
\]

To obtain the future bound, start first in $E$ and consider a cycle
beginning with $(x,y)\in\mathcal C$.  At its positive crossing,
$Q_a(x,y)\geq10x^2/h^2$.  The quadratic form is constant on each
constant-curvature block and increases at the positive-to-negative
switch.  Thus, until the next positive crossing,
$H_f\geq5(x/h)^2$.  At the next positive crossing use the new cycle
and Lemma~\ref{lower:cone}.  Since $x_0\geq h/10$, this proves
\[
H_f(\Phi^{\circ n}(z))\geq \chi^{2\lfloor n/p_{\mathrm{cyc}}\rfloor}/20,
\qquad z\in E.
\]
If $z=\Phi^{\circ s}(z_0)\in E_*$ with $0\leq s<p_{\mathrm{cyc}}$, apply the same
bound at time $n+s$ and use
$\lfloor(n+s)/p_{\mathrm{cyc}}\rfloor\geq\lfloor n/p_{\mathrm{cyc}}\rfloor$.
\end{proof}

\paragraph*{Rejection cutoff}
Return to the product target \eqref{lower:target}.
By \eqref{lower:curvatures}, for $\kappa\geq64$ its optimal
curvature bounds are $m=\kappa^{-1}$ and $L=1$.
Lemma~\ref{up:reversibility} applies to this target:
the kick and drift shears preserve volume, and the Verlet step is
reversible under the momentum flip. Hence $P$ is reversible for
$\pi$. Let $\Pi=\pi\otimes N(0,I_d)$ and write $H_0$
for the energy of a point drawn from $\Pi$.

Under $\Pi$ the $d_f$ position-velocity pairs $(x_i,v_i)$ are independent, so
Lemma~\ref{lower:phase-event} yields
\begin{equation}\label{lower:no-expander}
\Pi(\text{no fast coordinate belongs to }E_*)
\leq(1-q_*)^{d_f}\leq\kappa^{-6}.
\end{equation}
Each scalar position energy has the Gamma$(1/2,1)$ law.  For the
asymmetric coordinate this follows by changing variables from either
half-line to $\lambda_\pm x^2/2$: conditioning on either sign gives
the same Gamma law.  Each velocity energy has the same law,
independently.  Therefore $H_0$ has the Gamma$(d,1)$ law and
$\mathbb E_\Pi e^{H_0/2}=2^d$.
Set
\[
H_{\mathrm{cut}}=2d+12\log\kappa.
\]
Markov's inequality gives
\begin{equation}\label{lower:H-tail}
\Pi(H_0>H_{\mathrm{cut}})\leq2^d e^{-H_{\mathrm{cut}}/2}\leq\kappa^{-6}.
\end{equation}
Let $\mathcal B$ be the union of the two exceptional phase events in
\eqref{lower:no-expander} and \eqref{lower:H-tail}.  Then
\begin{equation}\label{lower:B-prob}
\Pi(\mathcal B)\leq2\kappa^{-6}.
\end{equation}

Define
\begin{equation}\label{lower:j0}
\ell_\kappa=\log\!\bigl(20[H_{\mathrm{cut}}+6\log\kappa]\bigr),\qquad
j_0=p_{\mathrm{cyc}}\left\lceil\frac{\ell_\kappa}{2\log\chi}\right\rceil.
\end{equation}
The corresponding cutoff for the physical integration time is $t_0=j_0h$.
If an incoming phase point is outside $\mathcal B$, at least one fast
coordinate satisfies \eqref{lower:future-energy}. For every $Jh\geq t_0$
its final energy is at least $H_{\mathrm{cut}}+6\log\kappa$. All other
coordinate energies are nonnegative, whereas the initial total energy
is at most $H_{\mathrm{cut}}$. Hence the change in Hamiltonian energy
$\Delta H$ satisfies
\begin{equation}\label{lower:alpha-after-cutoff}
\Delta H\geq6\log\kappa,\qquad 1\wedge e^{-\Delta H}\leq\kappa^{-6}
\quad\text{on }\mathcal B^c\cap\{J\geq j_0\}.
\end{equation}

A rejected transition leaves the position fixed, and the next attempt
redraws the momentum. We therefore need an acceptance bound conditional
on the position. For
$V\sim N(0,I_d)$, define
\[
p_{\mathrm{bad}}(x)=\mathbb P_V((x,V)\in\mathcal B),\qquad
\mathcal G=\{x:p_{\mathrm{bad}}(x)\leq\kappa^{-2}\}.
\]
From \eqref{lower:B-prob} and Markov's inequality,
\begin{equation}\label{lower:G}
\pi(\mathcal G^c)\leq2\kappa^{-4}.
\end{equation}
At every fixed $x\in\mathcal G$, a fresh attempt accepts with probability
at most
\begin{equation}\label{lower:a0}
a_0=\min\{1,\mathbb P(J<j_0)+\kappa^{-2}+\kappa^{-6}\}.
\end{equation}
The relevant comparison is between $t_0$ and the duration scale $T$.
Since $J\geq\lfloor\tau/h\rfloor$,
\[
\mathbb P(J<j_0)\leq\mathbb P(\tau<t_0)
\leq\min\{1,2t_0/T\}
\leq\min\{1,t_0/\sqrt\kappa\}.
\]
Each cycle takes physical time $p_{\mathrm{cyc}}h=3ph$.
The bounds $b_0<1$ and $\delta\leq(80p^2)^{-1}$ give
$0<\log\chi\leq1/32$. Since $3ph<22$ and $\ell_\kappa\geq1$,
\eqref{lower:j0} gives
\[
t_0+2
\leq\frac{3ph}{2}\frac{\ell_\kappa}{\log\chi}+3ph+2
\leq\left(11+\frac{24}{32}\right)\frac{\ell_\kappa}{\log\chi}
<12\,\frac{\ell_\kappa}{\log\chi}.
\]
The two exceptional probabilities in \eqref{lower:a0} sum to at most
$2/\sqrt\kappa$. Together with the short-duration bound, this gives
\begin{equation}\label{lower:accept-bound}
a_0\leq\min\left\{1,
       \frac{12\ell_\kappa}{\sqrt\kappa\log\chi}\right\}.
\end{equation}

\subsection{Mixing and gradient lower bounds}

Let $A_+=\{z\geq0\}$, where $z$ is the slow coordinate in
\eqref{lower:target}, and let
\begin{equation}\label{lower:initial}
\nu=\pi(\,\cdot\mid A_+)=2\mathbf1_{A_+}\pi.
\end{equation}
Its density with respect to $\pi$ has $L^2(\pi)$ norm $\sqrt2$, and
$\KL(\nu\Vert\pi)=\log2$.

By \eqref{lower:G} and $d\nu/d\pi\le2$,
\begin{equation}\label{lower:initial-bad}
\mathbb P_\nu(X_0\notin\mathcal G)\leq4\kappa^{-4}.
\end{equation}
Until the first acceptance, the position remains exactly $X_0$.
On $\{X_0\in\mathcal G\}$ the conditional chance that each transition
accepts is at most $a_0$, regardless of the earlier rejections,
because its momentum, duration, rounding variable, and uniform
variable for the Metropolis test are drawn independently at each attempt.

Fix an integer $n\geq0$. For $1\leq i\leq n$, let $F_i$ be the event
that $X_0\in\mathcal G$ and the first $i-1$ transitions all reject.
The probability that $F_i$ occurs and transition $i$ accepts is at most
$a_0\mathbb P(F_i)\leq a_0$. Summing over $i=1,\ldots,n$ and using
\eqref{lower:initial-bad} gives
\begin{equation}\label{lower:first-accept}
\mathbb P(\text{an acceptance occurs in the first }n\text{ transitions})
\leq4\kappa^{-4}+na_0.
\end{equation}
On the complementary event the output still belongs to $A_+$. Since
$\pi(A_+)=1/2$,
\begin{equation}\label{lower:TV}
\TV(\nu P^n,\pi)
\geq\frac12-4\kappa^{-4}-na_0.
\end{equation}
Consequently, for $\kappa\geq64$, TV accuracy $1/4$ requires
\begin{equation}\label{lower:attempts}
n\geq\frac1{5a_0}
\geq \frac15\max\left\{1,
        \frac{\sqrt\kappa\log\chi}{12\ell_\kappa}\right\}.
\end{equation}

\paragraph*{Gradient count}
The $i$th transition costs $J_i$ new gradient evaluations.
Since $\mathbb E\tau\geq T/3$ and $J\geq\tau/h-1$,
\begin{equation}\label{lower:EJ}
\mathbb E J\geq\frac{2\sqrt\kappa}{3h}-1
              \geq\frac{\sqrt\kappa}{2h},
\qquad \kappa\geq64.
\end{equation}
Consequently,
\begin{equation}\label{lower:gradient-count}
\mathbb E N_\nabla\geq\sum_{i=1}^n\mathbb E J_i=n\mathbb E J.
\end{equation}
Combining \eqref{lower:EJ}, \eqref{lower:attempts}, and
$2\uppi\leq ph<2\uppi+1$ gives
\[
\mathbb E N_\nabla\geq
\frac1{10(2\uppi+1)}
\max\left\{p\sqrt\kappa,
       \frac{p\kappa\log\chi}{12\ell_\kappa}\right\}.
\]
The two terms account for the cost of one attempt and the additional
attempts required by rejection. Smaller steps increase the first cost;
the following choice of perturbation quantifies the rejection cost for
larger steps. Choose $\delta=(80p^2)^{-1}$. By \eqref{lower:expansion},
\eqref{lower:B}, and $\log(1+t)\geq15t/16$ for $0\leq t\leq1/15$,
\[
\log\chi\geq\frac34\delta b_0
\geq\frac{3\uppi}{2560p^3}.
\]
Also, $d\leq 6q_*^{-1}\log\kappa+3$ gives
$H_{\mathrm{cut}}+6\log\kappa=2d+18\log\kappa
\leq (13/q_*)(1+\log\kappa)$. Hence
\[
\ell_\kappa\leq\log(260/q_*)+\log(1+\log\kappa)
\leq8\Lambda_\kappa,
\]
where the last inequality uses $\log(260/q_*)<13<7\Lambda_\kappa$
for $\kappa\geq64$.
Substituting these bounds gives
\begin{equation}\label{lower:cost-refined}
 \mathbb E N_\nabla\ge
 \max\left\{c_3p\sqrt\kappa,\frac{c_4\kappa}{p^2\Lambda_\kappa}\right\},
\end{equation}
where
\[
 c_3=\frac{1}{10(2\uppi+1)},\qquad
 c_4=\frac{\uppi}{81920}\,c_3.
\]
Taking the weighted
geometric mean of the two terms, with weights $2/3$ and $1/3$, cancels $p$ and gives
\eqref{lower:cost}.
This completes the proof of Theorem~\ref{lower:thm}.

The initial law $\nu$ can be sampled explicitly. To sample a fast
coordinate, choose the
positive sign with probability
\[
\frac{\sqrt{\lambda_-}}{\sqrt{\lambda_+}+\sqrt{\lambda_-}},
\]
then set the coordinate to $|Z|/\sqrt{\lambda_+}$ for the positive
sign or $-|Z|/\sqrt{\lambda_-}$ for the negative sign, where $Z$ is
standard normal. Draw $w$ standard normal
and $z=\sqrt\kappa\,|Z_0|$ with $Z_0$ standard normal. These independent
draws give exactly \eqref{lower:initial}.

\appendix
\section{Smoothing and passage to the limit}\label{app:regularization}

We complete the regularization argument in the proof of
Lemma~\ref{ent:free-flight}.

\begin{lemma}[Regularity of the smoothed laws]\label{app:smoothed-regularity}
Assume \ref{ass:rhmc-a1} and \ref{ass:rhmc-a2}. Let $p$ be a compactly
supported probability measure on $\R^{2d}$ with $\KL(p\|\Pi)<\infty$.
For $(X,V)\sim p$ and an independent $\xi\sim\gamma_d$, set
\begin{gather*}
 p^\iota=\Law(X+\sqrt\iota\,\xi,V),\qquad 0<\iota\le1,\\
 p_t^\iota=\Law(Y_t,V),\qquad Y_t=X+tV+\sqrt\iota\,\xi.
\end{gather*}
Let $n_t$ be the density of $\nu_t^\iota=\Law(Y_t)$, and let
$u_t(y)=\mathbb E[V\mid Y_t=y]$ and
$\Sigma_t(y)=\operatorname{Cov}(V\mid Y_t=y)$.
On each finite time interval $[0,b]$, the following hold.
\begin{enumerate}[label=(\alph*),ref=(\alph*)]
\item\label{app:regularity-moments}
The fields $u_t$ and $\Sigma_t$ are smooth in $(t,y)$ and, together
with their first spatial derivatives, bounded on $[0,b]\times\R^d$.
\item\label{app:regularity-score}
For some $C>0$,
\[
 |\nabla\log n_t(y)|\le C(1+|y|),\qquad t\in[0,b],\ y\in\R^d.
\]
\item\label{app:regularity-gaussian}
For some $c,C>0$,
\[
 c\exp(-C(1+|y|^2))\le n_t(y)\le C\exp(-c|y|^2),
 \qquad t\in[0,b],\ y\in\R^d.
\]
The constants in \ref{app:regularity-moments},
\ref{app:regularity-score} and \ref{app:regularity-gaussian} may depend
on $\iota$.
\item\label{app:regularity-entropies}
All entropies appearing in \eqref{ent:current-bound} are finite for
the smoothed laws, and
\[
 \sup_{\substack{0<\iota\le1\\0\le t\le b}}
 \left\{\KL(p_t^\iota\|\Pi)+\KL(\nu_t^\iota\|\pi)
       +E_v(p_t^\iota)
       +\bigl|2E_v(p_t^\iota)-\Delta\mathcal L(\nu_t^\iota)\bigr|\right\}<\infty.
\]
\end{enumerate}
\end{lemma}

\begin{proof}
The conditional law of $(X,V)$ given $Y_t=y$ has
weight proportional to $\exp(-|y-X-tV|^2/(2\iota))$ relative to $p$.
Differentiation of the normalized expectation of any bounded scalar
function $F(X,V)$ gives
\[
 \nabla_y\mathbb E[F\mid Y_t=y]
 =\iota^{-1}\operatorname{Cov}(F,X+tV\mid Y_t=y).
\]
Apply this with $F=V_i$ and $F=V_iV_j$. Boundedness of $(X,V)$
gives the claimed bounds on $u_t,\Sigma_t$ and their spatial
derivatives. These moments are smooth in $(t,y)$ by differentiation
on compact sets. The score is
$\nabla\log n_t(y)=\iota^{-1}(\mathbb E[X+tV\mid Y_t=y]-y)$,
which has linear growth. Finally, $n_t$ is a mixture of translates
of a Gaussian with covariance $\iota I_d$ and bounded centers on
compact time intervals. It is therefore bounded below by
$c\exp(-C(1+|y|^2))$ and above by $C\exp(-c|y|^2)$ there, with
$c,C>0$ allowed to depend on $\iota$ and the interval.

Write $\mathrm h$ for differential entropy. Since the energy is
bounded on the support of $p$, finite joint relative entropy gives
finite $\mathrm h(p)$. The shear has unit Jacobian and spatial
convolution increases entropy by concavity, so
$\mathrm h(p_t^\iota)=\mathrm h(p^\iota)\ge\mathrm h(p)$.
Since $\Pi$ has density $e^{-H}/Z_\Pi$ with $Z_\Pi$ a normalizing
constant, for any law $q$ on phase space
\[
 \KL(q\Vert\Pi)=-\mathrm h(q)+\mathbb E_qH+\log Z_\Pi.
\]
Apply this with $q=p_t^\iota$. The entropy term is bounded above by
$-\mathrm h(p)$. The energy term is bounded because $H$ has at most
quadratic growth by \ref{ass:rhmc-a2} and the smoothed laws have
uniformly bounded second moments for $0<\iota\le1$ and $0\le t\le b$.
Hence $\KL(p_t^\iota\Vert\Pi)$ is bounded uniformly on this range. By
the chain rule $E_{\mathrm{joint}}=E_{\mathrm{pos}}+E_v$ with both
terms nonnegative, so $E_{\mathrm{pos}}$ and $E_v$ are bounded by the
joint entropy. Finally, by \eqref{pre:deltaL},
$2E_v-\Delta\mathcal L=2E_v-E_{\mathrm{pos}}+U_\Delta$, and
$|U_\Delta|$ is bounded by the same growth and moment bounds, so this
integrand is uniformly bounded on $[0,b]$.
\end{proof}

For each fixed $t$, the joint law $p_t^\iota$ and its position marginal
$\nu_t^\iota$ converge in $W_2$ to $p_t$ and $\nu_t$, respectively.
To pass from \eqref{ent:current-bound} for the smoothed law to
\eqref{ent:current-bound} for the original law, we let $\iota\to0$ and
use the following convergence claims.
\begin{enumerate}[label=(\roman*),ref=(\roman*)]
\item\label{app:convergence-entropies}
The joint and position differential entropies converge. Concavity gives the
lower bound under convolution, and upper semicontinuity under $W_2$
convergence gives the reverse limiting inequality. The latter
follows from lower semicontinuity of relative entropy to a fixed
Gaussian and convergence of second moments.
\item\label{app:convergence-energies}
The energy expectations
converge by quadratic growth and continuity, so the right side of
\eqref{ent:current-bound} converges by dominated convergence.
\item\label{app:convergence-currents}
The endpoint currents $\mathcal J(p_a^\iota)$ and $\mathcal J(p_b^\iota)$
converge to $\mathcal J(p_a)$ and $\mathcal J(p_b)$ by the continuity
argument of Step~1 in the proof of Lemma~\ref{ent:free-flight}, applied
to $\iota\to0$ in place of $s\to t$.
\end{enumerate}
This proves \eqref{ent:current-bound} for the original law.
Finally, $\mathcal W'=\mathcal J(p_t)$, the chain rule and
\eqref{pre:deltaL} give \eqref{ent:free-second}.

\section{A lower bound on gradient evaluations for Gaussian targets}\label{app:gaussian-lower}

\begin{proposition}\label{gaussian:cost-lower}
Let $0<m\le L$ and
$U(x_1,x_2)=\tfrac12mx_1^2+\tfrac12Lx_2^2$.
Consider the position kernel $P$ with Verlet step size $h$ satisfying
$h\sqrt L\le1$. The integration times $J$ may have any law with
$0<\E J<\infty$, drawn afresh at each transition independently of the
position, velocity and Metropolis uniform.
Set $A=\{x_1\le0\}$ and $\nu_0=2\mathbf1_A\pi$.
For every integer $n\ge0$,
\begin{equation}\label{gaussian:gradient-lower}
 \TV(\nu_0P^n,\pi)\le\tfrac14
 \quad\Longrightarrow\quad
 \E C_{\nabla}\ge n\E J\ge\frac{3\uppi}{16h\sqrt m}\ge\frac{3\uppi}{16}\sqrt\kappa,
\end{equation}
where $C_{\nabla}$ is the gradient count in \eqref{up:exact-count}.
\end{proposition}

\begin{proof}
Put $a=h\sqrt m$, $r=\sqrt{1-a^2/4}$ and
$\theta=2\arcsin(a/2)$. Let $(x_1,x_2)$ be the position and $(v_1,v_2)$
the velocity drawn at the start of a transition. The two coordinates
evolve independently under Verlet, and solving the scalar recurrence in
the first coordinate gives its proposed value after $j$ steps as
$Y_j/\sqrt m$, where
\[
 Y_j=\cos(j\theta)X+\frac{\sin(j\theta)}r V,
 \qquad X=\sqrt m\,x_1,\quad V=v_1.
\]
Under $\Pi$, $X,V$ are independent standard normals. Let
$\vartheta_j\in[0,\uppi]$ be the angle between $(1,0)$ and
$(\cos(j\theta),\sin(j\theta)/r)$. Rotational symmetry gives
$\Pbb(X\le0,Y_j>0)=\vartheta_j/(2\uppi)$.
Since acceptance can only remove these crossings and $\pi(A)=1/2$,
\[
 \frac{\int_A K_j(x,A^c)\,\pi(\mathrm dx)}{\pi(A)}
 \le\frac{\vartheta_j}{\uppi}.
\]
The continuous polar angle of $(\cos t,\sin t/r)$ has derivative
$r/(r^2\cos^2t+\sin^2t)\le1/r$, so
$\vartheta_j\le j\theta/r$. Moreover,
$\theta\le a/r$ and $r^2\ge3/4$ because $a\le h\sqrt L\le1$.
Averaging over $J$ therefore yields
\begin{equation}\label{gaussian:crossing}
 \frac{\int_A P(x,A^c)\,\pi(\mathrm dx)}{\pi(A)}
 \le\frac{4}{3\uppi}h\sqrt m\,\E J.
\end{equation}

Each $K_j$ preserves $\pi$ by Lemma~\ref{up:reversibility}, hence so
does $P$. Since $\nu_0\le2\pi$ as measures and $P$ is positive and
$\pi$-invariant, $\nu_0P^i\le2\pi P^i=2\pi$ for every $i\ge0$.
For the position chain $(X_i)$ started from $\nu_0$, this gives
\[
 \Pbb(X_i\in A,X_{i+1}\in A^c)
 \le2\int_A P(x,A^c)\,\pi(\mathrm dx)
 \le\frac{4}{3\uppi}h\sqrt m\,\E J.
\]
Since $X_0\in A$, reaching $A^c$ by time $n$ requires at least one
such crossing:
\[
 \{X_n\in A^c\}\subseteq\bigcup_{i=0}^{n-1}\{X_i\in A,\,X_{i+1}\in A^c\},
\]
so a union bound gives
$\Pbb(X_n\in A^c)\le\sum_{i<n}\Pbb(X_i\in A,X_{i+1}\in A^c)
\le\frac{4}{3\uppi}h\sqrt m\,n\E J$. Testing the total variation
distance on the set $A^c$, which has $\pi(A^c)=1/2$ and
$\nu_0P^n(A^c)=\Pbb(X_n\in A^c)$, gives
\[
 \TV(\nu_0P^n,\pi)
 \ge\pi(A^c)-\nu_0P^n(A^c)
 =\tfrac12-\Pbb(X_n\in A^c)
 \ge\tfrac12-\frac{4}{3\uppi}h\sqrt m\,n\E J.
\]
This proves \eqref{gaussian:gradient-lower}, since
$\E C_{\nabla}\ge\E\sum_{i=1}^nJ_i=n\E J$ by
\eqref{up:exact-count}. Finally,
$h\sqrt m=(h\sqrt L)/\sqrt\kappa\le1/\sqrt\kappa$ gives the last bound.
\end{proof}

\section{Fourier bounds for the three duration laws}\label{app:duration-fourier}

\begin{proof}[Proof of Lemma~\ref{up:example-fourier}]
It suffices to consider $0<\theta<\uppi$. Put
$u=\theta/2$ and $B_N=\sin(Nu)/(N\tan u)$. Summing the weights
\eqref{up:triangular} gives
\[
 \operatorname{Re}\widehat q(\theta)
 =\frac{\sin^2(Nu)}{N^2\sin^2u}
       -\frac{\sin^2(Nu)}{N^2}=B_N^2.
\]
The endpoint values follow by continuity. Thus for the triangular law
$w=q$, and $c=\beta_P=1$ are valid. The lower spectral bound is
attained at $\theta=\uppi$.

For uniform durations, integration over the rounding intervals gives
\[
 \operatorname{Re}\widehat\mu_J(\theta)=
 \begin{cases}B_N,&N\text{ even},\\ \cos u\,B_N,&N\text{ odd}.
 \end{cases}
\]
For even $N$, subtracting $(1-B_N^2)/2$ from $1-B_N$ gives
$(1-B_N)^2/2$. For odd $N$, the corresponding difference is
$[(B_N-\cos u)^2+\sin^2u]/2$. Hence $c=1/2$ is valid for every $N$.
For even $N$, the ratio tends to $1/2$ as $\theta\to0$, so the comparison constant is optimal in this case.
Also $B_N=(u/\tan u)\operatorname{sinc}(Nu)$, where
$\operatorname{sinc}x=\sin x/x$ for $x\ne0$ and $\operatorname{sinc}0=1$.
The inequality $\operatorname{sinc}x\ge-1/4$ holds for $x\ge0$:
it is immediate on $[0,\uppi]$; on $[\uppi,4]$ use
$\sin x\ge-(x-\uppi)$ and $\uppi>3$; on $[4,\infty)$ use $\sin x\ge-1$.
Both $u/\tan u$ and $\cos u$ lie in $[0,1]$, proving that $\beta_P=3/4$ is valid.

For exponential durations, put $r=e^{-4/N}$ and $z=e^{i\theta}$.
The full characteristic function and its real part are
\[
 \widehat\mu_J(\theta)=\frac{(1+z)(1-r)}{2(1-rz)},\qquad
 \operatorname{Re}\widehat\mu_J(\theta)
 =\frac{\tanh^2(2/N)}{\tanh^2(2/N)+\tan^2u}.
\]
The real part is nonnegative and vanishes at $\theta=\uppi$, so $\beta_P=1$.
We first prove the full-law comparison
\begin{equation}\label{up:exponential-full-comparison}
 1-\operatorname{Re}\widehat\mu_J(\theta)\ge\frac58(1-B_N^2).
\end{equation}
Put $v=N\tan u$, $x=Nu$, and $y=v^2$. Since $N\tanh(2/N)\le2$,
the left side is at least $y/(4+y)$. For $y\ge20/3$ this is at least
$5/8$, which suffices. For $0\le y\le7$, we have $x\le v<\uppi$ and
\[
 \frac{x}{v}=\frac{\arctan(v/N)}{v/N}\ge1-\frac{y}{768},\qquad
 \operatorname{sinc}x\ge\operatorname{sinc}v
 \ge f(y):=1-\frac y6+\frac{y^2}{120}-\frac{y^3}{5040}>0.
\]
The first inequality uses $N\ge16$, and the last is the alternating
Taylor bound for sine. Thus $B_N^2\ge Q(y):=(1-y/768)^2f(y)^2$.
For the polynomial
\[
 \mathcal P(y)=1-\frac58(4+y)\frac{1-Q(y)}y,
\]
with its continuous value at zero, $\mathcal P(y)\ge1/250$ on $[0,7]$.
Indeed, on each interval $[j,j+1]$, $j=0,\dots,6$, the Bernstein
coefficients of degree eight of $\mathcal P$ are rational numbers that
can be computed exactly, and all of them exceed $1/250$; the smallest,
on $[3,4]$, is about $0.0048$. Since a polynomial is bounded below on an
interval by the minimum of its Bernstein coefficients, this proves
\eqref{up:exponential-full-comparison}, with overlapping ranges
$[0,7]$ and $[20/3,\infty)$.

It remains to condition on $\tau\le8/\sqrt m$. Put $t=e^{-8}$.
By memorylessness and $8/(h\sqrt m)=2N$,
\[
 \widehat w(\theta)=\widehat\mu_J(\theta)
                       \frac{1-t e^{2Ni\theta}}{1-t}.
\]
Using $1-\cos(a+b)\le2(1-\cos a)+2(1-\cos b)$ gives
\[
 1-\operatorname{Re}\widehat w(\theta)
 \ge\frac{(1-2t)(1-\operatorname{Re}\widehat\mu_J(\theta))
                 -2t(1-\cos(2N\theta))}{1-t}.
\]
For $x=N\theta/2$, we have
\[
 1-B_N^2\ge\frac27\min\{x^2,1\},\qquad
 1-\cos(2N\theta)=1-\cos(4x)\le\min\{8x^2,2\}.
\]
For the first inequality on $[0,1]$, use
$\operatorname{sinc}x\le1-19x^2/120$, which gives
$1-\operatorname{sinc}^2x\ge4199x^2/14400>2x^2/7$.
For $x\ge1$, $|\operatorname{sinc}x|\le\sin1$ gives the same bound at one.
Consequently the ratio of the second expression to $1-B_N^2$ is at
most $28$. We conclude that
\[
 \frac{1-\operatorname{Re}\widehat w(\theta)}
      {1-\operatorname{Re}\widehat q(\theta)}
 \ge\frac{(1-2t)5/8-56t}{1-t}\ge\frac35.
\]
The final numerical inequality follows from $t<1/2980$, certified by
$e^8\ge\sum_{j=0}^{30}8^j/j!>2980$. Also $e^{-8}<3/16$, so this cutoff
meets the omitted-mass requirement for $\delta_\rho=3/4$.
\end{proof}

\bigskip
\noindent\textbf{Acknowledgements.}
I thank Andre Wibisono for comments on a draft of this paper, including
the observation in Remark~\ref{up:talagrand-remark}. The ideas in this
paper are my own. I used publicly available large language models to
assist with calculations, proofs of technical lemmas, numerical
experiments, figures, drafting, and the Lean formalization.

\bibliographystyle{spr-ims-number}
\bibliography{rhmc_upper_lower}

\end{document}